\documentclass{carticle}

\title{\textbf{Numerical solutions of Hamiltonian \textsc{pde}s:\\a multi-symplectic integrator in light-cone coordinates}}
\date{\today}
\author{Hugo \textsc{Ricateau}}
\author{Leticia F. \textsc{Cugliandolo}}
\affil{Sorbonne Universit\'es, Universit\'e Pierre et Marie Curie - Paris VI, Laboratoire de Physique Th\'eorique et Hautes \'Energies, 4 Place Jussieu, 75252 Paris Cedex 05, France}
\abstractcontent{%
	We introduce a novel numerical method to integrate partial differential equations representing the Hamiltonian dynamics of field theories.
	It is a multi-symplectic integrator that locally conserves the stress-energy tensor with an excellent precision over very long periods.
	Its major advantage is that it is extremely simple (it is basically a centered box scheme) while remaining locally well defined.
	We put it to the test in the case of the non-linear wave equation (with quartic potential) in one spatial dimension, and we explain how to implement it in higher dimensions.
	A formal geometric presentation of the multi-symplectic structure is also given as well as a technical trick allowing to solve the degeneracy problem that potentially accompanies the multi-symplectic structure.
}

\begin{document}
	\maketitlepage
	\begingroup
	\begin{wrapfigure}[10]{l}{10pt}
	\end{wrapfigure}
	\sectionstar{Introduction}
		Partial differential equations (\textsc{pde}s) are involved in the description of many scientific problems of interest.
		Quite commonly, understanding the behavior of the latter starts by solving a \textsc{pde}.
		Unfortunately, the ability to find analytical solutions to such equations is the exception rather than the rule.
		
		Usual strategies to tackle \textsc{pde}s are to consider an approximation of the initial equation (by removing some non-leading terms, for example) or a particular domain of the parameter space (perturbative approaches, \dots).
		However, these partial pieces of information can be insufficient to understand the behavior of the system in a satisfactory way.
		
		It then becomes relevant to focus on approximate solutions, but this time, of the original equation and for the full range of variation of the parameters.
		This is exactly what one tries to achieve by using numerical methods.
		The question then arises as to how to control the numerical approximation.
		
		\endgroup
		To be more precise, let us take a time dependent process $\rho\lr{(t)}$, governed by a differential equation
		\begin{equation*}
			f\lr{(\rho,\rho',\rho'',\cdots)}=0\eqpc
		\end{equation*}
		where the $'$ indicates time derivative.
		In a finite-difference representation of this equation, the approximation process is quite well controlled and, at each time step, we know the order of the error made.
		Yet, \emph{a priori}, we cannot predict the accumulation of these errors over many time steps and we cannot control the approximation made on the time-dependent solution (especially in the long-time limit).
		
		The question can then be rephrased as, why should we trust a solution obtained with a numerical solver?
		
		To address this question the standard procedure is to test, as precisely as possible, all the known properties of the problem.
		Firstly:
		\begin{enumerate}[label = \emph{\roman{*}}., labelindent = 0em, leftmargin = *, widest* = 2, nosep]
			\item If a particular solution is available, we can easily check whether the numerical solution is in agreement with it.
			\item In the same spirit, we can compare a numerical solution to the exact one for some particular choices of the parameters (by turning off all the interaction terms, for example).
		\end{enumerate}
		Nevertheless, these two kinds of tests are not robust enough and nothing ensures that the numerical approximation will behave in the same way in a different regime (where no exact solutions are available).
		
		The second kind of test is based on symmetries and conservation laws:
		\begin{enumerate}[label = \emph{\roman{*}}., labelindent = 0em, leftmargin = *, widest* = 2, nosep]
			\item If the theory admits a symmetry group (not spontaneously broken) we expect the numerical solutions to be (as closely as possible) in agreement with the discrete analogue of this symmetry group (and, especially, the discrete part of it).
			\item Due to the symmetry group, the theory can exhibit some conserved quantities that the numerical solutions should preserve as closely as possible.
		\end{enumerate}
		Using such tools lets us recover the ability to keep control over the numerical approximations and authorizes us to trust (or not) the numerical solutions obtained.
		
		Obviously, the accessibility to such validations is closely related to the structure of the \textsc{pde} and the method we are going to introduce will only be applicable to the Hamiltonian ones (\emph{i.e.} arising from \textsc{De Donder}~-- \textsc{Weyl}~-- \textsc{Hamilton} equations) that will be defined later.
		Let us simply state for the moment that they are Lagrangian \textsc{pde}s (\emph{i.e.} arising from an \textsc{Euler}~-- \textsc{Lagrange} equation).
		
		We also have to stress that we will focus on the particular class of finite-difference methods.
		Other classes (spectral methods, finite-element methods, \dots) are not as natural as finite-difference ones to fulfil the constraints involved by the preservation of the multi-symplectic structure to be explained below.
		
		A finite-difference method is like a cooking recipe composed of two ingredients:
		\begin{enumerate}[label = \emph{\roman{*}}., labelindent = 0em, leftmargin = *, widest* = 2, nosep]
			\item Firstly, a lattice that samples a bounded region of the support (\emph{e.g.} the space-time manifold).
			\item Secondly, a set of discretization rules that translate the continuous quantities to their lattice analogues.
			The continuous unknowns, defined on the space-time manifold, are sampled through the lattice.
			The discretization rules specify how to combine these samples if we want to compute derivatives, force terms, \dots
		\end{enumerate}
		
		Applying these rules to the equation of motion toggles from a \textsc{pde} to a set of algebraic equations (governing the behavior of the quantities defined on the lattice).
		Solving these algebraic equations leads to a set of values on the lattice nodes.
		This is a sampling of the solution, and in adjunction with some interpolation rules, an approximate solution of the \textsc{pde} is thus constructed.
		However, it has to be noted that these samples are not necessarily exact and, both the samples and the interpolation process are responsible for dissimilarities with the exact solution.
		
		Many standard finite-difference schemes already exist and are often adequate.
		Each method has its own preferred application field.
		In the kind of problems we will be interested in, we need to control the very long time behavior (with respect to a characteristic time-scale in the system) and we need a procedure that minimises the error accumulated over a huge number of steps.
		For this reason, we need to develop our own numerical scheme that performs well over long time scales (even though we may have to make some compromise on its short time quality).
		
		Generally, most methods are able to behave rather correctly on short time-scales.
		Therefore, the simplest and the fastest the method, the better it is in this regime.
		However, the problem complicates at long times, since two phenomena conspire against the performance of most strategies:
		\begin{enumerate}[label = \emph{\roman{*}}., labelindent = 0em, leftmargin = *, widest* = 2, nosep]
			\item On the one hand, \emph{a priori}, the reduction of the time step improves the quality of the approximation and then the quality of the solution but this obviously inhibits reaching long times.
			\item On the other hand, at each step, some numerical truncation errors are induced by the finite precision of the numbers' representation in a computer.
			Such errors are generally inflated when the size of the step decreases and accumulate as the number of steps increases.
		\end{enumerate}
		It results that, for a given final time, the precision is bounded from above, hence the necessity to chose a numerical method designed to behave correctly whatever the number of steps to handle.
		
		Regarding mechanical systems (support is of dimension one, \emph{e.g.} just time), it exists a very particular class of finite-difference integrators: the symplectic ones.
		They are well known (especially by researchers in planetary evolution) because of their very good capability to preserve the energy of Hamiltonian systems with a high accuracy even over long times~\cite{MclachlanQuispel2006}.
		Such integrators are based on the conservation of a very important (even central) structure of mechanical systems: the symplecticity of the phase space.
		
		Generalizations to field theories (\textsc{pde}) appends some difficulties since the conservation of the energy is no longer rigid enough.
		Actually, the correct fundamental quantity to be preserved is now the stress-energy tensor.
		Its conservation is local (by opposition to the conservation of the energy, which is through a space integral) and hence more fundamental.
		Therefore, the symplectic structure is no longer adapted and needs to be generalized.
		
		Multi-symplectic numerical integrators, introduced by \textsc{Bridges} and \textsc{Reich} at the beginning of the 21\textsuperscript{st} century~\cite{Bridges1997,BridgesReich2001,BridgesReich2006}, generalize to \textsc{pde}s the concept of symplectic integrators.
		Applied to conservative \textsc{pde}s, multi-symplectic integrators exhibit excellent local conservation properties (especially of the stress-energy tensor) and a very stable behavior for long time integrations~\cite{FrankMooreReich2006}.
		
		In the past fifteen years the subject has been widely studied~\cite{HongLiuSun2006,IslasSchober2004,IslasSchober2005,IslasSchober2006,LiuZhang2006,MaKongHongCao2011,MclachlanRylandSun2014,MooreReich2003} and successfully applied to a broad variety of problems including the non-linear \textsc{Schr\"odinger} equation~\cite{ChenQin2001,ChenQinTang2002,YamingHuajunSonghe2011,HongLiuLi2007,IslasKarpeevSchober2001,ZhuChenSongHu2011}, the non-linear \textsc{Dirac} equation~\cite{HongLi2006}, the \textsc{Maxwell} equations~\cite{SunTse2011}, the \textsc{Klein}~-- \textsc{Gordon} equation~\cite{Chen2006}, the \textsc{Korteweg}~-- \textsc{de Vries} (\textsc{KdV}) equation~\cite{AscherMclachlan2004,AscherMclachlan2005,WeipengZichen2008,YushunBinXin2007}, the \textsc{Boussinesq} equation~\cite{Chen2005}, as well as the \textsc{Zakharov}~-- \textsc{Kuznetsov} (ZK) equation~\cite{Chen2003}.
		
		We introduce here a new finite-difference multi-symplectic method based on the centered box scheme.
		The latter was one of the first multi-symplectic schemes introduced~\cite{BridgesReich2001} and it has been proved that it is stable and possesses a number of desirable properties~\cite{FrankMooreReich2006}.
		However, it is not well defined locally~\cite{RylandMclachlan2008,RylandMclachlanFrank2007}, so it requires a global solver and hence is not scalable.
		The idea we introduce in this paper is to use a rotated lattice in the light-cone coordinates, that restores locality of the algorithm.
		
		More precisely, the aim of this work is to introduce a multi-symplectic integrator that is apt to efficiently obtain the long-time dynamics of a field theory with high precision.
		We organize the presentation in a practical way, first showing how the method performs compared to other ones in the market, and next discussing its theoretical justification.
		
		We also discuss, without any assumption on the dimension of space-time, the problem of the degeneracy of the multi-symplectic structure and we show how to solve it in the particular case of the non-linear wave equation.
		
		The outline of the paper is the following.
		
		In the first part of this paper, \cref{msilccvsother}, we will compare our multi-symplectic integrator in the light-cone coordinates (\textsc{msilcc}) to two standard methods:
		\begin{enumerate}[label = \emph{\roman{*}}., labelindent = 0em, leftmargin = *, widest* = 2, nosep]
			\item On the one hand, a very basic scheme based on the \textsc{Newton} approximation of derivatives (this method is widely used by a broad community and proves to be preserving the multi-symplectic structure too, although it is explicit and thus of lesser quality).
			\item On the other hand, the method proposed by \textsc{Boyanovsky}, \textsc{Destri} and \textsc{de Vega}~\cite{BoyanovskyDestriVega2004} constructed such that it exactly conserves the energy of the system (non-local conservation).
		\end{enumerate}
		The comparison will be performed using the so-called $\lambda\,\phi^4$ field theory in dimension $D=1+1$.
		We will, in particular, study the local conservation (or not) of the stress-energy tensor.
		This example will allow us to emphasize the strengths and weaknesses of our method.
		
		The second part of this paper, \cref{preliminaries}, will be devoted to the definition of the necessary concepts, to exhibit the multi-symplectic structure, to deduce from it the local conservation laws (as well as the global ones) and, finally, to present how to rewrite the equations to prepare the implementation of the \textsc{msilcc} method.
		Throughout this section, the concepts and results will be illustrated through the example of the non-linear wave equation (whose $\lambda\,\phi^4$ theory, used in \cref{msilccvsother} of the paper, is a particular case).
		
		Finally, in \cref{implementation} we will introduce our method in detail.
		We will also demonstrate the conservation properties as well as a review of the concrete solving methods of the algebraic equations involved in the numerical approximation (at this point, the necessity to work in light-cone coordinates will clearly appear).
		Again, the non-linear wave equation will be our working example.
		
		A short conclusions section will close the paper.
	\section{The multi-symplectic integrator in light-cone coordinates (\texorpdfstring{\textsc{msilcc}}{MSILCC}) versus the standard methods}\label{msilccvsother}
		This first section will be devoted to the comparison of our multi-symplectic integrator in light-cone coordinates (\textsc{msilcc}) to two standard methods.
		The first one is a very basic scheme based on the \textsc{Newton} approximation of derivatives, directly implemented in the Lagrangian formulation of the \textsc{pde}.
		This method is the simplest and, generally, the fastest to implement, so it is widely used and it is an unavoidable starting point.
		The second method, developed by \textsc{Boyanovsky}, \textsc{Destri} and \textsc{de Vega}~\cite{BoyanovskyDestriVega2004}, is designed such that the total energy (a non-local quantity) will be exactly conserved whatever the configuration of the field or the size of the integration step.
		
		These two methods will be presented in detail throughout this section, while the \textsc{msilcc} method will be detailed in the next sections.
		Now, first of all, let us introduce the model which will support the comparison.
		\subsection{The \texorpdfstring{$\lambda\,\phi^4$}{lambda phi4} theory in \texorpdfstring{$1+1$}{1+1} dimensions}
			\subsubsection{The equation of motion}
				The comparison will be preformed on the so-called $\lambda\,\phi^4$ model\footnote{which belong in the class of the non-linear wave equation with the potential
				\begin{equation*}
					V\lr{(\phi)}=\frac{r}{2}\phi^2+\frac{\lambda}{4}\phi^4\eqpd
				\end{equation*}} in dimension $D=d+1=1+1$.
				The unknown is the real dynamic field, $\phi\lr{(x,t)}$, governed by a second order, non-linear, \textsc{pde} (the equation of motion):
				\begin{equation}
					\square\phi={\partial_0}^2\phi-{\partial_1}^2\phi=-V'\lr{(\phi)}=-\phi\lr{(1+\phi^2)}\eqpc\label{msilccvsother.model.pde}
				\end{equation}
				where $x$ and $t$ are space and time respectively, $\partial_0=\nicefrac{\displaystyle{\partial}}{\displaystyle{\partial{ct}}}$, $\partial_1=\nicefrac{\displaystyle{\partial}}{\displaystyle{\partial{x}}}$, $c$ is a characteristic speed (\emph{e.g.} the speed of light) that we set to one, $c=1$, and the derivative of the potential, $V$, is given by $V'\lr{(\phi)}=\nicefrac{\displaystyle{\partial V}}{\displaystyle{\partial\phi}}$.
				The $\lambda$ appearing in the name of the model is the parameter of the non-linear term and it has been set to one.
				The other parameter in the potential, the one that accompanies the quadratic term, has also been set to one in such a way that the potential has only one absolute minimum at $\phi=0$.
				
				There is no exact general solution to this equation.
				Nevertheless, some particular solutions can be obtained in terms of \textsc{Jacobi} elliptic functions \cite{NIST2010} and they can be useful as a first check of the accuracy of a numerical integrator (see \cref{msilccvsother.comparison.conclusion}).
			\subsubsection{Boundary and initial conditions}\label{msilccvsother.model.icbc}
				As previously mentioned a finite-difference method can be decomposed in terms of two ingredients: the lattice and the discretization rules.
				
				The notion of lattice is a bit ambiguous and needs to be clarified.
				First, let us suppose it to be a regular tiling (since there is, \emph{a priori}, no reason to take a more complex structure).
				Moreover, the spatial part of the lattice should be finite.
				Otherwise the integrator would have to solve an infinite number of algebraic equations (with the same amount of unknowns), which is generally impossible.
				
				Since the spatial part of the support is bounded, the solutions need to be constrained on the boundaries.
				In the following we will impose periodic boundary conditions (\textsc{pbc}) (even though this is not a requirement for our method) with a period of length $L$:
				\begin{equation}
					\begin{split}
						\phi\lr{(x+L,t)}&=\phi\lr{(x,t)}\eqpc\\
						\partial_0\phi\lr{(x+L,t)}&=\partial_0\phi\lr{(x,t)}\eqpd
					\end{split}\label{msilccvsother.model.bc}
				\end{equation}
				
				We use an initial condition that complies with the \textsc{pbc}:
				\begin{equation}
					\begin{split}
						\phi\lr{(x,0)}&=A\sin\lr{(\frac{2\pi\,x}{L})}\eqpc\\
						\partial_0\phi\lr{(x,0)}&=0\eqpd
					\end{split}\label{msilccvsother.model.ic}
				\end{equation}
				Therefore, the total energy is
				\begin{align*}
					\begin{split}
						E_\text{exact}&=\int_0^L\d{x}\,\lr{[\frac{1}{2}{\lr{(\partial_0\phi\lr{(x,0)})}}^2+\frac{1}{2}{\lr{(\partial_1\phi\lr{(x,0)})}}^2\vphantom{+V\lr{(\phi\lr{(x,0)})}}.}\\
						&\qquad\qquad\qquad\qquad\qquad\qquad\quad\lr{.\vphantom{\frac{1}{2}{\lr{(\partial_0\phi\lr{(x,0)})}}^2+\frac{1}{2}{\lr{(\partial_1\phi\lr{(x,0)})}}^2}+V\lr{(\phi\lr{(x,0)})}]}\eqpc
					\end{split}\\
					&=A^2\lr{[\frac{\pi^2}{L}+\frac{L\lr{(8+3A^2)}}{32}]}\eqpd
				\end{align*}
				
				The initial amplitude, $A$, allows us to control the predominance of the non-linearity.
				For a sufficiently small amplitude, the non-linear part of the potential will be dominated by its quadratic part and the initial condition leads to a time-dependent solution that is close to the second eigenmode of the linear wave equation: $A\sin{\lr{(\nicefrac{\displaystyle{2\pi\,x}}{\displaystyle{L}})}}\cos{\lr{(\nicefrac{\displaystyle{2\pi\,t}}{\displaystyle{L}})}}$.
				Conversely, when $A$ increases, the non-linear term becomes predominant and the behavior of the solution turns out to be much more complex.
				
\begin{figure}[!htb]
					\begin{center}
						\includegraphics{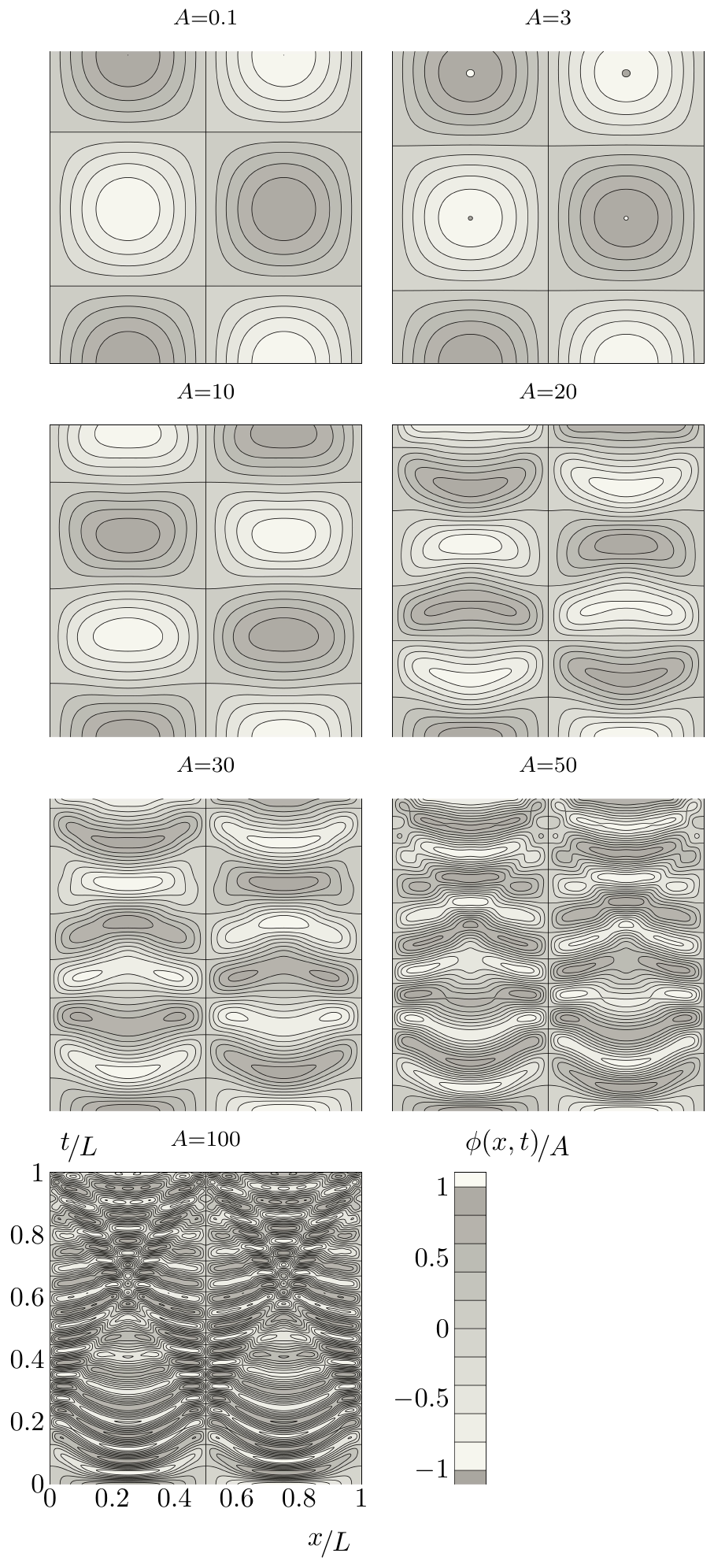}
					\end{center}
					\caption{Space-time plots of the solutions of \cref{msilccvsother.model.pde} with initial and boundary conditions given by \cref{msilccvsother.model.ic} and \cref{msilccvsother.model.bc} respectively.
					Different panels show data for different values of $A$, obtained with the \textsc{msilcc} method and $\nicefrac{L}{\sqrt{2}\,\delta}=1024$.
					Lines are iso-levels of the field while color is constant in between.
					\Cref{msilccvsother.model.crosssections} represents a cross-section of these space-time plots for the smallest values of $A$.}
					\label{msilccvsother.model.sols}
				\end{figure}

				\Cref{msilccvsother.model.sols} represents the short time behavior of the solution obtained using the \textsc{msilcc} method for different values of $A$.
				As expected, for $A=0.1$ the solution remains very close to the second eigenmode of the linear wave equation.
				For $A=3$, the solution evolves in two ways: its characteristic time-scale decreases, and the amplitude of the oscillations becomes a little bit bigger than $A$ (see in \cref{msilccvsother.model.sols}, $A=3$, the small circles at the center of the antinodes where the value of the field exceeds $A$).
				The impact of the non-linearity becomes significant.
				Then the larger is $A$, the shorter the characteristic time.
				The non-linearity is also destructing the structure of the eigenmode: when $A$ increases the solution is more and more distorted.
				
				The effect of $A$ is twofold, it will allow us to explore the influence of the non-linearity and the effect of decreasing the quality of the sampling when the typical variation scale of the field becomes closer and closer to the lattice spacing.
				
				\begin{figure}[!htb]
					\begin{center}
						\includegraphics{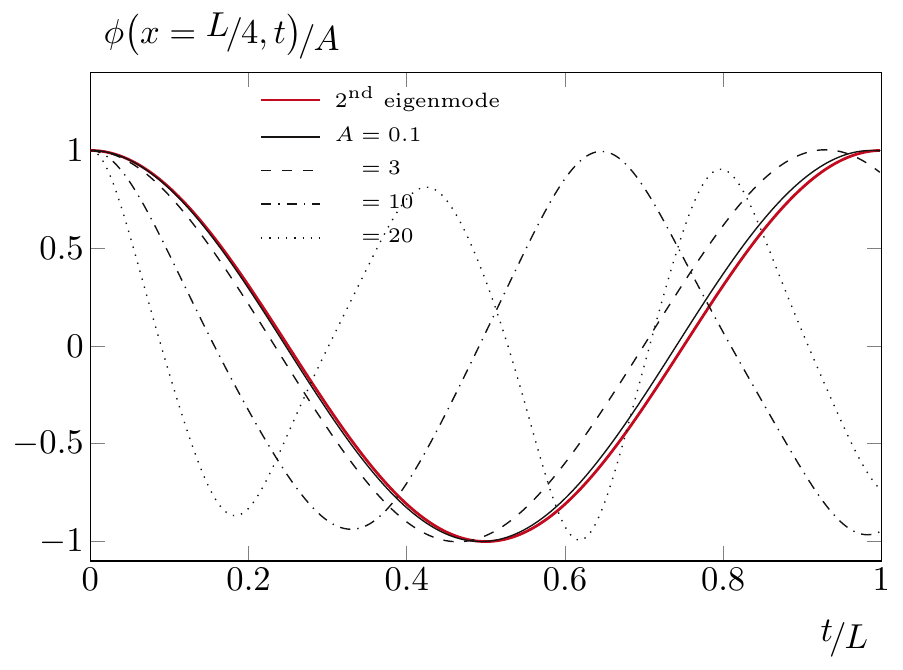}
					\end{center}
					\caption{Cross-sections of the space-time plots in \cref{msilccvsother.model.sols} along the axis $x=\nicefrac{\displaystyle{L}}{\displaystyle{4}}$.
					Red line is the cross-section of the second eigenmode of the linear theory: $A\cos\left(\nicefrac{\displaystyle{2\pi\,t}}{\displaystyle{L}}\right)$.}
					\label{msilccvsother.model.crosssections}
				\end{figure}
			\subsubsection{The stress-energy tensor, its conservation and the charges}
				As previously mentioned, the most fundamental quantity that the theory shall preserve is the stress-energy tensor (see \cref{preliminaries.dwh.stresstensorsection} for its definition, the one of the charges as well as the proof of their conservation).
				For the $\lambda\,\phi^4$ theory in $1+1$ dimensions the stress-energy tensor is symmetric and its four components read
				\begin{align*}
					\mathcal{T}^{00}&=\frac{1}{2}{\lr{(\partial_0\phi)}}^2+\frac{1}{2}{\lr{(\partial_1\phi)}}^2+\frac{1}{2}\phi^2+\frac{1}{4}\phi^4\eqpc\\
					\mathcal{T}^{01}&=\mathcal{T}^{10}=-\,\partial_0\phi\,\partial_1\phi\eqpc\\
					\mathcal{T}^{11}&=\frac{1}{2}{\lr{(\partial_0\phi)}}^2+\frac{1}{2}{\lr{(\partial_1\phi)}}^2-\frac{1}{2}\phi^2-\frac{1}{4}\phi^4\eqpd
				\end{align*}
				Its local conservation is given by
				\begin{subequations}
					\begin{align}
						\partial_0\mathcal{T}^{00}+\partial_1\mathcal{T}^{10}&=0\eqpc\label{msilccvsother.model.stresstensorconserv0}\\
						\partial_0\mathcal{T}^{01}+\partial_1\mathcal{T}^{11}&=0\eqpc\label{msilccvsother.model.stresstensorconserv1}
					\end{align}
				\end{subequations}
				or, in other words,
				\begin{align*}
					\partial_0\phi\,\lr{[\square\phi+\phi\lr{(1+\phi^2)}]}&=0\eqpc\\
					\partial_1\phi\,\lr{[\square\phi+\phi\lr{(1+\phi^2)}]}&=0\eqpc
				\end{align*}
				which are satisfied as long as the equation of motion (\ref{msilccvsother.model.pde}) holds.
				
				Conversely, the numerical equivalents of these local conservation laws (\cref{msilccvsother.model.stresstensorconserv0,msilccvsother.model.stresstensorconserv1}) will not be exactly satisfied.
				The violation comes from the fact that, in the discrete version of these equations, the term in brackets is not the discrete analogue of the equation of motion.
				This is precisely due to the fact that the discretization rules not always fulfil all the rules of differential calculus (\textsc{Leibniz}, \dots).
				
				Since these two quantities are non-vanishing they will allow us to control the quality of the numerical solution: a good numerical approximation should preserve, as closely as possible, the local conservation of the stress-energy tensor.
				These residues will be our first quantities of interest.
				
				Let us now define the charges as
				\begin{equation*}
					\mathcal{Q}^\nu=\int_0^L\mathcal{T}^{0\nu}\d{x}\eqpc
				\end{equation*}
				where $\nu$ is either $0$ or $1$.
				These are global quantities.
				Integrating over space, the local conservation of the stress-energy tensor leads to the conservation of the charges,
				\begin{equation*}
					\partial_0\mathcal{Q}^\nu=0\eqpd
				\end{equation*}
				Again, these quantities are not exactly conserved numerically and the resulting residues will be our second quantity of interest.
			\subsubsection{The testing conditions}\label{msilccvsother.preliminaries.testconditions}
				In the previous sub-section we highlighted the quantities allowing us to examine the quality of a numerical approximation (of a Hamiltonian \textsc{pde}).
				Let us now introduce in which context they will be observed.
				
				For each numerical method we will examine, through two situations, the error committed on the conservation of the stress-energy tensor (local) as well as the error committed on the evaluation of the charges (global).
				Firstly, we will have a look at how these errors behave as functions of $A$.
				We recall that $A$ has a twofold effect: it affects the weight of the non-linearity, but also the quality of the sampling (since when $A$ increases, the characteristic time-scale decreases while the time-step remains fixed).
				So, we expect the errors to be lower at small $A$ than at large $A$.
				Secondly, we will fix $A=10$ and observe how the errors behave as a function of time.
				We will explore both short and long time behavior.
				
				Let us now define the symbolic operator $\Delta$.
				When applied on a continuous equation, its effect is to extract the residue from the discrete analogue of the equation, and divide this residue by a characteristic quantity such that the result is not dimensional and the error can thus be compared to $1$.
				
				Before entering into the evaluation of the quality of the numerical methods let us present their construction in detail.
		\subsection{The \texorpdfstring{\textsc{Newton}}{Newton} method}
			The \textsc{Newton} method is probably the simplest finite-difference method one can develop.
			Its ease of use and its efficiency make it a classic.
			However, we will see that it can be inaccurate and even unstable.
			For now, let us describe its construction.
			\subsubsection{Sampling the space-time manifold}
				\begin{figure}[!htb]
					\begin{center}
						\includegraphics{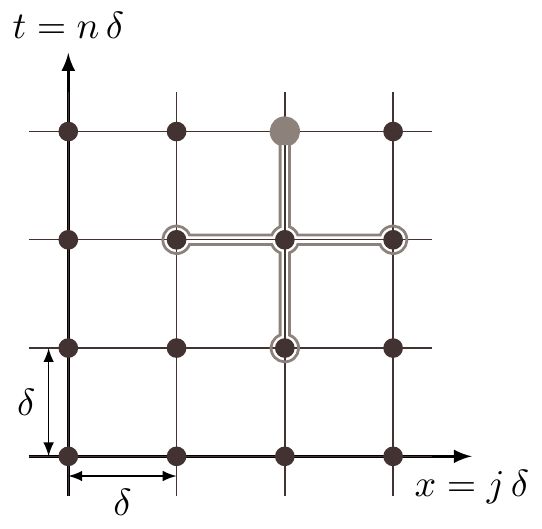}
					\end{center}
					\caption{Lattice description of the space-time manifold in the \textsc{Newton} method.
					The lattice spacing in space and time are chosen to be equal to preserve covariance.}
					\label{msilccvsother.newton.lattice}
				\end{figure}
				
				The support of the theory, $\mathcal{M}$, is a flat $1+1$ dimensional \textsc{Minkowski} space-time manifold.
				Taking into account the boundary conditions and the fact that the method will be used as an integrator, $\mathcal{M}$ becomes a flat half cylinder:
				\begin{equation*}
					\mathcal{M}=T^1\times\mathbb{R}_+\eqpc
				\end{equation*}
				where $T^1=S^1=\nicefrac{\displaystyle{\mathbb{R}}}{\displaystyle{L\,\mathbb{Z}}}$ is the flat one-dimensional torus of length $L$, and where, without loss of generality, the initial condition is supposed to be given at $t=0$.
				
				The lattice, $M$, will then be taken as a regular tiling of $\mathcal{M}$ with, as generator, a square of width $\delta$ aligned with the space and time coordinates.
				Therefore, the lattice is defined by
				\begin{equation*}
					M=\nicefrac{\displaystyle{\delta\,\mathbb{Z}}}{\displaystyle{N\,\mathbb{Z}}}\times\delta\,\mathbb{N}=\delta\,\mathbb{Z}_N\times\delta\,\mathbb{N}\eqpc
				\end{equation*}
				where $N\,\delta=L$.
				The geometry of $M$ is represented in \cref{msilccvsother.newton.lattice} and is nothing else than a square lattice.
				
				Now, the field, $\phi:\mathcal{M}\to\mathbb{R}$, can be sampled through the lattice as $\varphi:M\to\mathbb{R}$ such that
				\begin{equation*}
					\varphi_n^j=\phi\lr{(x=j\,\delta,t=n\,\delta)}\eqpc
				\end{equation*}
				where $n\in\mathbb{N}$ and $j\in\lr{\llbracket 0,N\llbracket}$.
				
				By this sampling process, at a given time, we have switched from the infinite number of degrees of freedom of the dynamic field to a representation with only a finite number ($N$) of degrees of freedom that can be used in a computer.
				This achieves the first step of the construction of the finite-difference approximation.
			\subsubsection{The \texorpdfstring{\textsc{Newton}}{Newton} scheme}
				In order to complete the construction of the finite-difference scheme, the second step is to provide the rules that will indicate how to combine the samples of the field (\emph{i.e.} the elements of $\varphi$) in order to obtain the physical quantities (and especially the equation of motion).
				
				The derivatives of the field will be approximated using the \textsc{Newton}'s rule which reads
				\begin{align*}
					\partial_0\phi\lr{(x=j\,\delta,t=n\,\delta)}&\approx {D_0^{_\pm}\varphi\,}_n^j\equiv \pm\frac{\varphi_{n\pm1}^j-\varphi_n^j}{\delta}\eqpc\\
					\partial_1\phi\lr{(x=j\,\delta,t=n\,\delta)}&\approx {D_1^{_\pm}\varphi\,}_n^j\equiv \pm\frac{\varphi_n^{j\pm1}-\varphi_n^j}{\delta}\eqpc
				\end{align*}
				where the $+$ (respectively $-$) stands for the forward (respectively backward) approximation.
				These two definitions ($+$ or $-$) are inequivalent and they express the fact that a finite difference can either represent the derivative at the end or at the beginning of the interval.
				This is exactly why \textsc{Newton}'s method is so easy to implement but is also the reason for its inaccuracy (a finite difference should only represent the derivative at the center of the interval).
				However, the centered \textsc{Newton}'s rule for the first order derivative (that combines $D^{_+}$ and $D^{_-}$ to involve the points $\#+1$ and $\#-1$) will not be used since it leads to an inconsistent approximation of the second order derivative\footnote{To show this fact, let us define the centered \textsc{Newton}'s rule as
					\begin{equation*}
						D^{_C}\varphi_n=\frac{D^{_+}\varphi_n+D^{_-}\varphi_n}{2}=\frac{\varphi_{n+1}-\varphi_{n-1}}{2\,\delta}\eqpd
					\end{equation*}
					The second order derivative, that reads
					\begin{equation*}
						D^{_C}D^{_C}\varphi_n=\frac{\varphi_{n+2}-2\varphi_n+\varphi_{n-2}}{4\,\delta^2}\eqpc
					\end{equation*}
					thus leads to two independent sub-lattices (the odd one and the even one).
				}.
				
				The equation of motion can be approximated in two ways: using forward then backward rules or vice versa, \emph{i.e.} $\partial^2\approx D^{_-}D^{_+}$ or $\partial^2\approx D^{_+}D^{_-}$ (either forward~-- forward or backward~-- backward leads to an inconsistent approximation of the second derivatives\footnote{since
					\begin{equation*}
						D^{_+}D^{_+}\varphi_n=\frac{\varphi_{n+2}-2\varphi_{n+1}+\varphi_n}{\delta^2}\eqpc
					\end{equation*}
					involves two unknowns ($\varphi_{n+2}$ and $\varphi_{n+1}$) and since
					\begin{equation*}
						D^{_-}D^{_-}\varphi_n=\frac{\varphi_n-2\varphi_{n-1}+\varphi_{n-2}}{\delta^2}\eqpc
					\end{equation*}
					requires to have solved the neighbouring equation in space to get $\varphi_n$ (which is incompatible with the periodic boundary conditions).
				}).
				In both cases the discrete version of the equation of motion at time $t=n\,\delta$ and position $x=j\,\delta$ reads
				\begin{equation*}
					\begin{split}
						\frac{\varphi_{n+1}^j-2\varphi_n^j+\varphi_{n-1}^j}{\delta^2}-\frac{\varphi_n^{j+1}-2\varphi_n^j+\varphi_n^{j-1}}{\delta^2}&=\\
						-\varphi_n^j&\lr{(1+{\varphi_n^j}^2)}\eqpd
					\end{split}
				\end{equation*}
				This algebraic equation is explicit in $\varphi_{n+1}^j$.
				The evolution of a given state can then be efficiently obtained using
				\begin{equation*}
					\varphi_{n+1}^j=\varphi_{n-1}^j+\varphi_n^{j-1}+\varphi_n^{j+1}-\delta^2\,\varphi_n^j\lr{(1+{\varphi_n^j}^2)}\eqpd
				\end{equation*}
				The nodes of the lattice that appear in this equation are highlighted in \cref{msilccvsother.newton.lattice}: the vertex in the left hand side of the equation is represented as
				\begin{tikzpicture}
					\fill[EE178] (0,0) circle(0.14);
				\end{tikzpicture}
				while the vertices involved in the right hand side of the equation are represented as
				\begin{tikzpicture}
					\fill[EE181] (0,0) circle(0.1);
					\draw[thick, EE178] (0,0) circle(0.14);
				\end{tikzpicture}
				.
				
				This concludes the definition of the \textsc{Newton} method.
				Let us now use these rules to obtain the discrete formulation of the stress-energy tensor.
			\subsubsection{The energy and the stress-energy tensor}
				As the derivatives can be approximated in two ways ($D^{_+}$ or $D^{_-}$), the stress-energy tensor can be defined in two ways too:
				\begin{align*}
					T_{^\pm}^{00}&=\frac{1}{2}{\lr{(D_0^{_\pm}\varphi)}}^2+\frac{1}{2}{\lr{(D_1^{_\pm}\varphi)}}^2+\frac{1}{2}{\varphi}^2+\frac{1}{4}{\varphi}^4\eqpc\\
					T_{^\pm}^{01}&=T_{^\pm}^{10}=-\,D_0^{_\pm}\varphi\,D_1^{_\pm}\varphi\eqpc\\
					T_{^\pm}^{11}&=\frac{1}{2}{\lr{(D_0^{_\pm}\varphi)}}^2+\frac{1}{2}{\lr{(D_1^{_\pm}\varphi)}}^2-\frac{1}{2}{\varphi}^2-\frac{1}{4}{\varphi}^4\eqpc
				\end{align*}
				where the space and time labels were omitted.
				These two definitions are inequivalent but both of them are valid and lead to a residue (again omitting the $n$ and $j$ indexes):
				\begin{align*}
					D_0^{_\mp}T_{^\pm}^{00}+D_1^{_\mp}T_{^\pm}^{10}&=\epsilon_{^\pm}^0\eqpc\\
					D_0^{_\mp}T_{^\pm}^{01}+D_1^{_\mp}T_{^\pm}^{11}&=\epsilon_{^\pm}^1\eqpd
				\end{align*}
				In practice the two definitions of these residues behave in the same way and the results will only present $\epsilon^0=\epsilon_{^+}^0$ and $\epsilon^1=\epsilon_{^+}^1$.
				
				Obviously, the same reasoning can be applied to the conservation of the charges but we do not detail it here.
			\subsubsection{Energy conservation}
\begin{figure}[!htb]
					\begin{center}
						\includegraphics{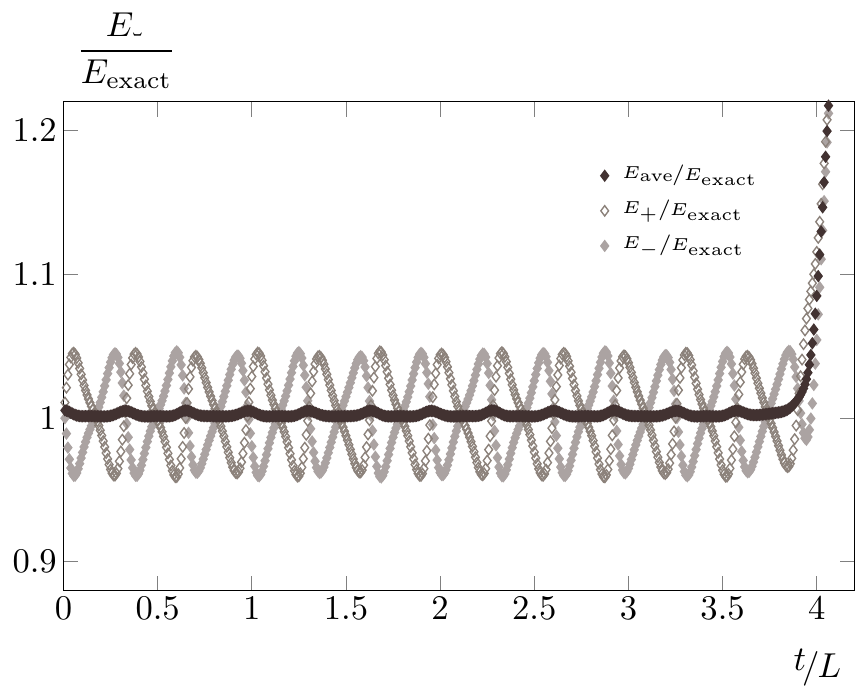}
					\end{center}
					\caption{Time evolution of the total energy for $A=10$ with the \textsc{Newton} integrator and $\nicefrac{L}{\delta}=128$.}
					\label{msilccvsother.newton.energy}
				\end{figure}

				The general treatment of the energy conservation will be exposed later, but let us briefly show how the total energy behaves.
				At time $t=n\,\delta$ (again, index $n$ will be omitted) it can be defined in two ways
				\begin{align*}
					E_{^+}&=Q_{^+}^0\eqpc\\
					E_{^-}&=Q_{^-}^0\eqpc
				\end{align*}
				where $Q^0_\pm$ is the charge defined as
				\begin{equation*}
					Q_{^\pm}^0=\delta\sum_{j=0}^{N-1}T_{^\pm}^{00}\eqpd
				\end{equation*}
				One can also envisage to combine these two definitions as
				\begin{equation*}
					E_{\text{ave}}=\frac{E_{^+}+E_{^-}}{2}\eqpd
				\end{equation*}
				These three definitions are represented as a function of time for $A=10$ in
				\cref{msilccvsother.newton.energy}.
				
				We observe that both $E_{^+}$ and $E_{^-}$ vary in time with an amplitude of the order of $5\%$ of their time-averaged value.
				The amplitude of the variations for $E_{\text{ave}}$ is reduced to $\sim1\%$ due to a compensation of the errors in $E_{^+}$ and $E_{^-}$.
				Nevertheless, $E_{\text{ave}}$ does not correspond to any discretization rule: it is the average of the energies obtained using different rules, which differs from the energy that would be obtained from the combination of the forward and backward rules (that, as already mentioned, leads to an incorrect approximation of the second order derivatives).
				
				We also notice that there is no apparent change in the amplitude of the deviations as time elapses up to a time-scale at which the energies rapidly diverge.
				The divergence of the energies is directly due to the divergence of the solution which has been destabilised by the integration method.
				
				We conclude that, although \textsc{Newton}'s method is straightforward to implement, it is inaccurate (the total energy conservation up to $1\%$ is not acceptable in most applications) and can even become unstable.
				Therefore, \textsc{Newton}'s method will not be a good choice to integrate a field theory over a long time.
				These observations will be confirmed by the study of the local conservation laws presented in the following sections.
				Before presenting the local analysis, let us first introduce another finite-difference method.
		\subsection{The \texorpdfstring{\textsc{Boyanovsky}~-- \textsc{Destri}~-- \textsc{de Vega}}{Boyanovsky - Destri - de Vega} (\texorpdfstring{\textsc{BDdV}}{BDdV}) method}
			The \textsc{Boyanovsky}~-- \textsc{Destri}~-- \textsc{de Vega} (\textsc{BDdV}) method~\cite{BoyanovskyDestriVega2004} that we present here has been developed such that it exactly preserves the total energy of the system, making it a good candidate for long time integrations.
			However, as previously said, the conservation of the total energy is not the most fundamental principle for a field theory that should foremost locally preserve the stress-energy tensor.
			\subsubsection{The lattice}\label{msilccvsother.bddv.latticesection}
				\begin{figure}[!htb]
					\begin{center}
						\includegraphics{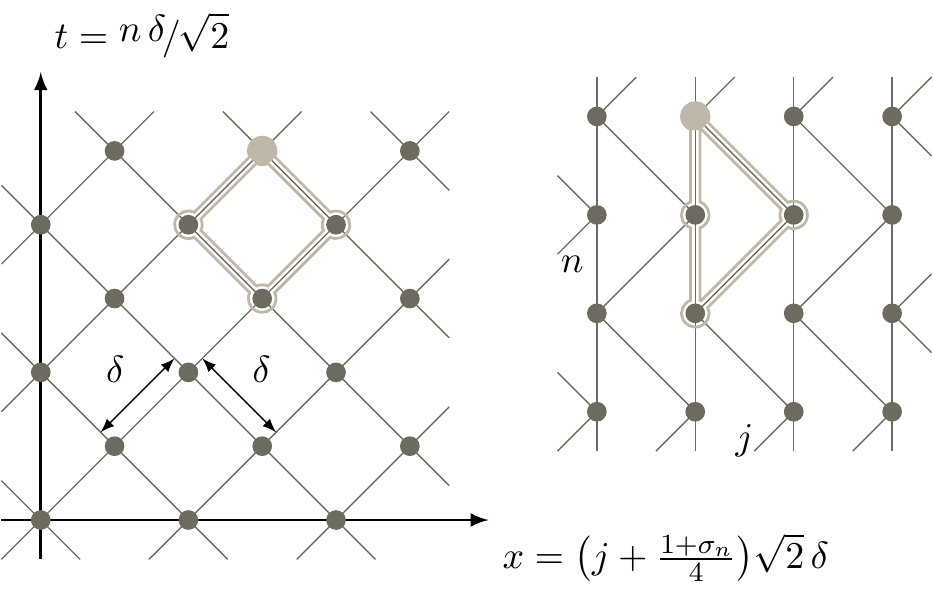}
					\end{center}
					\caption{Lattice description of the space-time manifold used in the \textsc{BDdV} method.}
					\label{msilccvsother.bddv.lattice}
				\end{figure}
				
				In the \textsc{BDdV} method the space-time manifold is rotated by $\nicefrac{\displaystyle{\pi}}{\displaystyle{4}}$.
				More precisely, the space-time manifold, $\mathcal{M}$, is unchanged, and the lattice, $M$, is still taken as a regular tilling of $\mathcal{M}$.
				The generator is still a square of width $\delta$, but aligned with the light-cone coordinates.
				Therefore, the lattice is defined by
				\begin{equation*}
					M=\lr{\{\lr{(\sqrt{2}\,\delta\,\lr{[j+\frac{1+\sigma_n}{4}]},\frac{n\,\delta}{\sqrt{2}})}\vphantom{n\in\mathbb{N}\,,j\in\mathbb{Z}_N}|}\lr{.\vphantom{\lr{(\lr{(j+\frac{1+\sigma_n}{4})}\sqrt{2}\,\delta\,,
					\frac{n\,\delta}{\sqrt{2}})}}n\in\mathbb{N}\,, \; j\in\mathbb{Z}_N\}}\eqpc
				\end{equation*}
				where
				\begin{equation*}
					\sigma_n=2\lr{(n\bmod 2)}-1\equiv\pm1\eqpd
				\end{equation*}
				$M$ is represented in \cref{msilccvsother.bddv.lattice} and is nothing else than a square lattice in the light-cone coordinate system which correctly respects the boundary conditions.
				
				Finally, in the same way as for the \textsc{Newton} method, the field, $\phi:\mathcal{M}\to\mathbb{R}$, can be sampled through the lattice as $\varphi:M\to\mathbb{R}$ such that
				\begin{equation*}
					\varphi_n^j=\phi\lr{(x=\sqrt{2}\,\delta\,\lr{[j+\frac{1+\sigma_n}{4}]},\,t=\frac{n\,\delta}{\sqrt{2}})}\eqpc
				\end{equation*}
				where $n\in\mathbb{N}$ and $j\in\lr{\llbracket 0,N\llbracket}$.
			\subsubsection{Exact energy preserving approximation}
				We now provide the rules that allow one to express the discrete analogues of the physical quantities, respecting the directions imposed by the lattice.
				Under this constraint, the derivatives are written along the light-cone coordinates as
				\begin{align*}
					\begin{split}
						\frac{\lr{(\partial_0-\partial_1)}\phi}{\sqrt{2}}&\lr{(x=\sqrt{2}\,\delta\,\lr{[j+\frac{1+3\,\sigma_n}{4}]},\,t=\frac{n\,\delta}{\sqrt{2}})}\\
						&\phantom{\lr{(\cdots)}}\approx {\check{D}_0^{_\pm}\varphi\,}_n^{j+\frac{\sigma_n}{2}}=\pm\frac{\varphi_{n\pm1}^j-\varphi_n^{j+\frac{\sigma_n\pm1}{2}}}{\delta}\eqpc
					\end{split}\\
					\frac{\lr{(\partial_0+\partial_1)}\phi}{\sqrt{2}}&\lr{(\cdots)}\approx {\check{D}_1^{_\pm}\varphi\,}_n^{j+\frac{\sigma_n}{2}}=\pm\frac{\varphi_{n\pm1}^j-\varphi_n^{j+\frac{\sigma_n\mp1}{2}}}{\delta}\eqpc
				\end{align*}
				where the $\cdots$ indicate that the field is evaluated at the same point on the space-time manifold as in the first equation.
				
				The method now differs from \textsc{Newton}'s since the discretization rules are not applied to the equation of motion but to the energy and the constraints imposed by its conservation are used to derive a modified discrete evolution equation.
				In the continuum limit this equation would be identical to the equation of motion, but in the discrete formulation it differs from the one we would have obtained had we directly applied the rules to the equation of motion.
				As a matter of fact, the idea introduced here is very deep but we will have the opportunity to come back to this later.
				
				The local energy density (the $00$ component of the stress-energy tensor) will then be approximated in two ways (following the same principle as for the two possible \textsc{Newton} approximations of the derivatives) and is given by
				\begin{equation}
					\begin{split}
						{T^{00}_{^\pm}\,}_n^{j+\frac{\sigma_n}{2}}&=\frac{1}{2}{\lr{({\check{D}_0^{_\pm}\varphi\,}_n^{j+\frac{\sigma_n}{2}})}}^2+\frac{1}{2}{\lr{({\check{D}_1^{_\pm}\varphi\,}_n^{j+\frac{\sigma_n}{2}})}}^2-\frac{1}{4}\\
						&\quad+\frac{1}{8}\lr{(1+{\varphi_{n\pm 1}^j}^2)}\lr{(2+{\varphi_n^j}^2+{\varphi_n^{j+\sigma_n}}^2)}\label{msilccvsother.bddv.setens.00component}\eqpc
					\end{split}
				\end{equation}
				where the $n$ and $j$ indices can no longer be omitted since they are not obvious.
				The difference between these two possible definitions reads
				\begin{equation*}
					{T^{00}_{^+}\,}_n^{j+\frac{\sigma_n}{2}}-{T^{00}_{^-}\,}_n^{j+\frac{\sigma_n}{2}}=\frac{\varphi_{n+1}^j-\varphi_{n-1}^j}{\delta^2}\,R_n^{j+\frac{\sigma_n}{2}}\eqpc
				\end{equation*}
				where
				\begin{equation*}
					\begin{split}
						R_n^{j+\frac{\sigma_n}{2}}&=\lr{(\varphi_{n+1}^j+\varphi_{n-1}^j)}\lr{[1+\frac{\delta^2}{8}\lr{(2+{\varphi_n^j}^2+{\varphi_n^{j+\sigma_n}}^2)}]}\\
						&\quad-\varphi_n^j-\varphi_n^{j+\sigma_n}\eqpd
					\end{split}
				\end{equation*}
				On the other hand, the total energy is given by
				\begin{equation*}
					{Q_{^\pm}^0}_n=\delta\sum_{j=0}^{N-1}{T_{^\pm}^{00}\,}_n^{j+\frac{\sigma_n}{2}}\eqpd
				\end{equation*}
				It can be shown (using periodic boundary conditions) that these two definitions are equivalent,
				\begin{equation*}
					{Q_{^+}^0}_n={Q_{^-}^0}_{n+1}=E_n\eqpc
				\end{equation*}
				defining the total energy at time $t=\nicefrac{\displaystyle{n\,\delta}}{\displaystyle{\sqrt{2}}}$ with no ambiguity.
				Now, this energy is exactly conserved if
				\begin{equation*}
					{Q_{^+}^0}_n=E_{n}=E_{n-1}={Q_{^-}^0}_n\eqpc
				\end{equation*}
				that will be satisfied as soon as
				\begin{equation*}
					{T^{00}_{^+}\,}_n^{j+\frac{\sigma_n}{2}}={T^{00}_{^-}\,}_n^{j+\frac{\sigma_n}{2}}\eqpc
				\end{equation*}
				that is to say, if
				\begin{equation*}
					R_n^{j+\frac{\sigma_n}{2}}=0\eqpd
				\end{equation*}
				Since $R$ involves the samples of the field at different times this equation is a sort of equation of motion.
				Moreover, since it is explicit in $\varphi_{n+1}^j$, the evolution of a given state can be efficiently followed using
				\begin{equation*}
					\varphi_{n+1}^j=-\varphi_{n-1}^j+\frac{\varphi_n^j+\varphi_n^{j+\sigma_n}}{\displaystyle{1+\frac{\delta^2}{8}\lr{(2+{\varphi_n^j}^2+{\varphi_n^{j+\sigma_n}}^2)}}}\eqpd
				\end{equation*}
				The nodes of the lattice that are involved in this equation are highlighted in \cref{msilccvsother.bddv.lattice}: the vertex in the left hand side of the equation is represented as
				\begin{tikzpicture}
					\fill[EE191] (0,0) circle(0.14);
				\end{tikzpicture}
				and the ones in the right hand side are represented as
				\begin{tikzpicture}
					\fill[EE194] (0,0) circle(0.1);
					\draw[thick, EE191] (0,0) circle(0.14);
				\end{tikzpicture}
				.
			\subsubsection{The stress-energy tensor}
				The $00$ component of the stress-energy tensor was defined in \cref{msilccvsother.bddv.setens.00component}.
				The two remaining independent components are
				\begin{align*}
					{T^{01}_{^\pm}\,}_n^{j+\frac{\sigma_n}{2}}&={T^{10}_{^\pm}\,}_\cdot^\cdot=\frac{1}{2}{\lr{({\check{D}_0^{_\pm}\varphi\,}_n^{j+\frac{\sigma_n}{2}})}}^2-\frac{1}{2}{\lr{({\check{D}_1^{_\pm}\varphi\,}_n^{j+\frac{\sigma_n}{2}})}}^2\eqpc\\
					\begin{split}
						{T^{11}_{^\pm}\,}_n^{j+\frac{\sigma_n}{2}}&=\frac{1}{2}{\lr{({\check{D}_0^{_\pm}\varphi\,}_n^{j+\frac{\sigma_n}{2}})}}^2+\frac{1}{2}{\lr{({\check{D}_1^{_\pm}\varphi\,}_n^{j+\frac{\sigma_n}{2}})}}^2+\frac{1}{4}\\
						&\quad-\frac{1}{8}\lr{(1+{\varphi_{n\pm 1}^j}^2)}\lr{(2+{\varphi_n^j}^2+{\varphi_n^{j+\sigma_n}}^2)}\eqpd
					\end{split}
				\end{align*}
				These two definitions ($+$ and $-$) are inequivalent (they only match once integrated over space) but both of them are valid and each one leads to two residues.
				\begin{align*}
					{\check{D}_0^{_\mp}\lr{(T_{^\pm}^{00}-T_{^\pm}^{10})}\,}_n^j+{\check{D}_1^{_\mp}\lr{(T_{^\pm}^{00}+T_{^\pm}^{10})}\,}_n^j&=\sqrt{2}\,{\epsilon_{^\pm}^0\,}_n^j\eqpc\\
					{\check{D}_0^{_\mp}\lr{(T_{^\pm}^{01}-T_{^\pm}^{11})}\,}_n^j+{\check{D}_1^{_\mp}\lr{(T_{^\pm}^{01}+T_{^\pm}^{11})}\,}_n^j&=\sqrt{2}\,{\epsilon_{^\pm}^1\,}_n^j\eqpc
				\end{align*}
				where
				\begin{align*}
					{\check{D}_0^{_\pm}T_{^\mp}^{\mu\nu}\,}_n^j&=\pm\frac{{T_{^\mp}^{\mu\nu}\,}_{n\pm1}^{j-\frac{\sigma_n}{2}}-{T_{^\mp}^{\mu\nu}\,}_n^{j\pm\frac{1}{2}}}{\delta}\eqpc\\
					{\check{D}_1^{_\pm}T_{^\mp}^{\mu\nu}\,}_n^j&=\pm\frac{{T_{^\mp}^{\mu\nu}\,}_{n\pm1}^{j-\frac{\sigma_n}{2}}-{T_{^\mp}^{\mu\nu}\,}_n^{j\mp\frac{1}{2}}}{\delta}\eqpc
				\end{align*}
				with both $\mu$ and $\nu$ being either $0$ or $1$.
				In practice, the two definitions of the residues behave in the same way and we will only present $\epsilon^0=\epsilon_{^+}^0$ and $\epsilon^1=\epsilon_{^+}^1$.
				
				The same reasoning can be applied to the conservation of the charges and will not be detailed.
			\subsubsection{Energy conservation}
\begin{figure}[!htb]
					\begin{center}
						\includegraphics{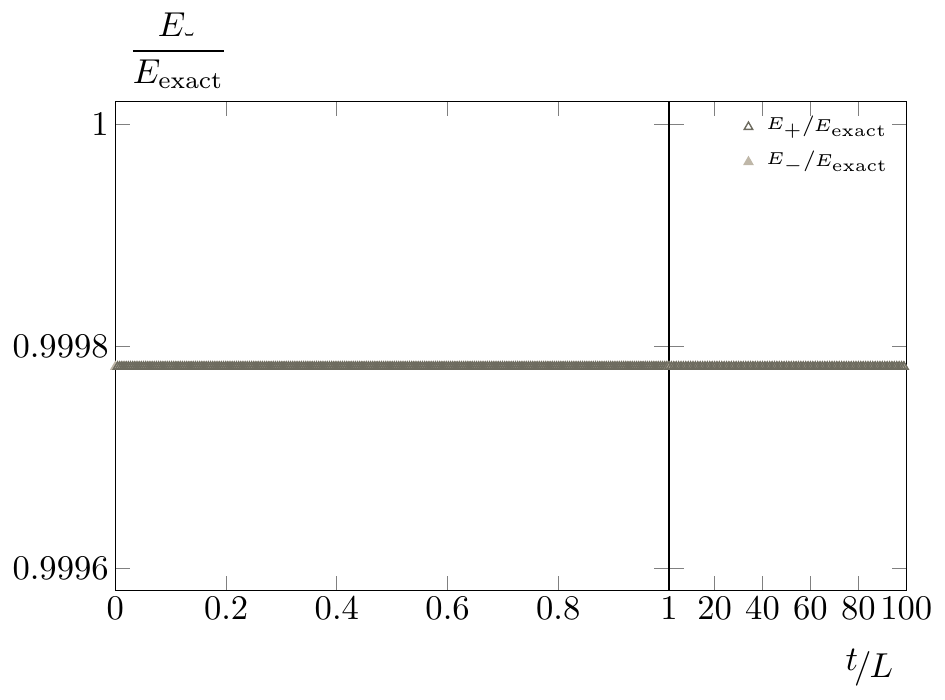}
					\end{center}
					\caption{Time evolution of the total energy for $A=10$ using the \textsc{BDdV} method with $\nicefrac{L}{\sqrt{2}\,\delta}=128$.
					Beyond $\nicefrac{\displaystyle{t}}{\displaystyle{L}}=1$ the horizontal axis is shown in a different linear scale and data at $255$ consecutive instants are skipped between two successive points.}
					\label{msilccvsother.bddv.energy}
				\end{figure}

				Although the energy is defined without ambiguity (since the two definitions of the stress-energy tensor are equivalent once integrated over space), we still define the total energy at time $t=\nicefrac{\displaystyle{n\,\delta}}{\displaystyle{\sqrt{2}}}$ in two ways
				\begin{equation*}
					{E_{^\pm}}_n={Q_{^\pm}^0}_n\eqpc
				\end{equation*}
				and we follow their time evolution independently.
				The numerical outcome is shown in \cref{msilccvsother.bddv.energy}.
				
				We observe, first of all, that the two definitions of the energy behave exactly in the same way, confirming that there is no ambiguity.
				Then we stress that the value of the energy differs from the exact one (the difference is of order $2\text{\textpertenthousand}$).
				This is not surprising and is due to the discretization process (what can be astonishing is that the difference is so small).
				We also observe that the total energy is exactly conserved, as expected.
				
				At this point, one could reasonably conclude that the \textsc{BDdV} method is a very good choice for the short and long time integration of conservative field theories.
				However, as we shall see in \cref{msilccvsother.comparison}, this conclusion would be premature.
				Unfortunately, the stress-energy tensor is not conserved locally as it is the total energy.
		\subsection{The \texorpdfstring{\textsc{msilcc}}{MSILCC} method: a short review of expected properties}
			The detailed description of the multi-symplectic integrator in light-cone coordinates (\textsc{msilcc}) will be given in the two last sections of the paper.
			
			We just want to stress here that the lattice is the same as the one used in the \textsc{BDdV} method (see \cref{msilccvsother.bddv.lattice}).
			In brief, the difference lies in the rules, which are closer to the ones employed in \textsc{Newton}'s method.
			
			The \textsc{msilcc} method is designed such that the discretization process exactly preserves the multi-symplectic structure of the phase space.
			It is also implemented in such a way to respect, as much as possible, the rules of differential calculus.
			The direct consequence is that the local conservation of the stress-energy tensor is remarkably good even on long time-scales.			
			However, the method is not engineered to conserve the global charges and we do not expect to have the same kind of ``magic'' compensation of local errors that ensures the \textsc{BDdV} method.
			
			\begin{figure}[!htb]
				\begin{center}
					\includegraphics{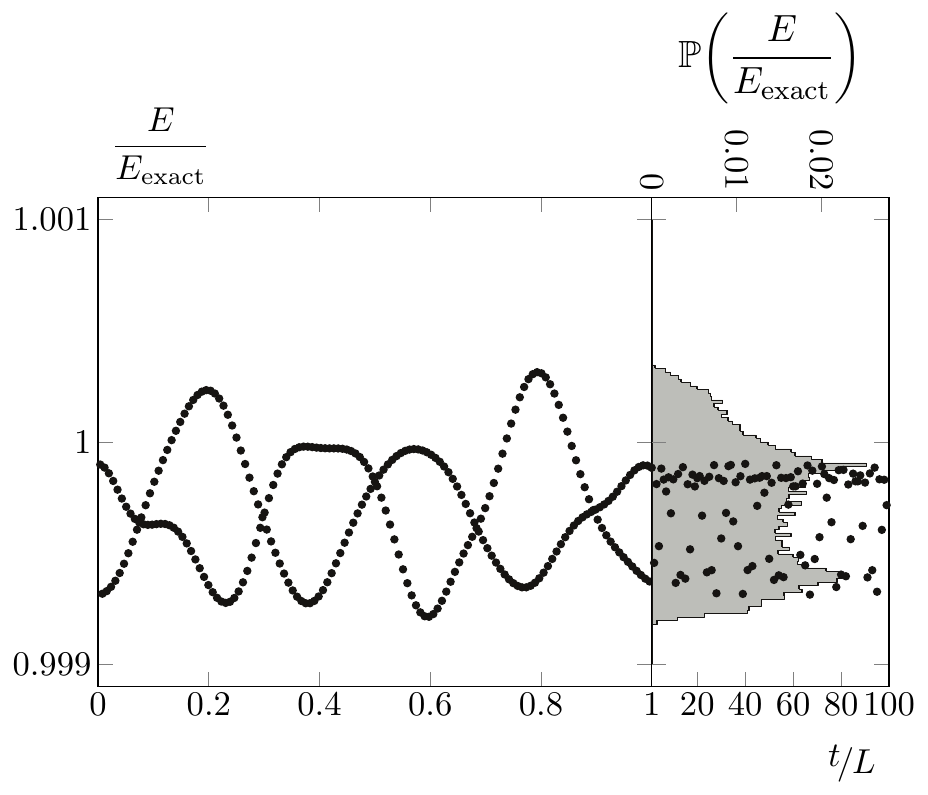}
				\end{center}
				\caption{Time evolution of the total energy for $A=10$ using the \textsc{msilcc} method with
				$\nicefrac{L}{\sqrt{2}\,\delta}=128$.
				The horizontal axis is the same as in \cref{msilccvsother.bddv.energy}.
				The second part of the graph also shows a histogram of the occurrences of the energy during the integration process (from $\nicefrac{\displaystyle{t}}{\displaystyle{L}}=0$ to $\nicefrac{\displaystyle{t}}{\displaystyle{L}}=100$) with $128$ bins uniformly distributed on the interval $\left[0.999,1.001\right]$.}
				\label{msilccvsother.msilcc.energy}
			\end{figure}
			
			Let us finish this very brief presentation of the expected properties of the \textsc{msilcc} method by showing how the total energy (which is here uniquely defined) behaves in time in \cref{msilccvsother.msilcc.energy}.
			We first observe that there seem to be two interlaced curves.
			Actually, this is not the case, there is only one energy that jumps from one carrier curve to the other.
			This ``double'' structure is due to the lattice geometry in combination with the discretization rules.
			More precisely, when the time index is odd there is a shift of the space index and hence the field is not sampled at the same places, leading to a different energy.
			Therefore, there are ``two curves'', one for odd times and the other one for even ones.
			
			Having clarified the effect of the time-discretization we now describe the actual time variation of the total energy.
			Firstly, over short time-scales, the deviations are around $1\text{\textperthousand}$ of its value.
			Secondly, there is no long term trend to increase this deviation.
			Accordingly, these two remarks allow us to promote the \textsc{msilcc} method as a good candidate for the long time integration of conservative field theories.
			In the following we analyse the local conservation laws.
		\subsection{Comparison}\label{msilccvsother.comparison}
			In this Section we compare the performance of the three numerical integrators discussed so far in a complementary way following what we have already discussed in \cref{msilccvsother.preliminaries.testconditions}.
			\subsubsection{Influence of the non-linearity}
\begin{figure}[!htb]
					\begin{center}
						\includegraphics{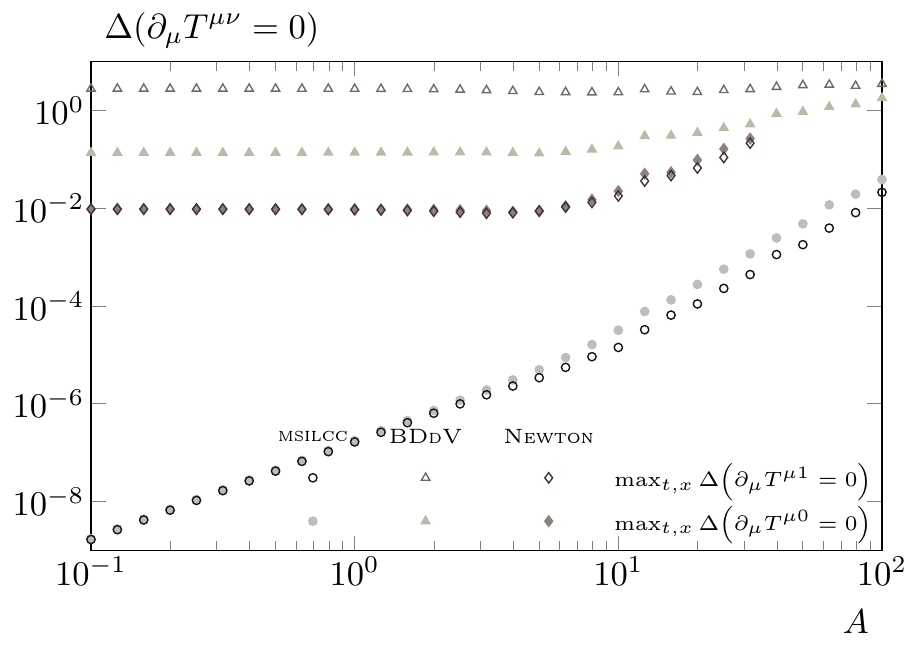}
					\end{center}
					\caption{Error committed by the different methods on the local conservation of the stress-energy tensor as a function of the initial amplitude (discretization is still on $128$ points).}
					\label{msilccvsother.comparison.amplitude.stressenergytensor}
				\end{figure}

\begin{figure}[!htb]
					\begin{center}
						\includegraphics{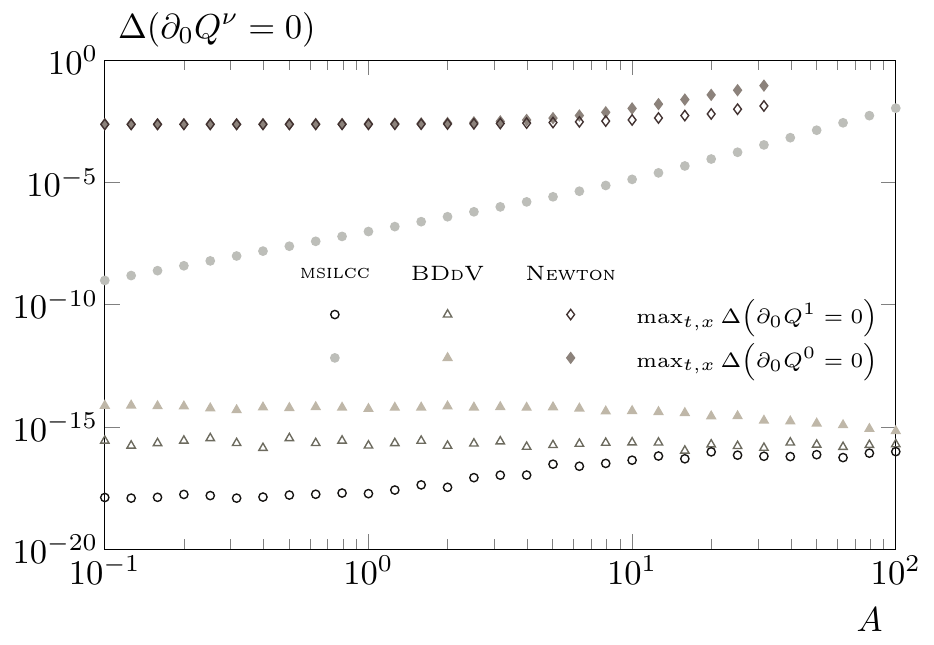}
					\end{center}
					\caption{Error committed by the different methods on the conservation of the charges as a function of the initial amplitude (discretization is still on $128$ points).}
					\label{msilccvsother.comparison.amplitude.charges}
				\end{figure}

				The first situation will explore the influence of the non-linearities (in coordination with the influence of the quality of the lattice spacing).
				\Cref{msilccvsother.comparison.amplitude.stressenergytensor,msilccvsother.comparison.amplitude.charges} represent the error committed on the conservation laws by the different methods as a function of $A$ (the amplitude of the initial condition).
				For this test, the system is integrated up to a time $\nicefrac{\displaystyle{t}}{\displaystyle{L}}=1$ (\emph{i.e.} the solution is obtained over a square), and the error, denoted $\Delta(y=0)$, is taken as the largest deviation from the identity $y=0$ ever encountered (in absolute value and divided by a characteristic quantity, as introduced in \cref{msilccvsother.preliminaries.testconditions}).
				
				Let us start by discussing \cref{msilccvsother.comparison.amplitude.stressenergytensor}.
				We first observe that all methods improve their performance for smaller $A$.
				
				We emphasize that there are some missing data-points for the \textsc{Newton} method (behind $A\sim20$).
				This is due to the fact that the approximation becomes unstable before the final integration time for too strong non-linearity.
				Beyond that point, the solution diverges and the errors as well.
				Behind this feature there is a first important remark: the largest the effect of the non-linearity, and the worse the quality of the sampling, the quickest the \textsc{Newton} approximation becomes unstable.
				This fact is worrying since the parameter region we want to explore is precisely the one in which the non-linearities are relevant.
				Concomitantly, we want to reach long times and it is not desirable to have to oversample the field in time with a too small time spacing.
				
				To pursue the remarks on the \textsc{Newton} method, we stress that it behaves quite well while it remains stable (local errors are between $10^{-2}$ and $10^{-1}$).
				However, it is at minimum $3$ orders of magnitude worse than the \textsc{msilcc} method.
				The difference can be attributed to the fact the \textsc{Newton}'s method is explicit, while the \textsc{msilcc} one is implicit (explicit schemes are known to be worse than implicit ones).
				
				Concerning the \textsc{BDdV} method, the violation of the local conservation laws is very important with an error that ranges between $10^{-1}$ and $10^{+1}$.
				Quite surprisingly, the exact conservation of the energy is only due to the compensation of these large errors once integrated over space.
				
				Finally, the \textsc{msilcc} method produces errors that range from $10^{-9}$ to $10^{-2}$.
				They disappear very abruptly when the non-linearity becomes negligible.
				This is actually due to the fact that the method is exact for a linear problem (as we will show in \cref{implementation.stresstensor.conservation}).
				So, the \textsc{msilcc} method appears, for now, as a very good choice to integrate conservative field theories over long times.
				
				Let us now look at how the errors on the charges behave and, in particular, the energy one (see \cref{msilccvsother.comparison.amplitude.charges}).
				First of all, we remark that the errors committed on the conservation of the energy are in agreement with what we observed earlier when we showed their evolution in time.
				The figure shows that the \textsc{BDdV} prescription is better (actually almost exact since $10^{-14}$ is of the order of the machine precision) than the \textsc{msilcc} method which is itself better than the \textsc{Newton} method.
				The conservation of the second charge is almost exact for both the \textsc{BDdV} and the \textsc{msilcc} schemes and is of the order of the conservation of the energy for the \textsc{Newton} method.
				
				In conclusion, \textsc{Newton}'s method presents a not so bad local conservation of the stress-energy tensor.
				However, by accumulation of these errors, the charges are not conserved and the violation of their conservation, though not too large, is not sufficiently good for high precision measurements.
				Concerning the \textsc{BDdV} method, it presents very poor local conservation properties that, quite surprisingly, lead to a very good conservation of the charges (due to a deceptive cancelation of the errors).
				Finally, the \textsc{msilcc} method behaves more like \textsc{Newton}'s but with much better local conservation features that yield an acceptable conservation of the charges.
				
			\subsubsection{Long time behavior}
\begin{figure}[!htb]
					\begin{center}
						\includegraphics{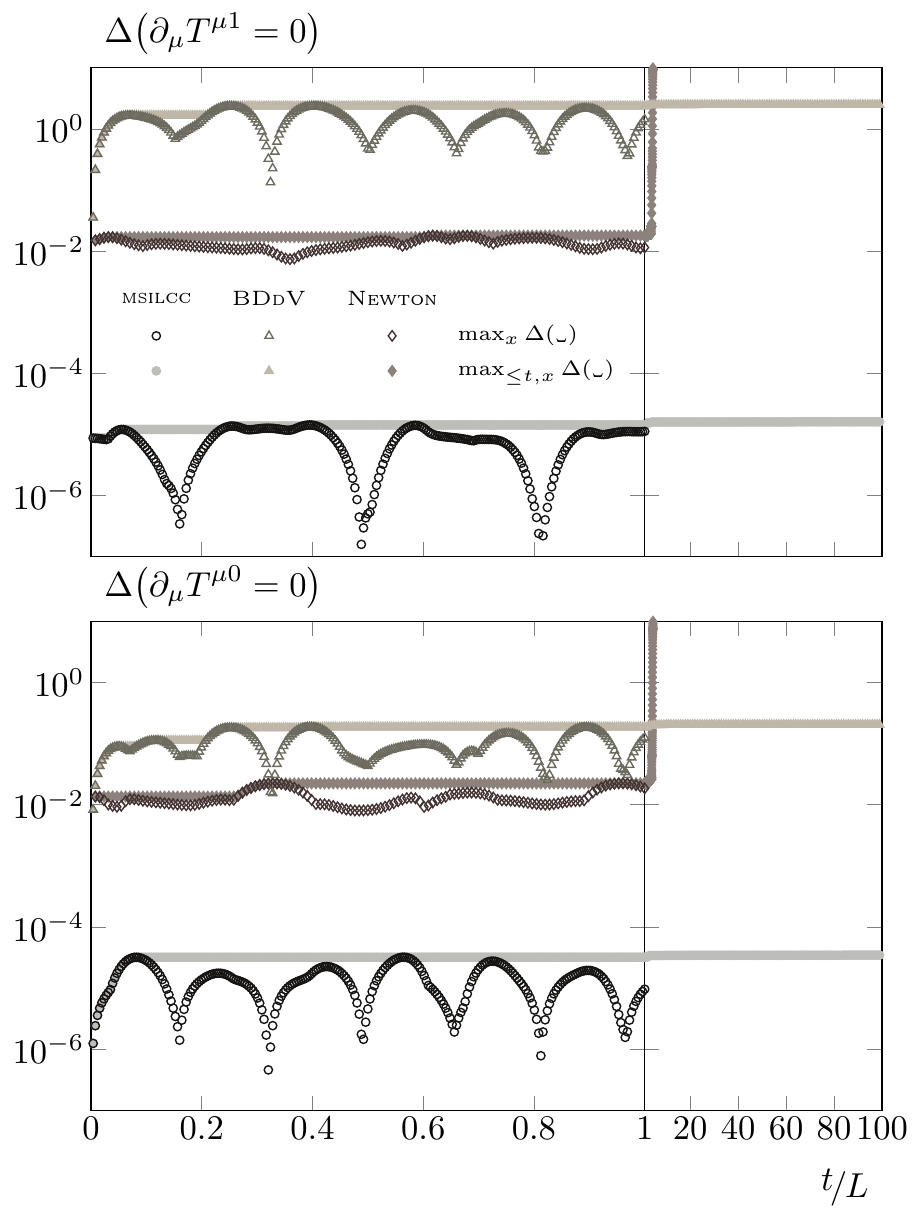}
					\end{center}
					\caption{Long time behavior of the error committed by the different methods on the local conservation of the stress-energy tensor for $A=10$ (discretization is still on $128$ points).
					The upper pair of curves are for the \textsc{BDdV} method (triangles), the intermediate ones for the \textsc{Newton} method (diamonds), and the lower ones for the \textsc{msilcc} method (circles).
					Open and closed symbols show different ways of measuring the error as defined in the text.
					Beyond $\nicefrac{\displaystyle{t}}{\displaystyle{L}}=1$ the horizontal axis is shown in a different linear scale and the curves with open symbols are not plotted since they vary too rapidly with respect to this new time-scale (these represent instantaneous errors that in any case are not relevant on this time-scale).}
					\label{msilccvsother.comparison.time.stressenergytensor}
				\end{figure}

\begin{figure}[!htb]
					\begin{center}
						\includegraphics{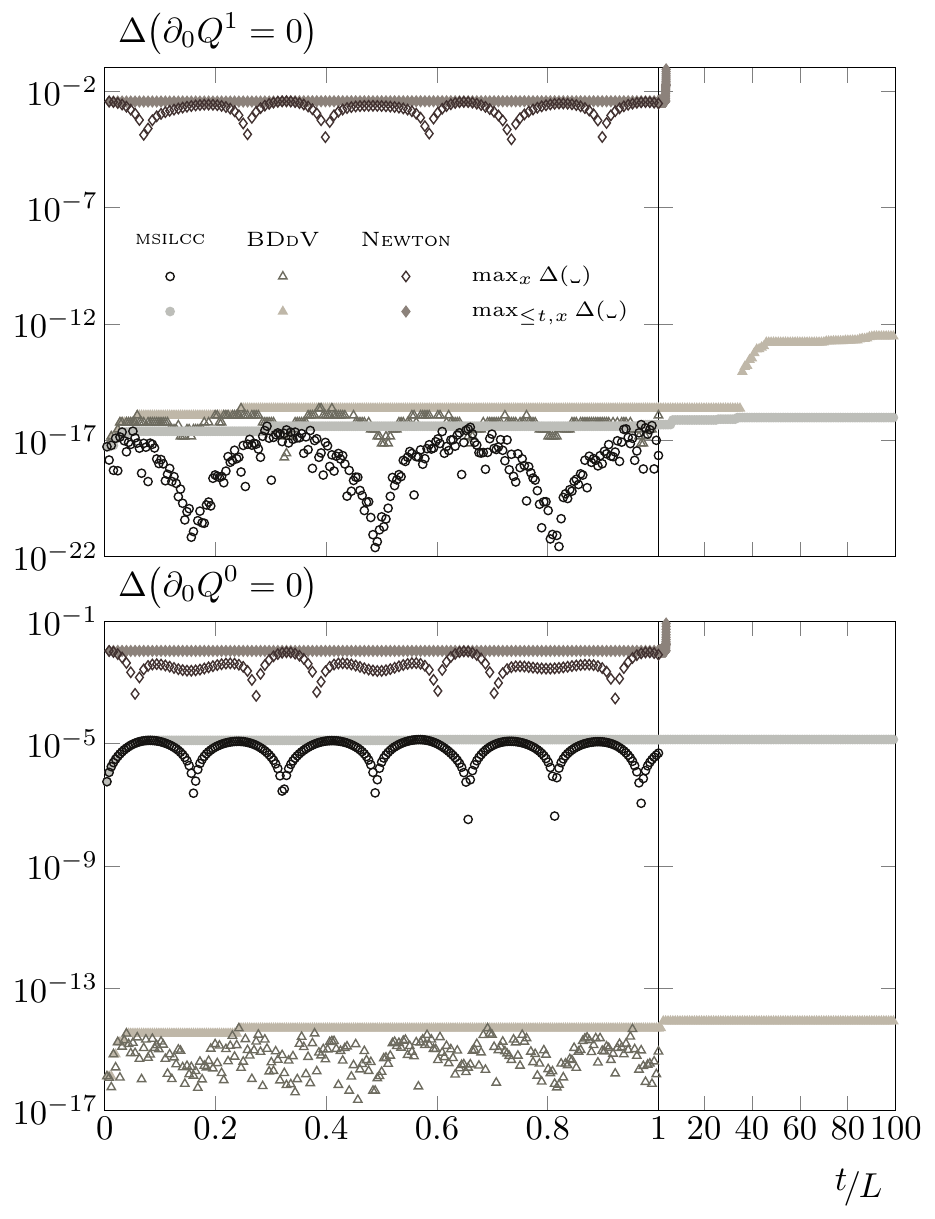}
					\end{center}
					\caption{Long time behavior of the error committed by the different methods on the conservation of the charges for $A=10$ (discretization is still on $128$ points).
					Same symbol convention and time axis as in \cref{msilccvsother.comparison.time.stressenergytensor}.}
					\label{msilccvsother.comparison.time.charges}
				\end{figure}

				We now explore the long time properties.
				\Cref{msilccvsother.comparison.time.stressenergytensor,msilccvsother.comparison.time.charges} represent the error committed on the conservation laws as a function of time.
				Note that for each method the integration is performed using $A=10$ and two errors are displayed: the first one is the largest deviation ever encountered (in absolute value), and the second one is the largest deviation at time $t$.
				The comments made on \cref{msilccvsother.comparison.amplitude.stressenergytensor,msilccvsother.comparison.amplitude.charges} still hold and we will only describe the time behavior here.
				
				Firstly, we observe that the \textsc{Newton} method rapidly becomes unstable (after $\nicefrac{\displaystyle{t}}{\displaystyle{L}}=4$) and is no longer able to describe the evolution of the field.
				
				Secondly, we remark that the instantaneous error evolves in time (over several orders of magnitude).
				So, it is preferable to consider, instead of the instantaneous error, the worst one ever encountered from the beginning.
				
				Finally, the most interesting comment that applies to the \textsc{BDdV} and the \textsc{msilcc} method as well is that the worse errors occur during the short time behavior: after a rapid evolution (of the order of the characteristic time-scale of the system, as we observe by comparing \cref{msilccvsother.comparison.amplitude.stressenergytensor,msilccvsother.comparison.amplitude.charges} with \cref{msilccvsother.model.sols} for $A=10$), the error stabilizes to a value which can be hopefully considered as definitive (of the order of $10^0$ for the \textsc{BDdV} method and $10^{-5}$ for the \textsc{msilcc} one).
				This is particularly true for the local conservation of the stress-energy tensor but less clear for the conservation of the charges even though they seem to reach a constant too.
				
			\subsubsection{Preliminary conclusion}\label{msilccvsother.comparison.conclusion}
				A first element of conclusion is that we need to be extremely wary of methods that possess remarkable properties on some observables but not necessarily the most fundamental objects of the theory.
				
				Within the three methods here presented, the \textsc{msilcc} is the only one that one could trust to integrate a conservative field theory over a long time interval.
				However, its implementation has a cost: the discrete equations of motion are implicit and more expensive  to solve (in terms of computational time) than the other two methods.
				Fortunately, the scheme remains well-defined locally (\emph{i.e.} there is no need to solve the set of algebraic equations globally) and it can be easily scaled to larger volumes and/or extended to theories defined on higher dimensions.
				
\begin{figure}[!htb]
					\begin{center}
						\includegraphics{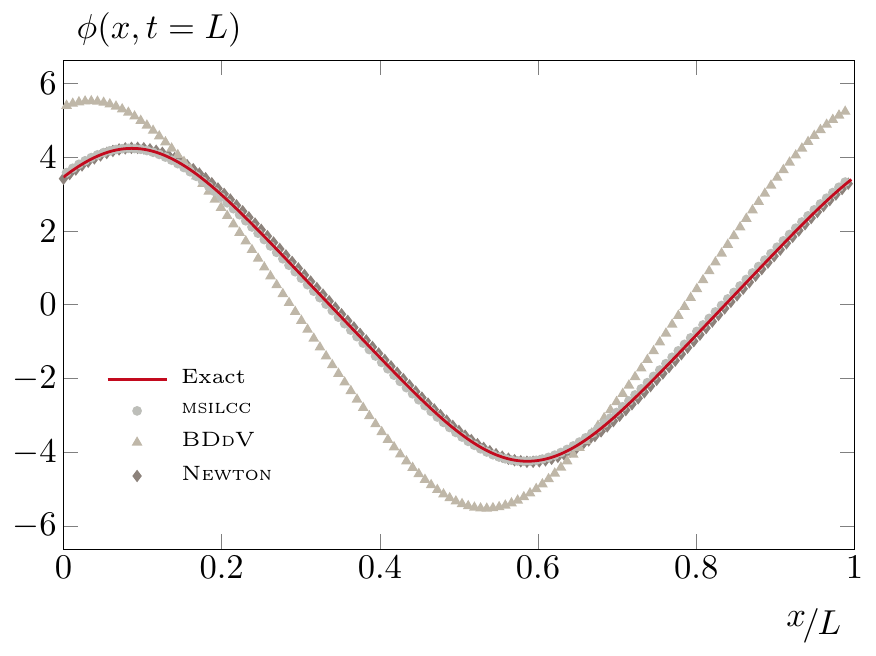}
					\end{center}
					\caption{The red line is the exact solution of \cref{msilccvsother.model.pde} in terms of \textsc{Jacobi} elliptic functions after a time $\nicefrac{\displaystyle{t}}{\displaystyle{L}}=1$ (a particular solution for different initial condition from the ones used so far).
					The datapoints represent the field evolved by the different numerical methods (triangles for \textsc{Newton}, diamonds for \textsc{BDdV} and circles for \textsc{msilcc}).
					On this scale we see no difference between the exact and the numerical solutions obtained with the \textsc{Newton} and \textsc{msilcc} methods.}
					\label{msilccvsother.conclusion.jacobi.evolution}
				\end{figure}

				Up to now, we have eluded the concrete results in term of the (numerical) solution of the \textsc{pde} and one can imagine that all these elements of conservation only have a negligible influence.
				This is not true.
				As an example, over an integration time as short as $\nicefrac{\displaystyle{t}}{\displaystyle{L}}=1$, we observe differences of the order of $1\%$ between the solutions obtained with the different methods (under the same conditions as described in \cref{msilccvsother.model.icbc} and for $A=10$).
				In some situations the differences can become dramatically larger (up to $20\%$) as we show on \cref{msilccvsother.conclusion.jacobi.evolution} that represents the field after an integration time $\nicefrac{\displaystyle{t}}{\displaystyle{L}}=1$, still with periodic boundary conditions, but with an initial state that corresponds to a particular solution of \cref{msilccvsother.model.pde} in term of a \textsc{Jacobi} elliptic function \cite{NIST2010}.
				
			\subsubsection{Symmetry breaking potential}
\begin{figure}[!htb]
					\begin{center}
						\includegraphics{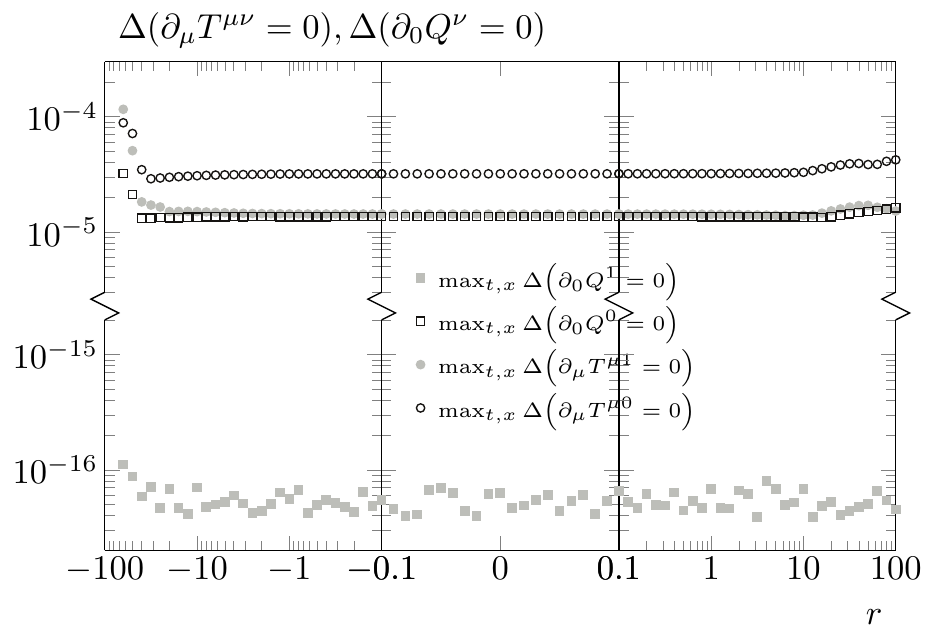}
					\end{center}
					\caption{Error committed by the \textsc{msilcc} method on the conservations of the charges and stress-energy tensor as a function of the parameter that accompanies the quadratic term in the potential: $r$.
					The initial condition follows \cref{msilccvsother.model.ic} with $A=10$ and the errors are accumulated over an integration time of $\nicefrac{\displaystyle{t}}{\displaystyle{L}}=1$ while we still have $\nicefrac{L}{\sqrt{2}\,\delta}=128$.
					The vertical axis is cut between $10^{-5}$ and $10^{-15}$ while the horizontal axis is in logarithmic scale from $-100$ to $-0.1$ and from $0.1$ to $100$ (scale is linear between $-0.1$ and $0.1$).
					Positive values of $r$ mean that the potential has only one minimum at $\phi=0$ while negative values of $r$ mean that the potential is a double well with two minimum at $\phi=\pm\sqrt{-r}$ (so far $r$ was setted to $1$).}
					\label{msilccvsother.conclusion.msilcc.quadfactor.both}
				\end{figure}

				So far we have not considered the influence of $r$ (the parameter that accompanies the quadratic term in the potential) and we want to show that the \textsc{msilcc} method behaves just as well for a potential in double well.
				\Cref{msilccvsother.conclusion.msilcc.quadfactor.both} shows the error committed on the conservations of the charges and the stress-energy tensor as a function of $r$.
				First of all, we observe that the \textsc{msilcc} method has the same conservation properties whatever the shape of the potential.
				Secondly, we remark that the errors do not depend on $r$ (except in the very large $r$ limit).
				This is a direct consequence of a feature that will be proved in \cref{implementation.stresstensor.conservation}: the deviations from the conservation of the stress-energy tensor only arise with the non-linear part of the Hamiltonian.
				Actually we no longer observe this feature in the large $r$ limit since for $A=10$ and for such values of $r$ the field can only ``explore'' the quadratic part of the potential.
				In other words, when $\lr{|r|}$ increases the non-linear part of the potential disappears and the errors are modified.
				
				Now, the rest of this paper will be devoted to the construction of the \textsc{msilcc} method in a more general setup which will be illustrated through the example of the non-linear wave equation.
				In the next section we will introduce a number of necessary concepts and we will discuss how the method is concretely applied, while in the last section we will finally present its construction.
	\section{Applicability and preliminaries}\label{preliminaries}
		This section introduces the concept of symplecticity (for mechanical systems) and its generalization to multi-symplecticity (for field theories).
		We then introduce the idea of Hamiltonian \textsc{pde}s through the \textsc{De Donder}~-- \textsc{Weyl} formalism.
		Afterwards, we show how to rewrite the equation of motion (when possible) in a way that allows one to use the \textsc{msilcc} method.
		Finally, we will exhibit (both in a general setup and with an example) the multi-symplectic structure, prove its conservation, and discuss the definition and properties of the stress-energy tensor.
		
		Most of the points presented in this section are just reminders except for two of them.
		Firstly, the link between the \textsc{De Donder}~-- \textsc{Weyl} formulation and multi-symplectic geometry is not so common (usually, the \textsc{De Donder}~-- \textsc{Weyl} formulation is treated through the formalism of the poly-symplectic geometry).
		Secondly, as far as we know, the discussion of the degeneracy of the multi-symplectic structure, and especially its resolution, is completely new.
		\subsection{From symplectic to multi-symplectic structure}
			\subsubsection{Emergence of the symplectic structure}
				Let us consider a mechanical problem, constituted by $N$ particles, in its Hamiltonian formulation (for an exhaustive presentation of the Hamiltonian formalism see, for example,~\cite{JoseSaletan1998}).
				The $i^{\text{th}}$ body is described by its position and momenta, or \textsc{Darboux} coordinates, $\lr{(q^i,p_i)}$, and the dynamics are fully characterized by the Hamiltonian, $H\lr{(\lr{\{q^i,p_i\}})}$.
				The phase space, $\Omega$, is a differential manifold locally parametrized by the union of the \textsc{Darboux} coordinates.
				
				The symplectic character of the phase space of a mechanical system is rather well known, so we will only remind its central aspects~\cite{JoseSaletan1998}.
				The symplectic manifolds are differentiable manifolds equipped with a closed non-degenerate 2-form.
				We now recall how this 2-form is constructed.
				
				Let us describe $\Omega$ using the generalized coordinates $\lr{\{\zeta^a\}}\equiv\lr{\{q^i,p_i\}}$ with $a\in\lr{\llbracket1,2N\rrbracket}$ (the union of the positions and momenta).
				We define $\lr{\{\bm{\partial}_a=\nicefrac{\displaystyle{\partial}}{\displaystyle{\partial\zeta^a}}\}}$ and $\lr{\{\db[a]{}=\d{\zeta^a}\}}$, respectively, as a basis for $T\Omega$, the tangent space of $\Omega$, and $T^*\Omega$, its dual (see~\cite{Nakahara2003} for the necessary concepts of geometry).
				
				Now, we provide to $\Omega$ a $2$-form
				\begin{equation*}
					\bm{\omega}=\omega_{ab}\,\db[a]{}\wedge\db[b]{}\eqpc
				\end{equation*}
				where we have used the \textsc{Einstein} summation convention, as we will do (unless explicitly stated) in the rest of this paper.
				In particular, for a mechanical system and in \textsc{Darboux} coordinates this $2$-form reads
				\begin{equation*}
					\bm{\omega}=\d{q^i}\wedge\d{p_i}\eqpd
				\end{equation*}
				It is obviously closed ($\d{\bm{\omega}}=0$), and it can be shown to be non-degenerate ($\det{\bm{\omega}}\neq0$, the proof will be given later).
				Therefore, $\lr{(\Omega,\bm{\omega})}$ defines a symplectic manifold where $\bm{\omega}$ is the symplectic structure.
				
				The conservation of the symplectic structure under a Hamiltonian flow leads to \textsc{Liouville}'s theorem (since the volume form on $\Omega$ is obtained from $\bm{\omega}$).
				
				The \textsc{Poisson} bracket is directly related to $\bm{\omega}$ since for any pair of differentiable functions of the phase space variables, $\ast\lr{(\lr{\{\zeta^a\}})}$ and $\diamond\lr{(\lr{\{\zeta^a\}})}$, we have
				\begin{equation}
					\lr{\{\ast,\diamond\}}=\bm{\omega}\lr{(\bm{\chi}_\ast,\bm{\chi}_\diamond)}=\omega^{\lr{[ab]}}\,\bm{\partial}_a\ast\bm{\partial}_b\diamond\eqpc\label{preliminaries.sympform.poissonbracket}
				\end{equation}
				where $\bm{\chi}_\ast$ is the vector field associated to $\ast$ by $\bm{\omega}$ defined as
				\begin{equation}
					\bm{\omega}\lr{(\bm{\chi}_\ast,\cdot)}=\d{\!\ast}\lr{[\cdot]}\eqpc\label{preliminaries.sympform.vectorfield}
				\end{equation}
				and where $[ab]$ denotes the antisymmetric combination of the indices
				\begin{equation}
					\omega^{\lr{[ab]}}=\frac{\omega^{ab}-\omega^{ba}}{2}\eqpd\label{preliminaries.sympform.antisymmetric2indices}
				\end{equation}
				Actually, \cref{preliminaries.sympform.poissonbracket,preliminaries.sympform.vectorfield} mean that the \textsc{Poisson} bracket (parametrized by $\omega^{ab}$) is the inverse of the symplectic $2$-form (parametrized by $\omega_{ab}$).
				In particular, for a mechanical system and in \textsc{Darboux} coordinates the \textsc{Poisson} bracket reads
				\begin{equation*}
					\lr{\{\ast,\diamond\}}=\frac{\partial\,\ast}{\partial p_i}\frac{\partial\,\diamond}{\partial q^i}-\frac{\partial\,\ast}{\partial q^i}\frac{\partial\,\diamond}{\partial p_i}\eqpd
				\end{equation*}
				
				Symplecticity is an essential structure of the phase space of mechanical systems and numerical methods that preserve it are very popular~\cite{MclachlanQuispel2006}.
				For field theories, the symplectic structure needs to be generalized if we want to develop the same class of methods for \textsc{pde}s.
			\subsubsection{Generalization to the multi-symplectic structure}
				A multi-symplectic manifold, $\lr{(\Omega,\lr{\{\bm{\omega}^\mu\}})}$, is a differential manifold equipped with several independent symplectic structures\footnote{Let us complete this definition with some vocabulary remarks.
				One can also encounter in the literature the concept of poly-symplectic (or $n$-plectic) manifold which is actually different: a $n$-plectic manifold is a differential manifold with a closed, non-degenerate, $\lr{(n+1)}$-form (poly-symplectic is for any $n>1$ while symplectic stands for $1$-plectic).
				So, a multi-symplectic manifold is poly-symplectic too since we can define $\bm{\varpi}=\bigwedge_\mu\bm{\omega}^\mu$, but the reverse is not necessarily true.
				Finally, we want to stress that the meaning of poly-symplectic and multi-symplectic can be exchanged depending on authors.} (in our case, one per space-time direction).
				As we shall see later, this structure is natural for a field theory subject to a conservation law that fits the invariance under the Hamiltonian flow.
				
				The next section will illustrate how the multi-symplectic aspect of the phase space of a field theory emerges.
				We will first recall how to construct a Hamiltonian formulation of a field theory.
				Then, we will identify the phase space, highlight the underlying multi-symplectic structure, and establish its conservation.
				Finally, we will discuss the properties of the stress-energy tensor.
		\subsection{The \texorpdfstring{\textsc{De Donder}~-- \textsc{Weyl}}{De Donder Weyl} (DW) Hamiltonian formulation of field theories}
			\subsubsection{From Lagrangian to DW Hamiltonian formulation}
				Let us start with a space-time Lorentzian manifold, $\mathcal{M}$, of dimension $D=1+d$.
				We assume $\mathcal{M}$ to be non dynamic (the metric is not subject to an equation of motion).
				Without lose of generality, $\mathcal{M}$ will be supposed to be flat with metric $\eta\equiv\diag{\lr{(1,-1,\cdots,-1)}}$.
				We define a local coordinate system $\lr{\{x^\mu\}}$, with $\lr{\{\partial_\mu=\nicefrac{\displaystyle{\partial}}{\displaystyle{\partial x^\mu}}\}}$ a basis of $T\mathcal{M}$, and $\mu\in\lr{\llbracket0,d\rrbracket}$.
				
				Next, we consider a field theory on $\mathcal{M}$ described by the action $\mathcal{S}\lr{[\lr{\{\phi^i\}}]}$, where $\lr{\{\phi^i\}}$ with $i\in\lr{\llbracket1,\mathcal{N}\rrbracket}$ is a collection of dynamic fields.
				This action originates in a Lagrangian density $\mathcal{L}$ (which is assumed to depend only on the field and its first derivatives) and reads
				\begin{equation*}
					\mathcal{S}=\int\d[D]{x}\;\mathcal{L}\lr{(\lr{\{\phi^i,\partial_\mu\phi^i\}})}\eqpc
				\end{equation*}
				where $\d[D]{x}$ is the measure over $\mathcal{M}$ and $\mu\in\lr{\llbracket0,d\rrbracket}$.
				The stationarity of $\mathcal{S}$ leads to the \textsc{Euler}~-- \textsc{Lagrange} equations (\emph{i.e.} the equations of motion) for the fields
				\begin{equation*}
					\frac{\delta\mathcal{S}}{\delta\phi^i}=0=
					\frac{\partial\mathcal{L}}{\partial\phi^i}-\partial_\mu\frac{\partial\mathcal{L}}{\partial\lr{(\partial_\mu\phi^i)}}\eqpd
				\end{equation*}
				This is the Lagrangian formulation of a field theory.
				
				The idea of the Hamiltonian formulation of classical mechanics is to substitute the generalized speed ($\dot{q}^i$) by a conjugate momentum ($p_i$).
				Obviously the same reasoning can be applied to a field theory but, unfortunately, it breaks the \textsc{Lorentz} covariance of the theory.
				The idea of \textsc{De Donder} and \textsc{Weyl} is to reestablish the covariance by introducing one conjugate momentum per direction of space-time (such that they are treated on an equal footing).
				Thus, we define
				\begin{equation*}
					{\psi_i}^\mu=\frac{\partial\mathcal{L}}{\partial\lr{(\partial_\mu\phi^i)}}\eqpc
				\end{equation*}
				as the conjugate momentum of the field $\phi^i$ along the $\mu^{\text{th}}$ direction of $\mathcal{M}$.
				Then, provided that the following \textsc{Legendre} transform is not singular,
				\begin{equation*}
					\mathcal{H}={\psi_i}^\mu\,\partial_\mu\phi^i-\mathcal{L}\eqpc
				\end{equation*}
				defines the \textsc{De Donder}~-- \textsc{Weyl} Hamiltonian density.
				Henceforth, the unknowns are the fields ($\lr{\{\phi^i\}}$) with their conjugate momenta in each direction of space-time ($\lr{\{{\psi_i}^\mu\}}$).
				Together, they are the local coordinates of a differential manifold, $\Omega$, which is the \textsc{De Donder}~-- \textsc{Weyl} definition of phase space (a multi-symplectic manifold as we will prove later).
				The \textsc{Hamilton} equations generalize to
				\begin{subequations}
					\begin{align}
						\partial_\mu{\psi_i}^\mu&=-\frac{\partial\mathcal{H}}{\partial\phi^i}\eqpc\label{preliminaries.dwh.hamilteq.moments}\\
						\partial_\mu\phi^i&=\frac{\partial\mathcal{H}}{\partial{\psi_i}^\mu}\eqpd\label{preliminaries.dwh.hamilteq.field}
					\end{align}
				\end{subequations}
				This is the DW Hamiltonian formulation of a field theory (see~\cite{Kanatchikov1993,PauflerRomer2002,vonRieth1984} for a detailed review of this formulation of classical field theories).
				
				\Cref{preliminaries.dwh.hamilteq.moments,preliminaries.dwh.hamilteq.field} can be rewritten in a more symmetrical way (as in the case of classical mechanics) by introducing a \textsc{Poisson} bracket for each space-time direction.
				This collection of \textsc{Poisson} brackets is intimately related to the multi-symplectic structure and each bracket is associated to a symplectic $2$-form.
				
				Let us now describe $\Omega$ in the generalized coordinates $\lr{\{\zeta^a\}}\equiv\lr{\{\phi^i\}}\cup\lr{\{{\psi_i}^\mu\}}$ where $a$ is an index conveniently chosen to sweep the collection.
				Following~\cite{Bridges1997,BridgesReich2001,BridgesReich2006}, we define on $T\Omega$ the vector $\bm{\zeta}=\zeta^a\bm{\partial}_a$ that allows one to rewrite \cref{preliminaries.dwh.hamilteq.moments,preliminaries.dwh.hamilteq.field} in the abstract form
				\begin{equation}
					\bm{M}^\mu\cdot\partial_\mu\bm{\zeta}=\bm{\nabla}\mathcal{H}\eqpc\label{preliminaries.dwh.hamilteq.abstractform}
				\end{equation}
				where $\lr{\{\bm{M}^\mu\}}$ is a set of constant skew-symmetric matrix of $T\Omega$.
				The gradient is defined as
				\begin{align*}
					\bm{\nabla}\colon\mathcal{C}^\infty\lr{(\Omega)}&\to T\Omega\\
					\bm{\nabla}\mathcal{H}&\mapsto {\lr{(\d{\mathcal{H}})}}^*=\delta^{ab}\,\bm{\partial}_b\mathcal{H}\,\bm{\partial}_a\eqpc
				\end{align*}
				and we recall that $\lr{\{\bm{\partial}_a\}}$ is a basis of $T\Omega$.
				
				The multi-symplectic structure will directly arise from \cref{preliminaries.dwh.hamilteq.abstractform}.
				We will come back to this point later but we already emphasize that, as long as a \textsc{pde} can be written in such a form, the \textsc{msilcc} method will apply.
				
				The next step is to exhibit the multi-symplectic structure and to prove its conservation.
				This will be done in the next section, but let us first clarify the procedure presented above on an example.
			\subsubsection{The non-linear wave equation}
				We consider the dynamics of a real scalar field, $\phi$, whose Lagrangian density is given by
				\begin{equation*}
					\mathcal{L}=\frac{1}{2}\,\partial_\mu\phi\,\partial^\mu\phi-V\lr{(\phi)}\eqpd
				\end{equation*}
				The \textsc{Euler}~-- \textsc{Lagrange} equation reads
				\begin{equation*}
					\partial_\mu\partial^\mu\phi+V'\lr{(\phi)}=\square\phi+V'\lr{(\phi)}=0\eqpd
				\end{equation*}
				We introduce now
				\begin{equation*}
					\psi^\mu=\frac{\partial\mathcal{L}}{\partial\lr{(\partial_\mu\phi)}}=\partial^\mu\phi\eqpc
				\end{equation*}
				which is the conjugate momentum of $\phi$ in the $\mu^{\text{th}}$ direction.
				Then, the DW Hamiltonian density reads
				\begin{equation*}
					\mathcal{H}=\psi^\mu\,\partial_\mu\phi-\mathcal{L}=\frac{1}{2}\,\psi_\mu\,\psi^\mu+V\lr{(\phi)}\eqpd
				\end{equation*}
				Defining
				\begin{equation*}
					\bm{\zeta}^{\text{T}}=\begin{bmatrix}
						\phi&\psi^0&\psi^1&\cdots&\psi^d
					\end{bmatrix}\eqpc
				\end{equation*}
				the equation of motion is given by \cref{preliminaries.dwh.hamilteq.abstractform} provided that
				\begin{equation*}
					{M^\mu}_{ab}=\delta_{a-1}^{\mu\vphantom{0}}\,\delta_b^{\vphantom{\mu}0}-\delta_{b-1}^{\mu\vphantom{0}}\,\delta_a^{\vphantom{\mu}0}\eqpc
				\end{equation*}
				where $a,b\in\lr{\llbracket0,D\rrbracket}$.
				
				In dimension $D=1+0$ ($d=0$), aliasing $q\equiv\phi$, $p\equiv\psi^0$ and $H\equiv\mathcal{H}$, we recover the expected \textsc{Hamilton} equation of motion for a mechanical problem
				\begin{equation*}
					\begin{bmatrix}
						0&-1\\
						1&0
					\end{bmatrix}\cdot\begin{bmatrix}
						\dot{q}\\
						\dot{p}
					\end{bmatrix}=\begin{bmatrix}
						\displaystyle{\frac{\partial H}{\partial q}}\\[1em]
						\displaystyle{\frac{\partial H}{\partial p}}
					\end{bmatrix}=\begin{bmatrix}
						V'\lr{(q)}\\
						p
					\end{bmatrix}\eqpd
				\end{equation*}
				
				In dimension $D=1+1$ ($d=1$) the $\bm{M}$ matrices read
					\begin{align*}
						\bm{M}^0&=\begin{bmatrix}
							0&-1&0\\
							1&0&0\\
							0&0&0
						\end{bmatrix} \eqpc\\
						\bm{M}^1&=\begin{bmatrix}
							0&0&-1\\
							0&0&0\\
							1&0&0
						\end{bmatrix}\eqpc
					\end{align*}
				and we stress that, for both of them, the eigenvalues are $0$ and $\pm i$.
				
				In the general case all the $\bm{M}$ matrices have the same eigenvalues: $\pm i$ and $0$ ($d$ times degenerate).
				As we will see in the next section this fact causes some difficulties.
				For the moment, we are going to highlight the multi-symplectic structure of $\Omega$ relying on \cref{preliminaries.dwh.hamilteq.abstractform}.
				Then, we will return to the vanishing eigenvalues and we will explain how to treat them in the particular case of the non-linear wave equation.
			\subsubsection{The multi-symplectic structure}
				This section will be devoted to the construction of the multi-symplectic structure on $\Omega$ and, again following~\cite{Bridges1997,BridgesReich2001,BridgesReich2006}, we are going to show how it directly emerges from \cref{preliminaries.dwh.hamilteq.abstractform}.
				
				In the previous section we have shown that for a single scalar field theory the \textsc{Hamilton} equation is fully characterized by a
				Hamiltonian density and a collection of $D$ (constant and skew-symmetric) matrices ($\lr{\{\bm{M}^\mu\}}$) that can be used to define $D$ $2$-forms $\lr{\{\bm{\omega}^\mu\}}$:
				\begin{equation*}
					\bm{\omega}^\mu=-\frac{1}{2}\,{M^\mu}_{ab}\,\db[a]{}\wedge\db[b]{}\eqpd
				\end{equation*}
				The action of these $2$-forms on a pair of vectors is
				\begin{equation*}
					\bm{\omega}^\mu\lr{(\bm{\ast},\bm{\diamond})}=\lr{<\bm{M}^\mu\cdot\bm{\ast},\bm{\diamond}>}=-\lr{<\bm{\ast},\bm{M}^\mu\cdot\bm{\diamond}>}\eqpc
				\end{equation*}
				for any $\bm{\ast},\bm{\diamond}\in T\Omega$, where $\lr{<\text{\textvisiblespace},\text{\textvisiblespace}>}$ is the scalar product on $T\Omega$.
				
				All the $\lr{\{\bm{\omega}^\mu\}}$ are closed ($\d{\bm{\omega}^\mu}=0$) since all the $\lr{\{\bm{M}^\mu\}}$ are independent of the fields.
				On the other hand, let us assume for the moment (we will come back to this at the end of this section) that they are all non-degenerate.
				Finally, since the $\lr{\{\bm{M}^\mu\}}$ are linearly independent (as we will shortly show) the $\lr{\{\bm{\omega}^\mu\}}$ are linearly independent as well.
				So, $\lr{(\Omega,\lr{\{\bm{\omega}^\mu\}})}$ is a multi-symplectic manifold.
				
				Let us first prove the independence of the $\lr{\{\bm{M}^\mu\}}$ matrices.
				In order to lighten the following computation, let us introduce space-time indices which behave differently under the \textsc{Einstein} summation rule.
				From now, the $\varrho$ index takes just one value and does not imply summation (even if repeated), while the $\sigma$ index ($\sigma \neq \varrho$) behaves in the standard way.
				All other indices are unaffected.
				Let us now start by supposing that one of the $\lr{\{\bm{M}^\mu\}}$ is linearly dependent on the others: $\bm{M}^\varrho=\alpha_\sigma\bm{M}^\sigma$.
				Then, in \cref{preliminaries.dwh.hamilteq.abstractform}, the operator on the left hand side becomes
				\begin{equation*}
					\begin{split}
						\bm{M}^\mu\cdot\partial_\mu&=\bm{M}^\varrho\cdot\partial_\varrho+\bm{M}^\sigma\cdot\partial_\sigma\\
						&=\alpha_\sigma\bm{M}^\sigma\cdot\partial_\varrho+\bm{M}^\sigma\cdot\partial_\sigma\\
						&=\bm{M}^\sigma\cdot\lr{(\alpha_\sigma\partial_\varrho+\partial_\sigma)}\\
						&=\bm{M}^\sigma\cdot\tilde{\partial}_\sigma\eqpc
					\end{split}
				\end{equation*}
				where the $\tilde{}$ refers to a different coordinate system.
				As the direction $\tilde{\partial}_\varrho$ has just disappeared from the differential operator, it means that the dynamics along this direction are trivial.
				Thus the $\lr{\{\bm{M}^\mu\}}$ are linearly independent.
				
				To complete the proof of the multi-symplecticity of $\Omega$ it remains to discuss the question of the degeneracy of the $\lr{\{\bm{\omega}^\mu\}}$.
				The $2$-forms will be non-degenerate as long as they all satisfy
				\begin{equation*}
					\det{\bm{\omega}^\mu}=\det{\bm{M}^\mu}\neq0\eqpd
				\end{equation*}
				In other words the $\lr{\{\bm{M}^\mu\}}$ matrices should not have any zero eigenvalue.
				Nevertheless we have seen in the previous section that for the example of the non-linear wave equation all the $\lr{\{\bm{M}^\mu\}}$ matrices have the same eigenvalues and especially a $d$ times degenerate zero.
				So, as it is, the phase space of the non-linear wave equation is multi-symplectic only for $d=0$ (that is to say, for the mechanical problem).
				
				We are now going to present the resolution of this problem on the particular example of the non-linear wave equation.
				The same reasoning can be applied to other theories.
				
				Let us first recall that the aim of the DW Hamiltonian formalism was to preserve covariance.
				Thus, it treats all the space-time directions on an equal footing by introducing a conjugate momentum for each.
				Still, time is not space, and the conjugate momentum in time will have a different status.
				We will call it the canonical one.
				
				Degeneracy comes from the existence of non-canonical conjugate fields.
				So, the idea is to ensure that all the conjugate momenta be canonical (\emph{i.e.} the conjugate momentum along time of a dynamic field).
				The solution is to add extra fields, interacting with each other, such that the new field theory allows a conjugate momentum to be shared by several fields and thus to be canonical for someone.
				Finally, we still want to preserve covariance and, for each direction of space-time, each momentum needs to be the conjugate of a field.
				
\begin{figure}[!htb]
					\begin{center}
						\includegraphics{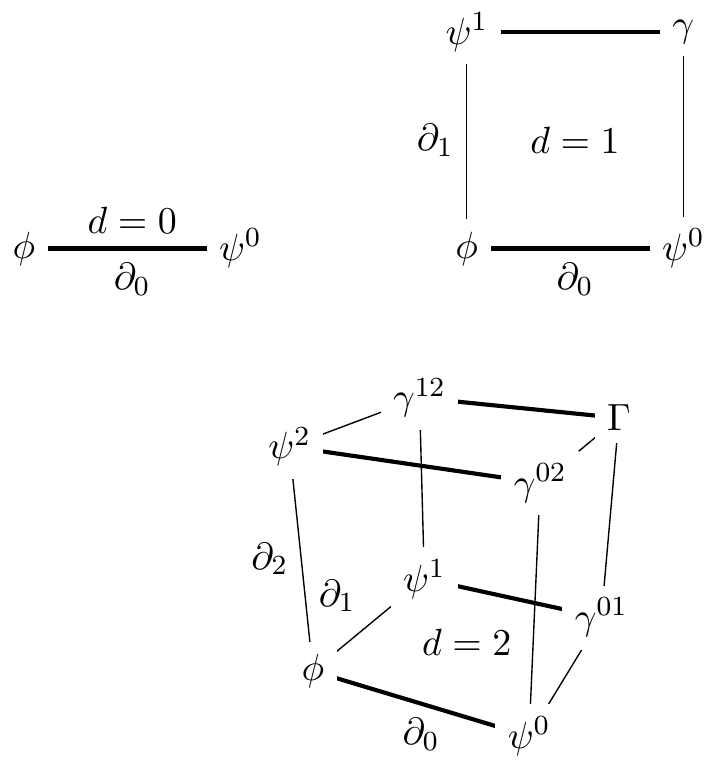}
					\end{center}
					\caption{Graphic representation of the construction of a non-degenerate theory in dimension $D=1+0$, $=1+1$ and $=1+2$ (the canonical direction is shown with a thick line while the others are not).
					The fields in the $d=1$ case will be introduced in \cref{preliminaries.1dfields.def.phi,preliminaries.1dfields.def.psi0,preliminaries.1dfields.def.psi1,preliminaries.1dfields.def.gamma}.}
					\label{preliminaries.dwh.nodegfields.graphic}
				\end{figure}

				Putting these elements together, we can modify the Lagrangian density such that it describes an equivalent, but non-degenerate, problem.
				This construction can be graphically represented by placing both the fields and their conjugate momenta on the vertices of a $D$-dimensional hypercube: each direction stands for a space-time direction and a line means that one of the fields on the edges is the conjugate momentum (along this direction) of the other.
				So, on any path, a dynamic field alternates with a conjugate field.
				We have drawn it for the three first dimensions in \cref{preliminaries.dwh.nodegfields.graphic}.
				
				Now, to construct the new field theory explicitly, we first consider the collection of fields
				\begin{equation}
					{\lr{\{{\lr{\{\Phi^{^{\lr{(i)}}}_{\mu_1\cdots\mu_i}=\Phi^{^{\lr{(i)}}}_{\lr{[\mu_1\cdots\mu_i]}}\}}}_{\mu_1,\cdots,\mu_i\in\lr{\llbracket0,d\rrbracket}}\}}}_{i\in\lr{\llbracket0,D\rrbracket}}\eqpc\label{preliminaries.dwh.nodegfields.collection}
				\end{equation}
				where the pair of square brackets denote the anti-symmetrization defined as
				\begin{equation*}
					\diamond_{\lr{[\mu_1\cdots\mu_i]}}=\frac{1}{i!}\sum_{p\in\mathcal{P}_i\mathrlap{\lr{(\lr{\{\mu_1\cdots\mu_i\}})}}}\sigma\lr{(p)}\diamond_p\eqpc
				\end{equation*}
				where $\mathcal{P}_i$ is the symmetric group of $i$ symbols and $\sigma\lr{(p)}$ is the signature of $p$. In particular, \cref{preliminaries.sympform.antisymmetric2indices} gives this definition for two indices.
				
				Before going further, we stress that the $\Phi^{^{\lr{(i)}}}$ will not be treated as a tensor field but as a collection of scalar fields, conveniently assembled in the same object.
				The collection (\ref{preliminaries.dwh.nodegfields.collection}) is separated in two (equal) parts: $\{\Phi^{^{\lr{(2i)}}}\}$ contains the dynamical fields, while $\{\Phi^{^{\lr{(2i+1)}}}\}$ contains the conjugate fields.
				As a final remark, this collection is composed of
				\begin{equation*}
					\sum_{i=0}^D\binom{D}{i}=2^D
				\end{equation*}
				elements and will indeed be suitable to populate the vertices of the $D$-dimensional hypercube introduced earlier as a graphic representation of this construction.
				
				Next, we consider the Lagrangian density
				\begin{equation*}
					\begin{split}
						\mathcal{L}_0&=\frac{1}{2}\sum_{j=0}^{\leq\nicefrac{D}{2}}\frac{1}{\lr{(2j)}!}\,\partial^{^{\vphantom{\lr{(2j)}}}}_{\nu_0}\Phi^{^{\lr{(2j)}}}_{\nu_1\cdots\nu_{2j}}\,\partial_{_{\vphantom{\lr{(2j)}}}}^{\nu_0}\Phi_{_{\lr{(2j)}}}^{\nu_1\cdots\nu_{2j}}\\
						&\quad-\sum_{j=0}^{\leq\nicefrac{D}{2}}\frac{1}{\lr{(2j)}!}\,V_{_{\lr{(2j)}}}^{\lr{[\nu_1\cdots\nu_{2j}]}}\lr{(\Phi^{^{\lr{(2j)}}}_{\nu_1\cdots\nu_{2j}})}\eqpc
					\end{split}
				\end{equation*}
				where the $\nu$ indices are in $\lr{\llbracket0,d\rrbracket}$, which is nothing else than the concatenation of $2^d$ independent theories\footnote{The anti-symmetrization of the indices in the potential are required since the fields are anti-symmetric objects and in particular $V_{_{\lr{(2)}}}^{01}\lr{(\Phi^{^{\lr{(2)}}}_{01})}$ should coincide with $V_{_{\lr{(2)}}}^{10}\lr{(\Phi^{^{\lr{(2)}}}_{10})}$.} (note that all these theories need to belong in the class of the non-linear wave equation but it is not required that they have the same potential).
				Then we add to this Lagrangian two vanishing coupling terms.
				First, we obviously have
				\begin{equation*}
					\begin{split}
						0=&+\frac{1}{2}\sum_{j=1}^{\leq\nicefrac{D}{2}}\frac{1}{\lr{(2j-1)}!}\,\Phi^{^{\lr{(2j)}}}_{\mu\nu_1\cdots\nu_{2j-1}}\,\partial_{_{\vphantom{\lr{(2j)}}}}^\mu\partial^{^{\vphantom{\lr{(2j)}}}}_\nu\Phi_{_{\lr{(2j)}}}^{\nu\nu_1\cdots\nu_{2j-1}}\\
						&-\frac{1}{2}\sum_{j=1}^{\leq\nicefrac{D}{2}}\frac{1}{\lr{(2j-1)}!}\,\Phi^{^{\lr{(2j)}}}_{\mu\nu_1\cdots\nu_{2j-1}}\,\partial_{_{\vphantom{\lr{(2j)}}}}^\mu\partial^{^{\vphantom{\lr{(2j)}}}}_\nu\Phi_{_{\lr{(2j)}}}^{\nu\nu_1\cdots\nu_{2j-1}}\eqpc
					\end{split}
				\end{equation*}
				and then, using an integration by parts (and ignoring the boundary terms), we claim that
				\begin{equation*}
					\begin{split}
						\mathcal{L}_+&=-\frac{1}{2}\sum_{j=1}^{\leq\nicefrac{D}{2}}\frac{1}{\lr{(2j-1)}!}\,\partial_{_{\vphantom{\lr{(2j)}}}}^\mu\Phi_{_{\lr{(2j)}}}^{\nu\nu_1\cdots\nu_{2j-1}}\,\partial^{^{\vphantom{\lr{(2j)}}}}_\nu\Phi^{^{\lr{(2j)}}}_{\mu\nu_1\cdots\nu_{2j-1}}\\
						&\quad+\frac{1}{2}\sum_{j=1}^{\leq\nicefrac{D}{2}}\frac{1}{\lr{(2j-1)}!}\,\partial^{^{\vphantom{\lr{(2j)}}}}_\nu\Phi_{_{\lr{(2j)}}}^{\nu\nu_1\cdots\nu_{2j-1}}\,\partial_{_{\vphantom{\lr{(2j)}}}}^\mu\Phi^{^{\lr{(2j)}}}_{\mu\nu_1\cdots\nu_{2j-1}}
					\end{split}
				\end{equation*}
				will not modify the equations of motion since it is just $0$ rewritten in a convenient way.
				Secondly, we have
				\begin{equation*}
					0=-\sum_{j=0}^{\mathclap{\leq\nicefrac{D}{2}-1}}\frac{1}{\lr{(2j)}!}\,\Phi^{^{\lr{(2j)}}}_{\nu_1\cdots\nu_{2j}}\,\partial^{^{\vphantom{\lr{(2j)}}}}_\mu\partial^{^{\vphantom{\lr{(2j)}}}}_\nu\Phi_{_{\lr{(2j+2)}}}^{\mu\nu\nu_1\cdots\nu_{2j}}\eqpc
				\end{equation*}
				since it contains the contraction of a symmetric tensor ($\partial_\mu\partial_\nu$) with an anti-symmetric one ($\Phi^{\mu\nu\cdots}$).
				Again, using an integration by parts (and forgetting the boundary terms), we define
				\begin{equation*}
					\mathcal{L}_-=\sum_{j=0}^{\mathclap{\leq\nicefrac{D}{2}-1}}\frac{1}{\lr{(2j)}!}\,\partial^{^{\vphantom{\lr{(2j)}}}}_\mu\Phi_{_{\lr{(2j+2)}}}^{\mu\nu\nu_1\cdots\nu_{2j}}\,\partial^{^{\vphantom{\lr{(2j)}}}}_\nu\Phi^{^{\lr{(2j)}}}_{\nu_1\cdots\nu_{2j}}
				\end{equation*}
				and this term does not affect the dynamics either.
				
				We now consider the theory described by the Lagrangian density
				\begin{equation*}
					\mathcal{L}=\mathcal{L}_0+\mathcal{L}_++\mathcal{L}_-\eqpc
				\end{equation*}
				which should be equivalent to the simultaneous treatment of $2^d$ independent problems (that belong in the class of the non-linear wave equation).
				Let us now construct the DW Hamiltonian formulation of this theory.
				We first introduce, for each dynamic field in the collection (\ref{preliminaries.dwh.nodegfields.collection}) and for each direction of space-time, a conjugate momentum
				\begin{align*}
					\Psi_{_{\lr{(2i+1)}}}^{\mu_0\cdots\mu_{2i}}&=\Psi_{_{\lr{(2i+1)}}}^{\mu_0\lr{[\mu_1\cdots\mu_{2i}]}}=\frac{\partial\mathcal{L}}{\partial\lr{(\partial^{^{\vphantom{\lr{(2i)}}}}_{\mu_0}\Phi^{^{\lr{(2i)}}}_{\mu_1\cdots\mu_{2i}})}}\\
					\begin{split}
						&=\partial_{_{\vphantom{\lr{(2i)}}}}^{\mu_0}\Phi_{_{\lr{(2i)}}}^{\mu_1\cdots\mu_{2i}}+\partial^{^{\vphantom{\lr{(2i+2)}}}}_\mu\Phi_{_{\lr{(2i+2)}}}^{\mu\mu_0\cdots\mu_{2i}}\\
						&\quad+2i\,\eta_{_{\vphantom{\lr{(2i)}}}}^{\mu_0\lr{[\mu_1\vphantom{\lr{|\mu|}\mu_2\cdots\mu_{2i}}.}}\partial^{^{\vphantom{\lr{(2i)}}}}_\mu\Phi_{_{\lr{(2i)}}}^{\lr{.\vphantom{\mu_1}\lr{|\mu|}\mu_2\cdots\mu_{2i}]}}\\
						&\quad-2i\,\partial_{_{\vphantom{\lr{(2i)}}}}^{\lr{[\mu_1\vphantom{\lr{|\mu_0|}\mu_2\cdots\mu_{2i}}.}}\Phi_{_{\lr{(2i)}}}^{\lr{.\vphantom{\mu_1}\lr{|\mu_0|}\mu_2\cdots\mu_{2i}]}}\\
						&\quad+2i\lr{(2i-1)}\,\eta_{_{\vphantom{\lr{(2i-2)}}}}^{\mu_0\lr{[\mu_1\vphantom{\mu_3\cdots\mu_{2i}}.}}\partial_{_{\vphantom{\lr{(2i-2)}}}}^{\mu_2}\Phi_{_{\lr{(2i-2)}}}^{\lr{.\vphantom{\mu_1}\mu_3\cdots\mu_{2i}]}}\eqpc
					\end{split}
				\end{align*}
				where we have used that
				\begin{equation*}
					\frac{\partial\big(\partial^{^{\vphantom{\lr{(j)}}}}_{\nu_0}\Phi^{^{\lr{(j)}}}_{\nu_1\cdots\nu_{j}}\big)}{\partial\lr{(\partial^{^{\vphantom{\lr{(i)}}}}_{\mu_0}\Phi^{^{\lr{(i)}}}_{\mu_1\cdots\mu_{i}})}}=j!\,\delta_i^j\,\delta_{\nu_0}^{\mu_0}\,\delta_{\lr{[\nu_1\vphantom{\cdots\nu_j}.}}^{\mu_1}\!\cdots\delta_{\lr{.\vphantom{\nu_1\cdots}\nu_j]}}^{\mu_i}\eqpd
				\end{equation*}
				
				Then, defining
				\begin{align*}
					\Phi_{_{\lr{(2i+1)}}}^{\mu_0\cdots\mu_{2i}}&=\Psi_{_{\lr{(2i+1)}}}^{\lr{[\mu_0\cdots\mu_{2i}]}}\\
					&=\lr{(2i+1)}\,\partial_{_{\vphantom{\lr{(2i)}}}}^{\lr{[\mu_0\vphantom{\cdots\mu_{2i}}.}}\Phi_{_{\lr{(2i)}}}^{\lr{.\vphantom{\mu_0}\mu_1\cdots\mu_{2i}]}}+\partial^{^{\vphantom{\lr{(2i)}}}}_\mu\Phi_{_{\lr{(2i+2)}}}^{\mu\mu_0\cdots\mu_{2i}}\eqpc
				\end{align*}
				we can prove\footnote{By substituting the expression of $\Phi_{_{\lr{(2i+1)}}}^{\mu_0\cdots\mu_{2i}}$ and $\Phi_{_{\lr{(2i-1)}}}^{\mu_2\cdots\mu_{2i}}$ and after a straightforward identification it remains to prove that
					\begin{equation*}
						\lr{(2i+1)}\,\partial_{_{\vphantom{\lr{(2i)}}}}^{\lr{[\mu_0\vphantom{\cdots\mu_{2i}}.}}\Phi_{_{\lr{(2i)}}}^{\lr{.\vphantom{\mu_0}\mu_1\cdots\mu_{2i}]}}=\partial_{_{\vphantom{\lr{(2i)}}}}^{\mu_0}\Phi_{_{\lr{(2i)}}}^{\mu_1\cdots\mu_{2i}}-2i\,\partial_{_{\vphantom{\lr{(2i)}}}}^{\lr{[\mu_1\vphantom{\lr{|\mu_0|}\mu_2\cdots\mu_{2i}}.}}\Phi_{_{\lr{(2i)}}}^{\lr{.\vphantom{\mu_1}\lr{|\mu_0|}\mu_2\cdots\mu_{2i}]}}\eqpd
					\end{equation*}
					On the other hand, we successively have
					\begin{align*}
						\lr{(2i+1)}&\,\partial_{_{\vphantom{\lr{(2i)}}}}^{\lr{[\mu_0\vphantom{\cdots\mu_{2i}}.}}\Phi_{_{\lr{(2i)}}}^{\lr{.\vphantom{\mu_0}\mu_1\cdots\mu_{2i}]}}=\nonumber\\*
						&=\lr{(2i+1)}\,\frac{1}{2i+1}\sum_{j=1}^{2i+1}\partial_{_{\vphantom{\lr{(2i)}}}}^{\lr{[\mu_j\vphantom{\mu_{j+1}\cdots\mu_{2i}\lr{|\mu_0|}\mu_1\cdots\mu_{j-1}}.}}\Phi_{_{\lr{(2i)}}}^{\lr{.\vphantom{\mu_j}\mu_{j+1}\cdots\mu_{2i}\lr{|\mu_0|}\mu_1\cdots\mu_{j-1}]}}\\
						&=\partial_{_{\vphantom{\lr{(2i)}}}}^{\mu_0}\Phi_{_{\lr{(2i)}}}^{\mu_1\cdots\mu_{2i}}+\sum_{j=1}^{2i}\partial_{_{\vphantom{\lr{(2i)}}}}^{\lr{[\mu_j\vphantom{\mu_{j+1}\cdots\mu_{2i}\lr{|\mu_0|}\mu_1\cdots\mu_{j-1}}.}}\Phi_{_{\lr{(2i)}}}^{\lr{.\vphantom{\mu_j}\mu_{j+1}\cdots\mu_{2i}\lr{|\mu_0|}\mu_1\cdots\mu_{j-1}]}}\\
						&=\partial_{_{\vphantom{\lr{(2i)}}}}^{\mu_0}\Phi_{_{\lr{(2i)}}}^{\mu_1\cdots\mu_{2i}}+\sum_{j=1}^{2i}\partial_{_{\vphantom{\lr{(2i)}}}}^{\lr{[\mu_j\vphantom{\lr{|\mu_0|}\mu_1\cdots\mu_{j-1}\mu_{j+1}\cdots\mu_{2i}}.}}\Phi_{_{\lr{(2i)}}}^{\lr{.\vphantom{\mu_j}\lr{|\mu_0|}\mu_1\cdots\mu_{j-1}\mu_{j+1}\cdots\mu_{2i}]}}\\
						&=\partial_{_{\vphantom{\lr{(2i)}}}}^{\mu_0}\Phi_{_{\lr{(2i)}}}^{\mu_1\cdots\mu_{2i}}-\sum_{j=1}^{2i}\partial_{_{\vphantom{\lr{(2i)}}}}^{\lr{[\mu_1\vphantom{\lr{|\mu_0|}\mu_j\mu_2\cdots\mu_{j-1}\mu_{j+1}\cdots\mu_{2i}}.}}\Phi_{_{\lr{(2i)}}}^{\lr{.\vphantom{\mu_1}\lr{|\mu_0|}\mu_j\mu_2\cdots\mu_{j-1}\mu_{j+1}\cdots\mu_{2i}]}}\\
						\begin{split}
							&=\partial_{_{\vphantom{\lr{(2i)}}}}^{\mu_0}\Phi_{_{\lr{(2i)}}}^{\mu_1\cdots\mu_{2i}}\\
							&-\sum_{j=1}^{2i}{\lr{(-1)}}^{j-2}\,\partial_{_{\vphantom{\lr{(2i)}}}}^{\lr{[\mu_1\vphantom{\lr{|\mu_0|}\mu_2\cdots\mu_{j-1}\mu_j\mu_{j+1}\cdots\mu_{2i}}.}}\!{\lr{(-1)}}^{j-2}\,\Phi_{_{\lr{(2i)}}}^{\lr{.\vphantom{\mu_1}\lr{|\mu_0|}\mu_2\cdots\mu_{j-1}\mu_j\mu_{j+1}\cdots\mu_{2i}]}}
						\end{split}\\
						&=\partial_{_{\vphantom{\lr{(2i)}}}}^{\mu_0}\Phi_{_{\lr{(2i)}}}^{\mu_1\cdots\mu_{2i}}-\sum_{j=1}^{2i}\partial_{_{\vphantom{\lr{(2i)}}}}^{\lr{[\mu_1\vphantom{\lr{|\mu_0|}\mu_2\cdots\mu_{2i}}.}}\Phi_{_{\lr{(2i)}}}^{\lr{.\vphantom{\mu_1}\lr{|\mu_0|}\mu_2\cdots\mu_{2i}]}}\\
						&=\partial_{_{\vphantom{\lr{(2i)}}}}^{\mu_0}\Phi_{_{\lr{(2i)}}}^{\mu_1\cdots\mu_{2i}}-2i\,\partial_{_{\vphantom{\lr{(2i)}}}}^{\lr{[\mu_1\vphantom{\lr{|\mu_0|}\mu_2\cdots\mu_{2i}}.}}\Phi_{_{\lr{(2i)}}}^{\lr{.\vphantom{\mu_1}\lr{|\mu_0|}\mu_2\cdots\mu_{2i}]}}\eqpd
					\end{align*}
					\hfill$\square$} that
				\begin{equation*}
					\Psi_{_{\lr{(2i+1)}}}^{\mu_0\lr{[\mu_1\cdots\mu_{2i}]}}=\Phi_{_{\lr{(2i+1)}}}^{\mu_0\cdots\mu_{2i}}+2i\,\eta_{_{\vphantom{\lr{(2i-1)}}}}^{\mu_0\lr{[\mu_1\vphantom{\cdots\mu_{2i}}.}}\,\Phi_{_{\lr{(2i-1)}}}^{\lr{.\vphantom{\mu_1}\mu_2\cdots\mu_{2i}]}}\eqpc
				\end{equation*}
				and in consequence, that all the conjugate momenta are indeed already defined in the collection (\ref{preliminaries.dwh.nodegfields.collection}).
				
				The \textsc{Euler}~-- \textsc{Lagrange} equation reads
				\begin{align*}
					\frac{\partial\mathcal{L}}{\partial\Phi^{^{\lr{(2i)}}}_{\mu_1\cdots\mu_{2i}}}&=\partial^{^{\vphantom{\lr{(2i)}}}}_{\mu_0}\frac{\partial\mathcal{L}}{\partial\lr{(\partial^{^{\vphantom{\lr{(2i)}}}}_{\mu_0}\Phi^{^{\lr{(2i)}}}_{\mu_1\cdots\mu_{2i}})}}=\partial^{^{\vphantom{\lr{(2i)}}}}_{\mu_0}\Psi_{_{\lr{(2i+1)}}}^{\mu_0\lr{[\mu_1\cdots\mu_{2i}]}}\\
					\begin{split}
						&=\partial^{^{\vphantom{\lr{(2i)}}}}_\mu\partial_{_{\vphantom{\lr{(2i)}}}}^\mu\Phi_{_{\lr{(2i)}}}^{\mu_1\cdots\mu_{2i}}+\partial^{^{\vphantom{\lr{(2i+2)}}}}_\mu\partial^{^{\vphantom{\lr{(2i+2)}}}}_\nu\Phi_{_{\lr{(2i+2)}}}^{\nu\mu\mu_1\cdots\mu_{2i}}\\
						&\quad+2i\,\partial_{_{\vphantom{\lr{(2i)}}}}^{\lr{[\mu_1\vphantom{\lr{|\mu|}\mu_2\cdots\mu_{2i}}.}}\partial^{^{\vphantom{\lr{(2i)}}}}_\mu\Phi_{_{\lr{(2i)}}}^{\lr{.\vphantom{\mu_1}\lr{|\mu|}\mu_2\cdots\mu_{2i}]}}\\
						&\quad-2i\,\partial_{_{\vphantom{\lr{(2i)}}}}^{\lr{[\mu_1\vphantom{\lr{|\mu|}\mu_2\cdots\mu_{2i}}.}}\partial^{^{\vphantom{\lr{(2i)}}}}_\mu\Phi_{_{\lr{(2i)}}}^{\lr{.\vphantom{\mu_1}\lr{|\mu|}\mu_2\cdots\mu_{2i}]}}\\
						&\quad+2i\lr{(2i-1)}\,\partial_{_{\vphantom{\lr{(2i-2)}}}}^{\lr{[\mu_1\vphantom{\mu_3\cdots\mu_{2i}}.}}\partial_{_{\vphantom{\lr{(2i-2)}}}}^{\mu_2}\Phi_{_{\lr{(2i-2)}}}^{\lr{.\vphantom{\mu_1}\mu_3\cdots\mu_{2i}]}}
					\end{split}\\
					&=\partial^{^{\vphantom{\lr{(2i)}}}}_\mu\partial_{_{\vphantom{\lr{(2i)}}}}^\mu\Phi_{_{\lr{(2i)}}}^{\mu_1\cdots\mu_{2i}}=\square\Phi_{_{\lr{(2i)}}}^{\mu_1\cdots\mu_{2i}}\eqpd
				\end{align*}
				
				Summarizing, we first consider the concatenation of $2^d$ independent theories (through the Lagrangian density $\mathcal{L}_0$).
				Based on what we introduced in the previous sections, we know that the DW Hamiltonian formulation of this theory requires the collection $\{\Phi^{^{\lr{(2i)}}},\Psi_{_{\lr{(2i+1)}}}\}$ of $2^d\lr{(D+1)}$ fields to construct the phase space ($2^d$ dynamical fields and $D\,2^d$ conjugate momenta).
				We also know that this phase space is not a multi-symplectic manifold since the multi-symplectic structure of each ``sub''-theory is, independently, degenerate.
				As we claimed earlier, this comes from the existence of a non-canonical conjugate field.
				To fix this issue, we introduce the additional coupling terms $\mathcal{L}_+$ and $\mathcal{L}_-$.
				Above all, these couplings do not affect the dynamics and each ``sub''-theory remains independent.
				However, these couplings have an interesting side effect: they allow a conjugate momentum to be shared by several dynamical fields and thus enable all the conjugate fields to be canonical.
				Thus, correctly chosen, $\mathcal{L}_+$ and $\mathcal{L}_-$ lead to the closed collection $\{\Phi^{^{\lr{(2i)}}},\Phi_{_{\lr{(2i+1)}}}\}=\{\Phi^{^{\lr{(i)}}}\}$ of $2^D$ fields, equality composed of dynamic fields and canonical conjugate momenta, and where ``closed'' has two meanings.
				On the one hand, all the dynamic fields have all their conjugate momenta in the collection.
				On the other hand, all the conjugate fields are, for all the directions of space-time, the conjugate momentum of a dynamic field that belong in the collection.
				
				Now, from the \textsc{Euler}~-- \textsc{Lagrange} equations and the definitions of the conjugate momenta, we get the \textsc{De Donder}~-- \textsc{Weyl}~-- \textsc{Hamilton} equations
				\begin{subequations}
					\begin{align}
						-\partial^{^{\vphantom{\lr{(2i+1)}}}}_\mu\Phi_{_{\lr{(2i+1)}}}^{\mu\mu_1\cdots\mu_{2i}}-2i\,\partial_{_{\vphantom{\lr{(2i-1)}}}}^{\lr{[\mu_1\vphantom{\mu_2\cdots\mu_{2i}}.}}\Phi_{_{\lr{(2i-1)}}}^{\lr{.\vphantom{\mu_1}\mu_2\cdots\mu_{2i}]}}&=\frac{\partial\mathcal{H}}{\partial\Phi^{^{\lr{(2i)}}}_{\mu_1\cdots\mu_{2i}}}\label{preliminaries.dwh.hamilteq.nodeg.moments}\\
						\lr{(2i+1)}\,\partial^{^{\vphantom{\lr{(2i)}}}}_{\lr{[\mu_0\vphantom{\mu_1\cdots\mu_{2i}}.}}\Phi^{^{\lr{(2i)}}}_{\lr{.\vphantom{\mu_0}\mu_1\cdots\mu_{2i}]}}+\partial_{_{\vphantom{\lr{(2i+2)}}}}^\mu\Phi^{^{\lr{(2i+2)}}}_{\mu\mu_0\cdots\mu_{2i}}&=\frac{\partial\mathcal{H}}{\partial\Phi_{_{\lr{(2i+1)}}}^{\mu_0\cdots\mu_{2i}}}\label{preliminaries.dwh.hamilteq.nodeg.field}
					\end{align}
				\end{subequations}
				where
				\begin{equation*}
					\begin{split}
						\mathcal{H}&=\frac{1}{2}\sum_{j=0}^{\leq\mathrlap{\nicefrac{\lr{(D-1)}}{2}}\phantom{\nicefrac{D}{2}}}\frac{1}{\lr{(2j+1)}!}\,\Phi_{_{\lr{(2j+1)}}}^{\nu_0\cdots\nu_{2j}}\,\Phi^{^{\lr{(2j+1)}}}_{\nu_0\cdots\nu_{2j}}\\
						&\quad+\sum_{j=0}^{\leq\nicefrac{D}{2}}\frac{1}{\lr{(2j)}!}\,V_{_{\lr{(2j)}}}^{\lr{[\nu_1\cdots\nu_{2j}]}}\lr{(\Phi^{^{\lr{(2j)}}}_{\nu_1\cdots\nu_{2j}})}\eqpc
					\end{split}
				\end{equation*}
				is the Hamiltonian density associated to $\mathcal{L}$.
				
				In order to identify the form of \cref{preliminaries.dwh.hamilteq.abstractform}, we need to flatten all these indices (\emph{i.e.} for every configuration of values of all these indices we associate one, and only one, index).
				So, if the list $\lr{(\mu_1\cdots\mu_i)}$ is sorted and free of duplicates, we define
				\begin{align*}
					\begin{split}
						\text{fl}\lr{(i,\mu_1\cdots\mu_i)}&=1+\sum_{j=0}^{i-1}\binom{D}{j}\\
						&\quad+\sum_{j=1}^i\sum_{\substack{\nu_j=\\\mu_{j-1}+1}}^{\mu_j-1}\sum_{\substack{\nu_{j+1}=\\\nu_j+1}}^{d-i+j+1}\cdots\sum_{\substack{\nu_i=\\\nu_{i-1}+1}}^{d}1
					\end{split}\\
					&\in\lr{\llbracket1,2^D\rrbracket}\eqpc
				\end{align*}
				with $\mu_0=-1$ and where the stacking of sums can be re-expressed in term of generalized harmonic numbers as well.
				
				In the following, we assume that both $\lr{(\alpha_1\cdots\alpha_i)}$ and $\lr{(\beta_1\cdots\beta_j)}$ are sorted and duplicate free.
				Thus, defining the vector state $\bm{\zeta}$ as
				\begin{align*}
					\zeta^a&=\zeta^{a=\text{fl}\lr{(2i,\alpha_1\cdots\alpha_{2i})}}=\Phi^{^{\lr{(2i)}}}_{\alpha_1\cdots\alpha_{2i}}\\
					&=\zeta^{a=\text{fl}\lr{(2i+1,\alpha_1\cdots\alpha_{2i+1})}}=\Phi_{_{\lr{(2i+1)}}}^{\alpha_1\cdots\alpha_{2i+1}}\eqpc
				\end{align*}
				\cref{preliminaries.dwh.hamilteq.nodeg.moments,preliminaries.dwh.hamilteq.nodeg.field} can be rewritten in the form of \cref{preliminaries.dwh.hamilteq.abstractform}, \emph{i.e.}
				\begin{equation*}
					{M^{\mu\,}}_{\phantom{a}b}^a\,\partial_\mu\zeta^b=\bm{\partial}^a\mathcal{H}\eqpc
				\end{equation*}
				provided that
				\begin{equation*}
					\begin{split}
						{M^{\mu\,}}_{\vphantom{b=\text{fl}\lr{(j,\beta_1\cdots\beta_j)}}}^{a=\text{fl}\lr{(i,\alpha_1\cdots\alpha_i)}}&{\vphantom{\bm{M}^\mu}}^{\vphantom{a=\text{fl}\lr{(i,\alpha_1\cdots\alpha_i)}}}_{b=\text{fl}\lr{(j,\beta_1\cdots\beta_j)}}=\\
						\lr{(i\bmod 2-1)}&\times\\
						\Bigg[\delta_j^{i+1}\sum_{k=1}^j&{\lr{(-1)}}^{k-1}\delta_{\vphantom{\beta_j}\beta_1}^{\vphantom{\alpha_{k+1}}\alpha_1}\cdots\delta_{\vphantom{\beta_j}\beta_{k-1}}^{\alpha_{k-1}}\delta_{\vphantom{\beta_j}\beta_k}^{\vphantom{\alpha_{k+1}}\mu}\delta_{\vphantom{\beta_j}\beta_{k+1}}^{\vphantom{\alpha_{k+1}}\alpha_k}\cdots\delta_{\beta_j}^{\vphantom{\alpha_{k+1}}\alpha_i}\\
						+\,\delta_j^{i-1}\sum_{k=1}^i&{\lr{(-1)}}^{k-1}\delta_{\vphantom{\beta_j}\beta_1}^{\vphantom{\alpha_{k+1}}\alpha_1}\cdots\delta_{\vphantom{\beta_j}\beta_{k-1}}^{\alpha_{k-1}}\eta_{\vphantom{\beta_j}}^{\vphantom{\alpha_{k+1}}\mu\alpha_k}\delta_{\vphantom{\beta_j}\beta_k}^{\alpha_{k+1}}\cdots\delta_{\beta_j}^{\vphantom{\alpha_{k+1}}\alpha_i}\Bigg]\hphantom{\eqpd}\\
						+\,\lr{(i\bmod 2)}&\times\\
						\Bigg[\delta_i^{j+1}\sum_{k=1}^i&{\lr{(-1)}}^{k-1}\delta_{\vphantom{\alpha_{k+1}}\alpha_1}^{\vphantom{\beta_j}\beta_1}\cdots\delta_{\alpha_{k-1}}^{\vphantom{\beta_j}\beta_{k-1}}\delta_{\vphantom{\alpha_{k+1}}\alpha_k}^{\vphantom{\beta_j}\mu}\delta_{\alpha_{k+1}}^{\vphantom{\beta_j}\beta_k}\cdots\delta_{\vphantom{\alpha_{k+1}}\alpha_i}^{\beta_j}\\
						+\,\delta_i^{j-1}\sum_{k=1}^j&{\lr{(-1)}}^{k-1}\delta_{\vphantom{\alpha_{k+1}}\alpha_1}^{\vphantom{\beta_j}\beta_1}\cdots\delta_{\alpha_{k-1}}^{\vphantom{\beta_j}\beta_{k-1}}\eta_{\vphantom{\alpha_{k+1}}}^{\vphantom{\beta_j}\mu\beta_k}\delta_{\vphantom{\alpha_{k+1}}\alpha_k}^{\vphantom{\beta_j}\beta_{k+1}}\cdots\delta_{\vphantom{\alpha_{k+1}}\alpha_i}^{\beta_j}\Bigg]\eqpd
					\end{split}
				\end{equation*}
				These $\lr{\{\bm{M}^\mu\}}$ matrices are skew-symmetric, linearly independent and non-degenerate (they all have two eigenvalues, $\pm i$, $2^d$ times degenerate).
				The phase space of the theory, $\Omega$, is now a multi-symplectic manifold and we have finally obtained a correct covariant Hamiltonian formulation of the non-linear wave equation.
				
				Let us illustrate how this construction works in the particular dimension $D=1+1$.
				We start by considering the collection (\ref{preliminaries.dwh.nodegfields.collection}) and for notational convenience we create aliases for these $2^2$ fields as
				\begin{subequations}
					\begin{align}
						\Phi^{^{\lr{(0)}}}&=\phi\label{preliminaries.1dfields.def.phi}\\
						\Psi_{_{\lr{(1)}}}^0=\Phi_{_{\lr{(1)}}}^0&=\psi^0\label{preliminaries.1dfields.def.psi0}\\
						\Psi_{_{\lr{(1)}}}^1=\Phi_{_{\lr{(1)}}}^1&=\psi^1\label{preliminaries.1dfields.def.psi1}\\
						\Phi^{^{\lr{(2)}}}_{01}=-\Phi^{^{\lr{(2)}}}_{10}&=\gamma\label{preliminaries.1dfields.def.gamma}\\
						\Psi_{_{\lr{(3)}}}^{001}=-\Psi_{_{\lr{(3)}}}^{010}=\Phi_{_{\lr{(1)}}}^1&=\psi^1\nonumber\\
						\Psi_{_{\lr{(3)}}}^{101}=-\Psi_{_{\lr{(3)}}}^{110}=\Phi_{_{\lr{(1)}}}^0&=\psi^0\nonumber\eqpd
					\end{align}
				\end{subequations}
				Here, we have introduced one new dynamic field, $\gamma$, independent of $\phi$.
				See \cref{preliminaries.dwh.nodegfields.graphic} for a graphical representation of these fields.
				Then we consider the concatenation of these two theories described by the Lagrangian density
				\begin{equation*}
					\begin{split}
						\mathcal{L}&=\frac{1}{2}{\lr{(\partial_0\phi)}}^2-\frac{1}{2}{\lr{(\partial_1\phi)}}^2-\frac{1}{2}{\lr{(\partial_0\gamma)}}^2+\frac{1}{2}{\lr{(\partial_1\gamma)}}^2\\
						&\quad+\partial_0\phi\,\partial_1\gamma-\partial_0\gamma\,\partial_1\phi-V\lr{(\phi)}-\tilde{V}\lr{(\gamma)}\eqpc
					\end{split}
				\end{equation*}
				where $\tilde{V}$ is the potential of the extra theory which can be freely chosen.
				The (canonical) conjugate momenta are defined by
					\begin{align*}
						\psi^0&=\frac{\partial\mathcal{L}}{\partial\lr{(\partial_0\phi)}}=\frac{\partial\mathcal{L}}{\partial\lr{(\partial_1\gamma)}}=\partial_0\phi+\partial_1\gamma \eqpc
						\\
						\psi^1&=\frac{\partial\mathcal{L}}{\partial\lr{(\partial_1\phi)}}=\frac{\partial\mathcal{L}}{\partial\lr{(\partial_0\gamma)}}=-\partial_0\gamma-\partial_1\phi\eqpd
					\end{align*}
				The Hamiltonian density is
				\begin{equation*}
					\mathcal{H}=\frac{1}{2}{\psi^0}^2-\frac{1}{2}{\psi^1}^2+V\lr{(\phi)}+\tilde{V}\lr{(\gamma)}\eqpc
				\end{equation*}
				and defining the vector state
				\begin{equation*}
					\bm{\zeta}^{\text{T}}=\begin{bmatrix}
						\phi&\psi^0&\psi^1&\gamma
					\end{bmatrix}\eqpc
				\end{equation*}
				the dynamics are fully described by \cref{preliminaries.dwh.hamilteq.abstractform} provided that
					\begin{align*}
						\bm{M}^0&=\begin{bmatrix}
							0&-1&0&0\\
							1&0&0&0\\
							0&0&0&1\\
							0&0&-1&0
						\end{bmatrix}\eqpc\\
						\bm{M}^1&=\begin{bmatrix}
							0&0&-1&0\\
							0&0&0&1\\
							1&0&0&0\\
							0&-1&0&0
						\end{bmatrix}\eqpd
					\end{align*}
				These two matrices are no longer degenerate (their eigenvalues are $+i$, $+i$, $-i$ and $-i$) and they define a multi-symplectic structure on the phase space.
				
				Recaping, we started with two theories (one for $\phi$ and one for $\gamma$). While isolated they break the multi-symplectic structure of phase space.
				Joined together they complete each other such that the theory of the coalition restores multi-symplecticity.
				
				Using this construction, we reduced the number of unknowns for each independent ``sub''-theory from $3$: $\lr{\{\phi,\psi^0,\psi^1\}}$ and $\lr{\{\gamma,\psi^1,\psi^0\}}$ to $2$: $\lr{\{\phi,\psi^0\}}$ and $\lr{\{\gamma,\psi^1\}}$.
				In a general settup we reduced the number of unknowns per independent ``sub''-theory from $D+1$ (the field and all its conjugate momentums) to $2$ (the field and its canonical conjugate momentum) as in the standard non-covariant Hamiltonian formulation of field theory.
				
				Finally, we want to stress that from the numerical point of view, to integrate a $D$-dimensional theory, we actually need to integrate $2^d$ $D$-dimensional theories.
				These ``extra'' theories can be used in two ways (or a mix of the two):
				\begin{enumerate}[label = \emph{\roman{*}}., labelindent = 0em, leftmargin = *, widest* = 2, nosep]
					\item By considering a theory the solution of which is known we get an error estimate of the integration process.
					This feature comes from the fact that the integration is performed through the conjugate momentum which is itself shared between different dynamic fields.
					If an error occurs during the integration of one of the fields, it will reverberate on the others and will be caught by the control field(s).
					\item They can be used to integrate, at the same time, several replicas of the theory (in a statistical approach for example) or even different theories.
				\end{enumerate}
				
				In the present section we have introduced a construction that leads to a phase space that is a multi-symplectic manifold.
				In the next sections, we will first prove the conservation of the multi-symplectic structure under the Hamiltonian flow.
				Then we will define the stress-energy tensor, the charges and discuss their properties.
			\subsubsection{Conservation of the multi-symplectic structure}
				To prove the conservation of the multi-symplectic structure, we consider the dual of \cref{preliminaries.dwh.hamilteq.abstractform} which reads
				\begin{equation}
					\begin{split}
						{M^\mu}_{ab}\,\partial_\mu\zeta^b\,\db[a]{}&=\bm{\partial}_a\mathcal{H}\,\db[a]{}\\
						\text{\emph{i.e.}}\quad\bm{\omega}^\mu\lr{(\partial_\mu\bm{\zeta},\cdot)}&=\d{\mathcal{H}}\eqpd
					\end{split}\label{preliminaries.dwh.formhamilteq.abstractform}
				\end{equation}
				Now, taking the exterior derivative of it, we successively have
				\begin{align*}
					\d{\lr{(\bm{\omega}^\mu\lr{(\partial_\mu\bm{\zeta},\cdot)})}}&=\d{\d{\mathcal{H}}}=0\\
					&=\d{\lr{({M^\mu}_{ab}\,\partial_\mu\zeta^b\,\db[a]{})}}\\
					&={M^\mu}_{ab}\,\d{\lr{(\partial_\mu\zeta^b)}}\wedge\db[a]{}\\
					&={M^\mu}_{ab}\,\lr{(\partial_\mu\db[b]{})}\wedge\db[a]{}\eqpc
				\end{align*}
				where we have used that $\d{}$ and $\partial_\mu$ commute since they act in different spaces ($\Omega$ does not depend on the position on $\mathcal{M}$).
				
				On the other hand, we have
				\begin{align*}
					\partial_\mu\bm{\omega}^\mu&=\frac{1}{2}{M^\mu}_{ab}\,\lr{(-\lr{(\partial_\mu\db[\vphantom{b}a]{})}\wedge\db[b]{}-\db[a]{}\wedge\lr{(\partial_\mu\db[b]{})})}\\
					&=\frac{1}{2}{M^\mu}_{ab}\,\lr{(-\lr{(\partial_\mu\db[\vphantom{b}a]{})}\wedge\db[b]{}+\lr{(\partial_\mu\db[b]{})}\wedge\db[a]{})}\\
					&={M^\mu}_{ab}\,\lr{(\partial_\mu\db[b]{})}\wedge\db[a]{}=0\eqpd
				\end{align*}
				Thus, we have proved (on-shell) the local conservation of the multi-symplectic structure
				\begin{equation}
					\partial_\mu\bm{\omega}^\mu=0\eqpd\label{preliminaries.dwh.locmultisympconserv}
				\end{equation}
				By definition a multi-symplectic integrator is a numerical method that exactly preserves the discrete version of \cref{preliminaries.dwh.locmultisympconserv}.
			\subsubsection{The stress-energy tensor, its conservation and the charges}\label{preliminaries.dwh.stresstensorsection}
				We define now the stress-energy tensor as the symmetric $2$-tensor
				\begin{align}
						\mathcal{T}^{\mu\nu}&=\frac{1}{2}\lr{(\bm{\omega}^\mu\lr{(\partial^\nu\bm{\zeta},\bm{\zeta})}+\bm{\omega}^\nu\lr{(\partial^\mu\bm{\zeta},\bm{\zeta})}+\eta^{\mu\nu}\bm{\omega}^\kappa\lr{(\bm{\zeta},\partial_\kappa\bm{\zeta})})}\nonumber\\*
						&\quad+\eta^{\mu\nu}\mathcal{H}\label{preliminaries.dwh.stresstensor}\\
					&=\mathcal{T}^{\nu\mu}\eqpd\nonumber
				\end{align}
				Before proving that it is subject to a local conservation law, we need to make some remarks.
				First of all, performing an integration by parts and provided that boundary terms vanish (space-time is unbounded, subject to periodic boundary conditions, \dots) we have
				\begin{equation}
					\int\d[D]{x}\,\bm{\omega}^\mu\lr{(\partial_\nu\partial^\nu\bm{\zeta},\bm{\zeta})}=-\int\d[D]{x}\,\bm{\omega}^\mu\lr{(\partial_\nu\bm{\zeta},\partial^\nu\bm{\zeta})}=0\eqpd\label{preliminaries.dwh.firstrequirement.locstressconserv}
				\end{equation}
				Since $\bm{\omega}^\mu$ is skew-symmetric in its two arguments $\bm{\omega}^\mu\lr{(\partial_\nu\bm{\zeta},\partial^\nu\bm{\zeta})}$ identically vanishes.
				We stress that \cref{preliminaries.dwh.firstrequirement.locstressconserv} holds for any $\bm{\zeta}$ (\emph{i.e.} off-shell), so the integrands are equal, and hence
				\begin{equation*}
					\bm{\omega}^\mu\lr{(\partial_\nu\partial^\nu\bm{\zeta},\bm{\zeta})}=-\bm{\omega}^\mu\lr{(\partial_\nu\bm{\zeta},\partial^\nu\bm{\zeta})}=0\eqpd
				\end{equation*}
				Secondly, following the same reasoning we have
				\begin{align*}
					\bm{\omega}^\nu\lr{(\bm{\zeta},\partial_\nu\partial^\mu\bm{\zeta})}&=-\bm{\omega}^\nu\lr{(\partial^\mu\bm{\zeta},\partial_\nu\bm{\zeta})}\\
					&=-\bm{\omega}^\nu\lr{(\partial_\nu\bm{\zeta},\partial^\mu\bm{\zeta})}=\bm{\omega}^\nu\lr{(\partial^\mu\bm{\zeta},\partial_\nu\bm{\zeta})}\\
					&=0\eqpd
				\end{align*}
				We now have all the necessary ingredients to prove the local conservation of the stress-energy tensor.
				Using the relations established above, and \cref{preliminaries.dwh.formhamilteq.abstractform,preliminaries.dwh.locmultisympconserv}, it is straightforward to deduce
				\begin{equation}
					\partial_\nu\mathcal{T}^{\mu\nu}=\partial_\nu\mathcal{T}^{\nu\mu}=0\eqpd\label{preliminaries.dwh.locstressconserv}
				\end{equation}
				Hence, the stress-energy tensor is locally conserved (on-shell).
				
				Let us now sketch why we have defined the stress-energy tensor as \cref{preliminaries.dwh.stresstensor}.
				We suppose that the action of our theory,
				\begin{equation*}
					\mathcal{S}=\int\d[D]{x}\,\mathcal{L}\lr{(\lr{\{\phi^i,\partial_\mu\phi^i\}})}\eqpc
				\end{equation*}
				is invariant under the local transformation
				\begin{equation}
					x^\mu\to x'^\mu=x^\mu+\epsilon\,\ell^\mu\eqpc\label{preliminaries.dwh.localtransform.stresstensor}
				\end{equation}
				where $\epsilon\ll1$ is constant while $\ell$ may depend on $x$.
				The local variation of the field is
				\begin{equation*}
					\epsilon\,\delta\phi^i\lr{(x)}=\phi'^i\lr{(x')}-\phi^i\lr{(x)}\eqpc
				\end{equation*}
				while the global variation is
				\begin{equation*}
					\epsilon\,\Delta\phi^i\lr{(x)}=\phi'^i\lr{(x)}-\phi^i\lr{(x)}\eqpd
				\end{equation*}
				The local variation only depends on the nature of the field (scalar, vector, tensor, \dots) while the global variation includes the effect of the transformation.
				They are related by
				\begin{equation*}
					\Delta=\delta-\ell^\mu\,\partial_\mu+\mathcal{O}\lr{(\epsilon)}\eqpd
				\end{equation*}
				On the other hand, the measure is affected by a Jacobian and reads
				\begin{equation*}
					\d[D]{x'}=\lr{|\frac{\partial x'}{\partial x}|}\,\d[D]{x}=\lr{(1+\epsilon\,\partial_\mu\ell^\mu+\mathcal{O}\lr{(\epsilon^2)})}\d[D]{x}\eqpd
				\end{equation*}
				Then, the variation of the action is
				\begin{align*}
					\delta\mathcal{S}&=\int\d[D]{x}\,\delta\mathcal{L}+\int\mathcal{L}\,\delta\d[D]{x}\\
					&=\int\d[D]{x}\,\lr{(\delta\mathcal{L}+\mathcal{L}\,\partial_\mu\ell^\mu)}\\
					&=\int\d[D]{x}\,\lr{(\Delta\mathcal{L}+\partial_\mu\lr{(\mathcal{L}\,\ell^\mu)})}\eqpd
				\end{align*}
				The global variation of the Lagrangian density reads
				\begin{align*}
					\Delta\mathcal{L}&=\frac{\partial\mathcal{L}}{\partial\phi^i}\,\Delta\phi^i+\frac{\partial\mathcal{L}}{\partial\lr{(\partial_\mu\phi^i)}}\,\Delta\partial_\mu\phi^i\\
					&=\partial_\mu\frac{\partial\mathcal{L}}{\partial\lr{(\partial_\mu\phi^i)}}\,\Delta\phi^i+\frac{\partial\mathcal{L}}{\partial\lr{(\partial_\mu\phi^i)}}\,\partial_\mu\Delta\phi^i\\
					&=\partial_\mu\lr{(\frac{\partial\mathcal{L}}{\partial\lr{(\partial_\mu\phi^i)}}\,\Delta\phi^i)}\eqpd
				\end{align*}
				(using the \textsc{Euler}~-- \textsc{Lagrange} equation and the fact that the global variation commutes with the space-time derivatives).
				Hence, the variation of the action
				\begin{equation*}
					\delta\mathcal{S}=0=\int\d[D]{x}\,\partial_\mu\lr{(\frac{\partial\mathcal{L}}{\partial\lr{(\partial_\mu\phi^i)}}\,\Delta\phi^i+\mathcal{L}\,\ell^\mu)}\eqpc
				\end{equation*}
				defines a conserved \textsc{Noether} current associated to the symmetry of the action under transformation (\ref{preliminaries.dwh.localtransform.stresstensor}):
				\begin{align*}
					j^\mu&=\frac{\partial\mathcal{L}}{\partial\lr{(\partial_\mu\phi^i)}}\,\Delta\phi^i+\mathcal{L}\,\ell^\mu\\
					&=\frac{\partial\mathcal{L}}{\partial\lr{(\partial_\mu\phi^i)}}\,\delta\phi^i-\lr{(\frac{\partial\mathcal{L}}{\partial\lr{(\partial_\mu\phi^i)}}\,\partial^\nu\phi^i-\eta^{\mu\nu}\mathcal{L})}\ell_\nu\eqpd
				\end{align*}
				Now, if we restrict ourselves to global transformations only ($\ell$ becomes constant), the \textsc{Noether} current associated with the global translational invariance along $\ell$ is
				\begin{equation*}
					j^\mu=\lr{(\eta^{\mu\nu}\mathcal{L}-\frac{\partial\mathcal{L}}{\partial\lr{(\partial_\mu\phi^i)}}\,\partial^\nu\phi^i)}\ell_\nu\eqpc
				\end{equation*}
				since all the fields (whatever their nature) have no local variation (the Jacobian is the identity).
				Now, the stress-energy tensor is naturally defined as the collection of the \textsc{Noether} currents, each associated with the global translational invariance along a direction of space-time.
				It is then defined as
				\begin{equation*}
					\Theta^{\mu\nu}=\eta^{\mu\nu}\mathcal{L}-\frac{\partial\mathcal{L}}{\partial\lr{(\partial_\mu\phi^i)}}\,\partial^\nu\phi^i\eqpc
				\end{equation*}
				and is, by construction, subject to a local conservation law
				\begin{equation*}
					\partial_\mu\Theta^{\mu\nu}=0\eqpd
				\end{equation*}
				Note that this tensor is not necessarily symmetric but it can be symmetrized following the popular \textsc{Belinfante} procedure.
				See, \emph{e.g.}~\cite{PeterUzan2009}.
				
				The stress-energy tensor is expressed in terms of the Hamiltonian density (up to a sign) as
				\begin{equation*}
					\Theta^{\mu\nu}=\eta^{\mu\nu}\mathcal{H}+{\psi_i}^\mu\partial^\nu\phi^i-\eta^{\mu\nu}{\psi_i}^\kappa\partial_\kappa\phi^i\eqpd
				\end{equation*}
				$\bm{\omega}^\mu\lr{(\partial^\nu\bm{\zeta},\bm{\zeta})}/2$ produces terms of the form $\lr{({\psi_i}^\mu\partial^\nu\phi^i-\partial^\nu\phi^i{\psi_i}^\mu)}/2$, which, after integration by parts give ${\psi_i}^\mu\partial^\nu\phi^i$.
				Then, adding a term to restore its symmetry, the stress-energy tensor is indeed given by \cref{preliminaries.dwh.stresstensor}.
				
				To conclude this section, let us now define the charges and prove their conservation.
				The charges are the conserved quantities associated to the \textsc{Noether} currents.
				Since the stress-energy tensor is a collection of $D$ \textsc{Noether} currents, it will define $D$ charges as
				\begin{equation*}
					\mathcal{Q}^\mu=\int\d[d]{\bm{x}}\,\mathcal{T}^{0\mu}\eqpd
				\end{equation*}
				The stress-energy tensor is locally conserved (\cref{preliminaries.dwh.locstressconserv}) and hence the charges are subject to a global conservation law (using \textsc{Stokes}' theorem and again assuming there are no boundary terms):
				\begin{equation*}
					\frac{\d{}\mathcal{Q}^\mu}{\d{x^0}}=\int\d[d]{\bm{x}}\,\partial_0\mathcal{T}^{0\mu}=-\int\d[d]{\bm{x}}\,\partial_j\mathcal{T}^{j\mu}=0\eqpc
				\end{equation*}
				where $j$ runs on $\lr{\llbracket1,d\rrbracket}$.
			\subsubsection{Summary}
				Throughout this section, we have first introduced the notion of symplecticity and highlighted its central r\^ole in Hamiltonian mechanics.
				We have then presented the \textsc{De Donder}~-- \textsc{Weyl} covariant Hamiltonian formulation of field theory and we have shown that it leads to a multi-symplectic phase space.
				We have also introduced (relying on the particular example of the non-linear wave equation) a construction that allows one to obtain a theory, equivalent to the original one, but which does not lead to a degeneracy of the multi-symplectic structure of the phase space.
				Finally, we have proved the local (on-shell) conservation of multi-symplecticity.
				We have also constructed the stress-energy tensor and the charges, and we proved their local and global conservation (on-shell), respectively.
				
				To do so, following~\cite{Bridges1997,BridgesReich2001,BridgesReich2006} we introduced the general form (\ref{preliminaries.dwh.hamilteq.abstractform}) of the \textsc{De Donder}~-- \textsc{Weyl}~-- \textsc{Hamilton} equations.
				Now we claim that, whenever a \textsc{pde} can be written in the form (\ref{preliminaries.dwh.hamilteq.abstractform}), and has a non-degenerate multi-symplectic structure, the \textsc{msilcc} method can be applied.
				We are going to present how to implement it in the next section.
	\section{Construction and properties}\label{implementation}
		The \textsc{msilcc} method is a centered box (finite-difference) scheme~\cite{BridgesReich2001} except that we do not implement it on the traditional hypercubic lattice.
		Instead, we use a lattice based on light-cone coordinates, that has the advantage to restore the locality of the method (\emph{i.e.} there is the same number of unknowns as equations in each cell).
		
		In this section we will present the implementation of the \textsc{msilcc} method as well as a review of some of its interesting properties.
		We will illustrate it on the example of the $\lambda\,\phi^4$ theory mainly in $1+1$ dimensions.
		\subsection{The lattice: sampling the space-time manifold}
			Let us consider the new coordinate system
			\begin{equation}
				\begin{split}
					\check{x}^0&=\frac{x^0}{\sqrt{2}}-\frac{1}{\sqrt{2}}\sum_{j=1}^dx^j \eqpc\\
					\check{x}^j&=\sqrt{2}\,x^j+\delta^j_d\,\check{x}^0\eqpc
				\end{split}\label{implementation.lattice.lcc}
			\end{equation}
			the inverse of which is
			\begin{align*}
				x^0&=\frac{1}{\sqrt{2}}\sum_{\mu=0}^d\check{x}^\mu \eqpc\\
				x^j&=\frac{1}{\sqrt{2}}\lr{(\check{x}^j-\delta^j_d\,\check{x}^0)}\eqpd
			\end{align*}
			Defining $\check{\partial}_\mu=\nicefrac{\displaystyle{\partial}}{\displaystyle{\partial\check{x}^\mu}}$, the associated vector basis is
			\begin{align*}
				\check{\partial}_0&=\frac{\partial_0-\partial_d}{\sqrt{2}} \eqpc\\
				\check{\partial}_j&=\frac{\partial_0+\partial_j}{\sqrt{2}}\eqpc
			\end{align*}
			and its inverse
			\begin{align*}
				\partial_0&=\frac{\check{\partial}_0+\check{\partial}_d}{\sqrt{2}} \eqpc\\
				\partial_j&=\sqrt{2}\,\check{\partial}_j-\partial_0\eqpd
			\end{align*}
			One can remark that (\ref{implementation.lattice.lcc}) is not the usual light-cone coordinate system (except in dimension $D=1+1$).
			The difference is mainly that the set $\lr{\{\check{\partial}_\mu\}}$ is not an orthogonal basis (while it is for the usual definition).
			
			Let us now sample the space-time manifold $\mathcal{M}$ using the lattice
			\begin{equation}
				M=\lr{\{n=o+\delta\,\check{n}^\mu\,\check{\partial}_\mu\vphantom{\check{n}^\mu\in\mathbb{Z},n\in\mathcal{D}}|}\lr{.\vphantom{n=x_o+\delta\,\check{n}^\mu\,\check{\partial}_\mu}\check{n}^\mu\in\mathbb{Z},n\in\mathcal{M}\}}\eqpc\label{implementation.lattice.lattice}
			\end{equation}
			where $o$ is arbitrary (chosen such that $M$ respects as much as possible the boundaries of $\mathcal{M}$) and $\delta$ is the lattice spacing.
			
			At each point on the lattice ($n\in M$) we define the elementary cell (the definition can be extended to each point $n^C\in M+\delta^C$)
			\begin{equation*}
				\text{cell}\lr{(n)}=\lr{\{n+\delta\,\frac{\partial_0+\sigma\,\partial_\rho}{\sqrt{2}}\vphantom{\sigma=\pm1,\rho\in\lr{\llbracket0,d\rrbracket}}|}\lr{.\vphantom{x+\delta\,\frac{\partial_0+\sigma\,\partial_\rho}{\sqrt{2}}}\sigma=\pm1,\rho\in\lr{\llbracket0,d\rrbracket}\}}\eqpc
			\end{equation*}
			if all the vertices of the cell belong in $M$ (possibly using periodic boundary conditions).
			
\begin{figure}[!htb]
				\begin{center}
					\includegraphics{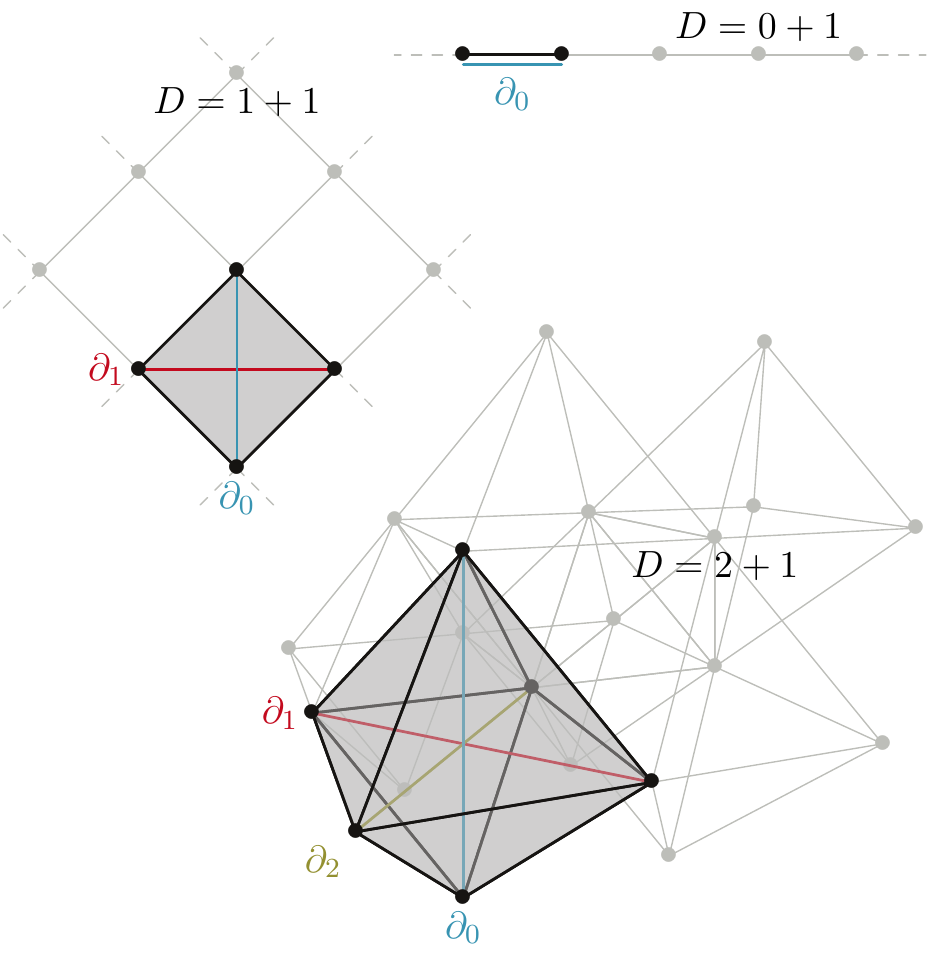}
				\end{center}
				\caption{A portion of the lattice $M$ in $D=0+1$, $1+1$ and $2+1$ dimensions where we have highlighted the elementary cell.}
				\label{implementation.lattice.representation}
			\end{figure}

			\Cref{implementation.lattice.representation} represents how the lattice looks like in low dimensions.
			On this lattice, the approximation of the derivatives along the original directions of space-time (highlighted in color) will simply be obtained using the midpoint rule.
			Hence, they only involve two points of $M$ and this makes the method local.
			
			Let us now highlight that the direction $\check{\partial}_\mu$ selects in $\text{cell}\lr{(n)}$ one (and only one) square, of width $\delta$, with vertices
			\begin{align*}
				n^R_{-|\check{\partial}_\mu}&=n+\delta^R_{-|\check{\partial}_\mu}=n\\
				n^R_{+|\check{\partial}_\mu}&=n+\delta^R_{+|\check{\partial}_\mu}=n^R_{-|\check{\partial}_\mu}+\delta\,\check{\partial}_\mu=n+\delta\,\frac{\partial_0+\sigma\,\partial_\rho}{\sqrt{2}}\\
				n^L_{-|\check{\partial}_\mu}&=n+\delta^L_{-|\check{\partial}_\mu}=n^R_{+|\check{\partial}_\mu}-\sigma\,\delta\,\sqrt{2}\,\partial_\rho=n+\delta\,\frac{\partial_0-\sigma\,\partial_\rho}{\sqrt{2}}\\
				n^L_{+|\check{\partial}_\mu}&=n+\delta^L_{+|\check{\partial}_\mu}=n^L_{-|\check{\partial}_\mu}+\delta\,\check{\partial}_\mu=n+\delta\,\sqrt{2}\,\partial_0\eqpc
			\end{align*}
			where the second equation is used to determine $\sigma=\pm1$ and $\rho\in\lr{\llbracket0,d\rrbracket}$.
			The center of this square (which is the center of the cell as well),
			\begin{equation}
				n^C=n+\delta^C=n+\delta\,\frac{\partial_0}{\sqrt{2}}\eqpc\label{implementation.lattice.cellcenter}
			\end{equation}
			is the point where all the approximations are made.
		\subsection{The numerical approximation scheme}
			\subsubsection{Definition}
				The approximation rules are constructed by a concatenation of the centered box scheme, itself a concatenation of midpoint rules, applied on a square of the elementary cell. Indeed, this scheme is the simplest way to construct a symplectic integrator.
				Therefore, by combining midpoint rules in the most symmetric possible way, we expect to conserve this property and thus obtain a simple multi-symplectic integrator (locally well defined thanks to the lattice).
				
				For $\bm{\zeta}$ or one of its components and for $f$ smooth enough, the approximation rules are given by
				\begin{align*}
					f\lr{(\lr{\{\zeta^a\}}\vphantom{n^C})}\lr{(n^C)}&\approx f\lr{(\lr{\{\lr{<\zeta^a>}\lr{(n^C)}\}})}\eqpc\\
					\check{\partial}_\mu f\lr{(\lr{\{\zeta^a\}}\vphantom{n^C})}\lr{(n^C)}&\approx\check{D}_\mu f\lr{(\lr{\{\zeta^a\}}\vphantom{n^C})}\lr{(n^C)}\eqpc
				\end{align*}
				with $n^C$ defined in \cref{implementation.lattice.cellcenter},
				\begin{align}
					\lr{<\zeta^a>}\lr{(n^C)}&=\frac{1}{2D}\sum_{m\in\mathrlap{\text{cell}\lr{(n)}}\phantom{m}}\zeta^a\lr{(m)}
					\eqpc[and]\label{implementation.rules.average}\\
					\begin{split}
						\check{D}_\mu f\lr{(\lr{\{\zeta^a\}}\vphantom{n^C})}\lr{(n^C)}&=\\
						&\mkern-18mu\frac{1}{\delta}\lr{[f\lr{(\lr{\{\frac{1}{2}\lr{[\zeta^a\lr{(n^L_{+|\check{\partial}_\mu})}+\zeta^a\lr{(n^R_{+|\check{\partial}_\mu})}]}\}})}\vphantom{-f\lr{(\lr{\{\frac{1}{2}\lr{[\zeta^a\lr{(n^L_{-|\check{\partial}_\mu})}+\zeta^a\lr{(n^R_{-|\check{\partial}_\mu})}]}\}})}}.}\\
						&\mkern-18mu\lr{.\vphantom{f\lr{(\lr{\{\frac{1}{2}\lr{[\zeta^a\lr{(n^L_{+|\check{\partial}_\mu})}+\zeta^a\lr{(n^R_{+|\check{\partial}_\mu})}]}\}})}}-f\lr{(\lr{\{\frac{1}{2}\lr{[\zeta^a\lr{(n^L_{-|\check{\partial}_\mu})}+\zeta^a\lr{(n^R_{-|\check{\partial}_\mu})}]}\}})}]}\eqpd
					\end{split}\nonumber
				\end{align}
				For the moment we do not know whether these approximation rules respect the rules of differential calculus (we will explore this issue in the following).
				A priori, the algebraic manipulations done in the continuous formulation will not be equivalent to the ones done on the discrete representation.
				Hence, the \textsc{msilcc} scheme should be applied only in the light-cone coordinate system (\emph{i.e.} all the derivatives $\partial_\mu$ have to be re-expressed in term of the derivatives $\check{\partial}_\mu$ before applying the scheme).
				When and only when directly applied on a field, we have
				\begin{align*}
					&\mathrlap{\partial_0\zeta^a\lr{(n^C)}=\frac{\check{\partial}_0+\check{\partial}_d}{\sqrt{2}}\zeta^a\lr{(n^C)}\approx\lr{[\cdots]}\eqpc}\hphantom{\quad\frac{1}{\sqrt{2}\,\delta}\lr{[\zeta^a\lr{(n^C+\delta\,\frac{\partial_\mu}{\sqrt{2}})}-\zeta^a\lr{(n^C-\delta\,\frac{\partial_\mu}{\sqrt{2}})}]}\eqpd}\\
					&\mathrlap{\partial_j\zeta^a\lr{(n^C)}=\frac{2\,\check{\partial}_j-\check{\partial}_0-\check{\partial}_d}{\sqrt{2}}\zeta^a\lr{(n^C)}\approx\lr{[\cdots]}\eqpd}\hphantom{\quad\frac{1}{\sqrt{2}\,\delta}\lr{[\zeta^a\lr{(n^C+\delta\,\frac{\partial_\mu}{\sqrt{2}})}-\zeta^a\lr{(n^C-\delta\,\frac{\partial_\mu}{\sqrt{2}})}]}}
				\end{align*}
				After some straightforward algebraic manipulations, we obtain
				\begin{align}
						&\partial_\mu\zeta^a\lr{(n^C)}\approx D_\mu\zeta^a\lr{(n^C)}=\label{implementation.rules.partialfield}\\*
						&\quad\frac{1}{\sqrt{2}\,\delta}\lr{[\zeta^a\lr{(n^C+\delta\,\frac{\partial_\mu}{\sqrt{2}})}-\zeta^a\lr{(n^C-\delta\,\frac{\partial_\mu}{\sqrt{2}})}]}\eqpd\nonumber
				\end{align}
				\Cref{implementation.rules.partialfield} defines the derivatives of the field along the original directions of space-time as nothing else than the midpoint rule.
				Nevertheless, remember that this is true only for a linear function of the field, otherwise it is necessary to return to $\check{D}_\mu$.
				
				The discrete analogue of the equation of motion (\ref{preliminaries.dwh.hamilteq.abstractform}) in $\text{cell}\lr{(n)}$ at $n^C$ is
				\begin{equation*}
					\begin{split}
						&\bm{M}^\mu\cdot\lr{[\bm{\zeta}\lr{(n+\delta\,\frac{\partial_0+\partial_\mu}{\sqrt{2}})}-\bm{\zeta}\lr{(n+\delta\,\frac{\partial_0-\partial_\mu}{\sqrt{2}})}]}=\\
						&\qquad\sqrt{2}\,\delta\,\bm{\nabla}\mathcal{H}\bigg(\frac{1}{2D}\sum_{\sigma=\pm1}\sum_{\rho\in\lr{\llbracket0,d\rrbracket}}\bm{\zeta}\lr{(n+\delta\,\frac{\partial_0+\sigma\,\partial_\rho}{\sqrt{2}})}\bigg)\eqpd
					\end{split}
				\end{equation*}
				As expected, the approximation of the equation of motion is indeed a concatenation of midpoint rules.
				Let us illustrate how the procedure works with an example.
			\subsubsection{Application to the \texorpdfstring{$\lambda\,\phi^4$}{lambda phi4} theory in \texorpdfstring{$0+1$}{0+1} dimension}
				The mechanical problem is described by the two unknowns $q$ and $p$ that only depend on time.
				We first sample them through $M$
				\begin{align*}
					q_n&=q\lr{(t=n\,\delta)}\eqpc\\
					p_n&=p\lr{(t=n\,\delta)}\eqpd
				\end{align*}
				Then, applying the \textsc{msilcc} scheme, we get the discrete version of the equations of motion
				\begin{align*}
					p_n-p_{n+1}&=\delta\,\frac{q_{n+1}+q_n}{2}\lr{[1+{\lr{(\frac{q_{n+1}+q_n}{2})}}^2]}\eqpc\\
					q_{n+1}-q_n&=\delta\,\frac{p_{n+1}+p_n}{2}\eqpd
				\end{align*}
				In this particular case, one could write down the explicit expressions for $\lr{(q_{n+1},p_{n+1})}$ as functions of $\lr{(q_n,p_n)}$.
				But these expressions are quite cumbersome and it is not worth presenting them here.
				
				This is a symplectic approximation.
			\subsubsection{The \texorpdfstring{$\lambda\,\phi^4$}{lambda phi4} theory in \texorpdfstring{$1+1$}{1+1} dimensions}
				The lattice is now the same as for the \textsc{BDdV} method (see \cref{msilccvsother.bddv.lattice}).
				Again, defining
				\begin{equation*}
					\sigma_n=2\lr{(n\bmod 2)}-1\equiv\pm1\eqpc
				\end{equation*}
				we sample the fields through $M$ as
				\begin{equation*}
					{\zeta^a\,}_n^j=\zeta^a\lr{(x=\sqrt{2}\,\delta\,\lr{[j+\frac{1+\sigma_n}{4}]},\,t=\frac{n\,\delta}{\sqrt{2}})}\eqpc
				\end{equation*}
				where $n\in\mathbb{N}$ and $j\in\lr{\llbracket 0,N\llbracket}$.
				Therefore the discrete version of the equation of motion is given by the set of algebraic equations
				\begin{align*}
					\begin{split}
						{\psi^0\,}_{n+1}^j-{\psi^0\,}_{n-1}^j+\sigma_n\,{\phantom{\psi_0\,}\mathllap{\psi^1\,}}_n^{j+\sigma_n}-\sigma_n\,{\phantom{\psi_0\,}\mathllap{\psi^1\,}}_n^j=\qquad\qquad\qquad&\\
						-\sqrt{2}\,\delta\,\frac{{\phi\,}_{n-1}^j+{\phi\,}_n^j+{\phi\,}_n^{j+\sigma_n}+{\phi\,}_{n+1}^j}{4}\qquad\qquad&\\
						\lr{[1+{\lr{(\frac{{\phi\,}_{n-1}^j+{\phi\,}_n^j+{\phi\,}_n^{j+\sigma_n}+{\phi\,}_{n+1}^j}{4})}}^2]}\mathrlap{\eqpc}&
					\end{split}\\
					\begin{split}
						{\phantom{\psi_0\,}\mathllap{\phi\,}}_{n+1}^j-{\phantom{\psi_0\,}\mathllap{\phi\,}}_{n-1}^j-\sigma_n\,{\phantom{\psi_0\,}\mathllap{\gamma\,}}_n^{j+\sigma_n}+\sigma_n\,{\phantom{\psi_0\,}\mathllap{\gamma\,}}_n^j=\qquad\qquad\qquad&\\
						\sqrt{2}\,\delta\,\frac{{\psi^0\,}_{n-1}^j+{\psi^0\,}_n^j+{\psi^0\,}_n^{j+\sigma_n}+{\psi^0\,}_{n+1}^j}{4}\mathrlap{\eqpc}\qquad&
					\end{split}\\
					\begin{split}
						{\phantom{\psi_0\,}\mathllap{\gamma\,}}_{n+1}^j-{\phantom{\psi_0\,}\mathllap{\gamma\,}}_{n-1}^j+\sigma_n\,{\phantom{\psi_0\,}\mathllap{\phi\,}}_n^{j+\sigma_n}-\sigma_n\,{\phantom{\psi_0\,}\mathllap{\phi\,}}_n^j=\qquad\qquad\qquad&\\
						-\sqrt{2}\,\delta\,\frac{{\psi^1\,}_{n-1}^j+{\psi^1\,}_n^j+{\psi^1\,}_n^{j+\sigma_n}+{\psi^1\,}_{n+1}^j}{4}\mathrlap{\eqpc}\qquad&
					\end{split}\\
					{\phantom{\psi_0\,}\mathllap{\psi^1\,}}_{n+1}^j-{\phantom{\psi_0\,}\mathllap{\psi^1\,}}_{n-1}^j-\sigma_n\,{\psi^0\,}_n^{j+\sigma_n}+\sigma_n\,{\psi^0\,}_n^j=\mathrlap{\,0\eqpd}\qquad\qquad\qquad&
				\end{align*}
				These are the equations used to integrate the $\lambda\,\phi^4$ theory with the \textsc{msilcc} method in the first section.
				Again, they are not implicit but too much complicated to write down in an explicit form.
				Hence they were treated as implicit equations and solved using the \textsc{Levenberg}~-- \textsc{Marquardt} algorithm for non-linear least squares~\cite{NocedalWright2006}.
				
				In this way we generated the data exposed in \cref{msilccvsother}.
		\subsection{Conservation properties}
			\subsubsection{\texorpdfstring{\textsc{Leibniz}}{Leibniz}'s product rule for quadratic forms}
				Let us now explore how the approximation rules behave with respect to the rules of differential calculus.
				We first apply the discrete derivative to a quadratic form.
				After a straightforward but tedious calculation we obtain
				\begin{align}
					&\check{\partial}_\mu\zeta^a\zeta^b\lr{(n^C)}\approx\check{D}_\mu\zeta^a\zeta^b\lr{(n^C)}=\nonumber\\
					&\quad={\lr{<\zeta^a>}}_{\check{\partial}_\mu}\lr{(n^C)}\,\check{D}_\mu\zeta^b\lr{(n^C)}+{\langle\zeta^b\rangle}_{\check{\partial}_\mu}\lr{(n^C)}\,\check{D}_\mu\zeta^a\lr{(n^C)}\nonumber\\
					&\quad\,\hphantom{={\lr{<\zeta^a>}}_{\check{\partial}_\mu}\lr{(n^C)}\,\check{D}_\mu\zeta^b\lr{(n^C)}}\mathllap{\approx\zeta^a\check{\partial}_\mu\zeta^b\lr{(n^C)}}+\zeta^b\check{\partial}_\mu\zeta^a\lr{(n^C)}\eqpc\label{implementation.Leibniz.fieldpartialfield}
				\end{align}
				where the average value on the square selected by $\check{\partial}_\mu$ is
				\begin{equation}
					\begin{split}
						{\lr{<\zeta^a>}}_{\check{\partial}_\mu}\lr{(n^C)}&=\frac{1}{4}\lr{[\zeta^a\lr{(n^L_{+|\check{\partial}_\mu})}+\zeta^a\lr{(n^R_{+|\check{\partial}_\mu})}\vphantom{+\zeta^a\lr{(n^L_{-|\check{\partial}_\mu})}+\zeta^a\lr{(n^R_{-|\check{\partial}_\mu})}}.}\\
						&\quad\,\lr{.\vphantom{\zeta^a\lr{(n^L_{+|\check{\partial}_\mu})}+\zeta^a\lr{(n^R_{+|\check{\partial}_\mu})}}+\zeta^a\lr{(n^L_{-|\check{\partial}_\mu})}+\zeta^a\lr{(n^R_{-|\check{\partial}_\mu})}]}\eqpd
					\end{split}\label{implementation.Leibniz.orientedaverage}
				\end{equation}
				First of all, \cref{implementation.Leibniz.fieldpartialfield} defines the approximation rule for $\zeta^a\check{\partial}_\mu\zeta^b$ such that the \textsc{Leibniz}'s product rule for quadratic forms holds (actually the \textsc{msilcc} scheme was designed for that purpose since it is a simple way to construct an approximation that preserves the multi-symplectic structure).
				As a second remark, the \textsc{Leibniz}'s product rule remains valid on the discrete scheme for quadratic forms only.
				Finally, this is not true for $\partial_\mu$ (except in $D=1+1$ dimensions since ${\lr{<>}}_{\check{\partial}_\mu}$ coincides with $\lr{<>}$ by definition).
				Hence, the necessity to work in the light-cone coordinate system (all the derivative have to be re-expressed in terms of $\check{\partial}$ before making approximations).
				
				To conclude this section, we stress\footnote{When applied on objects of the collection $\lr{\{z^k\}}$ (defined in \cref{implementation.crossderivatives}) the average value on the square selected by $\check{\partial}_\mu$ as well as the full average value are defined on $M+\delta^C$ at $n$ as
					\begin{align*}
						\begin{split}
							{\langle z^k\rangle}_{\check{\partial}_\mu}\lr{(n)}&=\frac{1}{4}\lr{[z^k\lr{(n^L_{+|\check{\partial}_\mu}-\delta^C)}+z^k\lr{(n^R_{+|\check{\partial}_\mu}-\delta^C)}\vphantom{+z^k\lr{(n^L_{-|\check{\partial}_\mu}-\delta^C)}+z^k\lr{(n^R_{-|\check{\partial}_\mu}-\delta^C)}}.}\\
							&\quad\,\lr{.\vphantom{z^k\lr{(n^L_{+|\check{\partial}_\mu}-\delta^C)}+z^k\lr{(n^R_{+|\check{\partial}_\mu}-\delta^C)}}+z^k\lr{(n^L_{-|\check{\partial}_\mu}-\delta^C)}+z^k\lr{(n^R_{-|\check{\partial}_\mu}-\delta^C)}]}
						\end{split}\\
						\langle z^k\rangle\lr{(n)}&=\frac{1}{2D}\sum_{m\in\mathrlap{\text{cell}\lr{(n-\delta^C)}}\phantom{m}}z^k\lr{(m)}\eqpd
					\end{align*}} that the same relations hold on the lattice $M+\delta^C$.
			\subsubsection{Preservation of cross derivatives equality}\label{implementation.crossderivatives}
				Now, we define the collection $\lr{\{z^k\}}$.
				Each $z^k$ lives on $M+\delta^C$ and is linear in the field.
				So the collection $\lr{\{z^k\}}$ is limited to
				\begin{equation*}
					\lr{\{\lr{<\zeta^a>},\lr{\{{\lr{<\zeta^a>}}_{\check{\partial}_\mu},\check{D}_\mu\zeta^a\}}\}}\eqpd
				\end{equation*}
				The average values were defined in the previous section. We add here the definition of the derivatives
				\begin{equation*}
					\begin{split}
						&\check{\partial}_\mu f\textstyle{\lr{(\lr{\{z^k\}})}}\lr{(n)}\approx\check{D}_\mu f\textstyle{\lr{(\lr{\{z^k\}})}}\lr{(n)}=\\
						&\quad\frac{1}{\delta}\lr{[f\lr{(\lr{\{\frac{1}{2}\lr{[z^k\lr{(n^L_{+|\check{\partial}_\mu}-\delta^C)}+z^k\lr{(n^R_{+|\check{\partial}_\mu}-\delta^C)}]}\}})}\vphantom{-f\lr{(\lr{\{\frac{1}{2}\lr{[z^k\lr{(n^L_{-|\check{\partial}_\mu}-\delta^C)}+z^k\lr{(n^R_{-|\check{\partial}_\mu}-\delta^C)}]}\}})}}.}\\
						&\quad\lr{.\vphantom{f\lr{(\lr{\{\frac{1}{2}\lr{[z^k\lr{(n^L_{+|\check{\partial}_\mu}-\delta^C)}+z^k\lr{(n^R_{+|\check{\partial}_\mu}-\delta^C)}]}\}})}}-f\lr{(\lr{\{\frac{1}{2}\lr{[z^k\lr{(n^L_{-|\check{\partial}_\mu}-\delta^C)}+z^k\lr{(n^R_{-|\check{\partial}_\mu}-\delta^C)}]}\}})}]}\eqpd
					\end{split}
				\end{equation*}
				Following the same reasoning as when the field was $\zeta^a$ defined on $M$, we have (again, only when directly applied on a field)
				\begin{equation*}
					\begin{split}
						&\partial_\mu z^k\lr{(n)}\approx D_\mu z^k\lr{(n)}=\\
						&\quad\frac{1}{\sqrt{2}\,\delta}\lr{[z^k\lr{(n+\delta\,\frac{\partial_\mu}{\sqrt{2}})}-z^k\lr{(n-\delta\,\frac{\partial_\mu}{\sqrt{2}})}]}\eqpd
					\end{split}
				\end{equation*}
				Using these definitions, one can give meaning to the second derivative of the field.
				We have
				\begin{align*}
					&D_\mu D_\nu\zeta^a\lr{(n)}=\nonumber\\*
					&\quad=\frac{1}{\sqrt{2}\,\delta}\lr{[D_\nu\zeta^a\lr{(n+\delta\,\frac{\partial_\mu}{\sqrt{2}})}-D_\nu\zeta^a\lr{(n-\delta\,\frac{\partial_\mu}{\sqrt{2}})}]}\\
					\begin{split}
						&\quad=\frac{1}{2\,\delta^2}\lr{[\zeta^a\lr{(n+\delta\,\frac{\partial_\mu+\partial_\nu}{\sqrt{2}})}-\zeta^a\lr{(n+\delta\,\frac{\partial_\mu-\partial_\nu}{\sqrt{2}})}\vphantom{-\zeta^a\lr{(n-\delta\,\frac{\partial_\mu-\partial_\nu}{\sqrt{2}})}+\zeta^a\lr{(n-\delta\,\frac{\partial_\mu+\partial_\nu}{\sqrt{2}})}}.}\\
						&\quad\qquad\,\,\lr{.\vphantom{\zeta^a\lr{(n+\delta\,\frac{\partial_\mu+\partial_\nu}{\sqrt{2}})}-\zeta^a\lr{(n+\delta\,\frac{\partial_\mu-\partial_\nu}{\sqrt{2}})}}-\zeta^a\lr{(n-\delta\,\frac{\partial_\mu-\partial_\nu}{\sqrt{2}})}+\zeta^a\lr{(n-\delta\,\frac{\partial_\mu+\partial_\nu}{\sqrt{2}})}]}
					\end{split}\\
					\begin{split}
						&\quad=\frac{1}{2\,\delta^2}\lr{[\zeta^a\lr{(n+\delta\,\frac{\partial_\nu+\partial_\mu}{\sqrt{2}})}-\zeta^a\lr{(n+\delta\,\frac{\partial_\nu-\partial_\mu}{\sqrt{2}})}\vphantom{-\zeta^a\lr{(n-\delta\,\frac{\partial_\nu-\partial_\mu}{\sqrt{2}})}+\zeta^a\lr{(n-\delta\,\frac{\partial_\nu+\partial_\mu}{\sqrt{2}})}}.}\\
						&\quad\qquad\,\,\lr{.\vphantom{\zeta^a\lr{(n+\delta\,\frac{\partial_\nu+\partial_\mu}{\sqrt{2}})}-\zeta^a\lr{(n+\delta\,\frac{\partial_\nu-\partial_\mu}{\sqrt{2}})}}-\zeta^a\lr{(n-\delta\,\frac{\partial_\nu-\partial_\mu}{\sqrt{2}})}+\zeta^a\lr{(n-\delta\,\frac{\partial_\nu+\partial_\mu}{\sqrt{2}})}]}
					\end{split}\\
					&\quad=\frac{1}{\sqrt{2}\,\delta}\lr{[D_\mu\zeta^a\lr{(n+\delta\,\frac{\partial_\nu}{\sqrt{2}})}-D_\mu\zeta^a\lr{(n-\delta\,\frac{\partial_\nu}{\sqrt{2}})}]}\\
					&\quad=D_\nu D_\mu\zeta^a\lr{(n)}\eqpc
				\end{align*}
				proving the identity of the cross derivatives in discrete space-time.
				Using the relation between $D_\mu$ and $\check{D}_\mu$, we find that the same applies on the light-cone coordinates:
				\begin{equation*}
					\check{D}_\mu\check{D}_\nu\zeta^a\lr{(n)}=\check{D}_\nu\check{D}_\mu\zeta^a\lr{(n)}\eqpd
				\end{equation*}
			\subsubsection{Exact conservation of the multi-symplectic structure}
				Let us now prove the conservation of the multi-symplectic structure.
				We first perform the change of coordinates in the left hand side operator of the equation of motion (\ref{preliminaries.dwh.hamilteq.abstractform})
				\begin{equation*}
					\check{\bm{M}}^\mu\cdot\check{\partial}_\mu=\bm{M}^\mu\cdot\partial_\mu\eqpc
				\end{equation*}
				and we obtain the set of skew-symmetric matrices in the new coordinate system
				\begin{align*}
					\check{\bm{M}}^0&=\frac{1}{\sqrt{2}}\,\bm{M}^0-\frac{1}{\sqrt{2}}\sum_{j=1}^d\bm{M}^j \eqpc\\
					\check{\bm{M}}^j&=\sqrt{2}\,\bm{M}^j+\delta_d^j\,\check{\bm{M}}^0\eqpd
				\end{align*}
				From $\check{\bm{M}}^\mu$ we define $\check{\bm{\omega}}^\mu$ that actually behaves as a component of a $D$-vector in $\mathcal{M}$
				\begin{equation*}
					\check{\bm{\omega}}^\mu=\frac{\partial\check{x}^\mu}{\partial x^\rho}\,\bm{\omega}^\rho\eqpd
				\end{equation*}
				The set $\lr{\{\check{\bm{\omega}}^\mu\}}$ defines the multi-symplectic structure in the light-cone coordinate system and is subject to the same conservation law
				\begin{equation*}
					\check{\partial}_\mu\check{\bm{\omega}}^\mu=\frac{\partial x^\rho}{\partial\check{x}^\mu}\,\partial_\rho\lr{(\frac{\partial\check{x}^\mu}{\partial x^\sigma}\,\bm{\omega}^\sigma)}=\partial_\rho\bm{\omega}^\rho=0\eqpd
				\end{equation*}
				By taking the exterior derivative of the equation of motion (\ref{preliminaries.dwh.hamilteq.abstractform}) we have
				\begin{equation*}
					{\lr{.\check{M}^\mu.}\!}_{ab}\,\check{\partial}_\mu\db[b]{}=\bm{\partial}_a\bm{\partial}_b\mathcal{H}\lr{(\bm{\zeta})}\,\db[b]{}\eqpd
				\end{equation*}
				Then, the local conservation of multi-symplecticity is, numerically,
				\begin{align*}
					&\check{\partial}_\mu\check{\bm{\omega}}^\mu\approx\check{D}_\mu\check{\bm{\omega}}^\mu=\\
					&\quad=-\frac{1}{2}{\lr{.\check{M}^\mu.}\!}_{ab}\,\lr{(\check{D}_\mu\db[a]{}\wedge{\langle\db[b]{}\rangle}_{\check{\partial}_\mu}+{\lr{<\db[a]{}>}}_{\check{\partial}_\mu}\wedge\check{D}_\mu\db[b]{})}\\
					&\quad={\lr{.\check{M}^\mu.}\!}_{ab}\,\check{D}_\mu\db[b]{}\wedge{\lr{<\db[a]{}>}}_{\check{\partial}_\mu}\\
					&\quad=\bm{\partial}_a\bm{\partial}_b\mathcal{H}\lr{(\lr{<\bm{\zeta}>})}\;\langle\db[b]{}\rangle\wedge{\lr{<\db[a]{}>}}_{\check{\partial}_\mu}\\
					&\quad=0\eqpc
				\end{align*}
				since the contraction of the symmetric object $\bm{\partial}_a\bm{\partial}_b$ with the skew-symmetry of the wedge product vanishes.
				
				So, the multi-symplectic structure is indeed exactly preserved by the \textsc{msilcc} scheme which is hence a multi-symplectic integrator.
		\subsection{Conservation of the stress-energy tensor}
			In this section we investigate the effect of the \textsc{msilcc} scheme on the stress-energy tensor.
			\subsubsection{Local approximate conservation of the stress-energy tensor}\label{implementation.stresstensor.conservation}
				Let us first obtain two preliminary results.
				On the one hand, we have the commutativity of the oriented average value (\cref{implementation.Leibniz.orientedaverage}) with itself
				\begin{align*}
					\begin{split}
						&{\big\langle{\lr{<\zeta^a>}}_{\check{\partial}_\mu}\big\rangle}_{\check{\partial}_\nu}\lr{(n)}=\\
						&\quad=\frac{1}{4}\sum_{\sigma=\pm}\sum_{X=\mathrlap{\lr{\{L,R\}}}\phantom{X}}{\lr{<\zeta^a>}}_{\check{\partial}_\mu}\lr{(n-\delta^C+\delta^X_{\sigma|\check{\partial}_\nu})}
					\end{split}\\
					&\quad=\frac{1}{16}\sum_{\substack{\phantom{\sigma'}\mathllap{\sigma}=\pm\\\sigma'=\pm}}\sum_{\substack{\phantom{X'}\mathllap{X}=\mathrlap{\lr{\{L,R\}}}\phantom{X'}\\ X'=\mathrlap{\lr{\{L,R\}}}\phantom{X'}}}\zeta^a\lr{(n-2\,\delta^C+\delta^X_{\sigma|\check{\partial}_\nu}+\delta^{X'}_{\sigma'|\check{\partial}_\mu})}\\
					&\quad=\frac{1}{4}\sum_{\sigma'=\pm}\sum_{X'=\mathrlap{\lr{\{L,R\}}}\phantom{X}}{\lr{<\zeta^a>}}_{\check{\partial}_\nu}\lr{(n-\delta^C+\delta^{X'}_{\sigma'|\check{\partial}_\mu})}\\
					&\quad={\big\langle{\lr{<\zeta^a>}}_{\check{\partial}_\nu}\big\rangle}_{\check{\partial}_\mu}\lr{(n)}\eqpd
				\end{align*}
				On the other hand, we have the commutativity of the oriented average value with the derivative
				\begin{align*}
					\begin{split}
						&{\big\langle\check{D}_\mu\zeta^a\big\rangle}_{\check{\partial}_\nu}\lr{(n)}=\\
						&\quad=\frac{1}{4}\sum_{\sigma=\pm}\sum_{X=\mathrlap{\lr{\{L,R\}}}\phantom{X}}\check{D}_\mu\zeta^a\lr{(n-\delta^C+\delta^X_{\sigma|\check{\partial}_\nu})}
					\end{split}\\
					&\quad=\frac{1}{8\,\delta}\sum_{\substack{\phantom{\sigma'}\mathllap{\sigma}=\pm\\\sigma'=\pm}}\sum_{\substack{\phantom{X'}\mathllap{X}=\mathrlap{\lr{\{L,R\}}}\phantom{X'}\\ X'=\mathrlap{\lr{\{L,R\}}}\phantom{X'}}}\sigma'\,\zeta^a\lr{(n-2\,\delta^C+\delta^X_{\sigma|\check{\partial}_\nu}+\delta^{X'}_{\sigma'|\check{\partial}_\mu})}\\
					&\quad=\frac{1}{2\,\delta}\sum_{\sigma'=\pm}\sum_{X'=\mathrlap{\lr{\{L,R\}}}\phantom{X}}\sigma'\,{\lr{<\zeta^a>}}_{\check{\partial}_\nu}\lr{(n-\delta^C+\delta^{X'}_{\sigma'|\check{\partial}_\mu})}\\
					&\quad=\check{D}_\mu{\lr{<\zeta^a>}}_{\check{\partial}_\nu}\lr{(n)}\eqpd
				\end{align*}
				Now we consider the non-symmetrized part of the stress-energy tensor
				\begin{equation*}
					\mathcal{T}^{\mu\nu}=\frac{1}{2}\lr{[\bm{\omega}^\nu\lr{(\partial^\mu\bm{\zeta},\bm{\zeta})}-\eta^{\mu\nu}\bm{\omega}^\kappa\lr{(\partial_\kappa\bm{\zeta},\bm{\zeta})}]}+\eta^{\mu\nu}\mathcal{H}\lr{(\bm{\zeta})}\eqpd
				\end{equation*}
				Since $\mathcal{T}$ is a tensor ($\eta$ is a tensor and $\bm{\omega}$ and $\partial$ are vectors), we have
				\begin{align*}
					\check{\mathcal{T}}^{\mu\nu}&=\frac{\partial\check{x}^\mu}{\partial x^\rho}\frac{\partial\check{x}^\nu}{\partial x^\sigma}\,\mathcal{T}^{\rho\sigma}\\
					&=\frac{1}{2}\lr{[\check{\bm{\omega}}^\nu\lr{(\check{\partial}^\mu\bm{\zeta},\bm{\zeta})}-\check{\eta}^{\mu\nu}\check{\bm{\omega}}^\kappa\lr{(\check{\partial}_\kappa\bm{\zeta},\bm{\zeta})}]}+\check{\eta}^{\mu\nu}\mathcal{H}\lr{(\bm{\zeta})}\\
					&\approx\check{T}^{\mu\nu}\eqpc
				\end{align*}
				where
				\begin{equation*}
					\check{\eta}^{\mu\nu}=\partial_\rho\check{x}^\mu\,\partial^\rho\check{x}^\nu\eqpd
				\end{equation*}
				The numerical version of $\check{\mathcal{T}}$ is defined (using the approximation rules introduced earlier) as
				\begin{equation*}
					\begin{split}
						\check{T}^{\mu\nu}&=\frac{1}{2}\lr{[\check{\bm{\omega}}^\nu\lr{(\check{D}^\mu\bm{\zeta},{\lr{<\bm{\zeta}>}}_{\check{\partial}^\mu})}-\check{\eta}^{\mu\nu}\check{\bm{\omega}}^\kappa\lr{(\check{D}_\kappa\bm{\zeta},{\lr{<\bm{\zeta}>}}_{\check{\partial}_\kappa})}]}\\
						&\qquad+\check{\eta}^{\mu\nu}\mathcal{H}\lr{(\lr{<\bm{\zeta}>})}\eqpd
					\end{split}
				\end{equation*}
				Now, we use the exact conservation of the multi-symplectic structure
				\begin{align*}
					&\check{D}_\nu\check{\bm{\omega}}^\nu\lr{(\check{D}^\mu\bm{\zeta},{\lr{<\bm{\zeta}>}}_{\check{\partial}^\mu})}=0=\\
					&\quad=\check{\bm{\omega}}^\nu\lr{(\check{D}_\nu\check{D}^\mu\bm{\zeta},{\big\langle{\lr{<\bm{\zeta}>}}_{\check{\partial}^\mu}\big\rangle}_{\check{\partial}_\nu})}+\check{\bm{\omega}}^\nu\lr{({\big\langle\check{D}^\mu\bm{\zeta}\big\rangle}_{\check{\partial}_\nu},\check{D}_\nu{\lr{<\bm{\zeta}>}}_{\check{\partial}^\mu})}\eqpc
				\end{align*}
				and the dual of the equation of motion (\ref{preliminaries.dwh.hamilteq.abstractform}) to prove that
				\begin{align*}
					\check{\bm{\omega}}^\nu\lr{(\check{D}_\nu\check{D}^\mu\bm{\zeta},{\big\langle{\lr{<\bm{\zeta}>}}_{\check{\partial}^\mu}\big\rangle}_{\check{\partial}_\nu})}&=\check{\bm{\omega}}^\nu\lr{(\check{D}_\nu{\lr{<\bm{\zeta}>}}_{\check{\partial}^\mu},{\big\langle\check{D}^\mu\bm{\zeta}\big\rangle}_{\check{\partial}_\nu})}\\
					&=\db{\mathcal{H}}\lr{(\big\langle\lr{<\bm{\zeta}>}\big\rangle)}\lr{[\check{D}^\mu\lr{<\bm{\zeta}>}]}\eqpd
				\end{align*}
				Then, using all the preliminary results of this section, the approximation of the local conservation of the stress-energy tensor reads
				\begin{align}
					\check{D}_\nu\check{T}^{\mu\nu}&=\frac{1}{2}\lr{[0-\check{D}^\mu\check{\bm{\omega}}^\kappa\lr{(\check{D}_\kappa\bm{\zeta},{\lr{<\bm{\zeta}>}}_{\check{\partial}_\kappa})}]}+\check{D}^\mu\mathcal{H}\lr{(\lr{<\bm{\zeta}>})}\nonumber\\
					\begin{split}
						&=\check{D}^\mu\mathcal{H}\lr{(\lr{<\bm{\zeta}>})}-\frac{1}{2}\lr{[\check{\bm{\omega}}^\kappa\lr{(\check{D}^\mu\check{D}_\kappa\bm{\zeta},{\big\langle{\lr{<\bm{\zeta}>}}_{\check{\partial}_\kappa}\big\rangle}_{\check{\partial}^\mu})}\vphantom{+\check{\bm{\omega}}^\kappa\lr{({\big\langle\check{D}_\kappa\bm{\zeta}\big\rangle}_{\check{\partial}^\mu},\check{D}^\mu{\lr{<\bm{\zeta}>}}_{\check{\partial}_\kappa})}}.}\\
						&\qquad\qquad\qquad\qquad\lr{.\vphantom{\check{\bm{\omega}}^\kappa\lr{(\check{D}^\mu\check{D}_\kappa\bm{\zeta},{\big\langle{\lr{<\bm{\zeta}>}}_{\check{\partial}_\kappa}\big\rangle}_{\check{\partial}^\mu})}}+\check{\bm{\omega}}^\kappa\lr{({\big\langle\check{D}_\kappa\bm{\zeta}\big\rangle}_{\check{\partial}^\mu},\check{D}^\mu{\lr{<\bm{\zeta}>}}_{\check{\partial}_\kappa})}]}
					\end{split}\nonumber\\
					\begin{split}
						&=\check{D}^\mu\mathcal{H}\lr{(\lr{<\bm{\zeta}>})}-\frac{1}{2}\lr{[\check{\bm{\omega}}^\kappa\lr{(\check{D}_\kappa\check{D}^\mu\bm{\zeta},{\big\langle{\lr{<\bm{\zeta}>}}_{\check{\partial}^\mu}\big\rangle}_{\check{\partial}_\kappa})}\vphantom{+\check{\bm{\omega}}^\kappa\lr{(\check{D}_\kappa{\lr{<\bm{\zeta}>}}_{\check{\partial}^\mu},{\big\langle\check{D}^\mu\bm{\zeta}\big\rangle}_{\check{\partial}_\kappa})}}.}\\
						&\qquad\qquad\qquad\qquad\lr{.\vphantom{\check{\bm{\omega}}^\kappa\lr{(\check{D}_\kappa\check{D}^\mu\bm{\zeta},{\big\langle{\lr{<\bm{\zeta}>}}_{\check{\partial}^\mu}\big\rangle}_{\check{\partial}_\kappa})}}+\check{\bm{\omega}}^\kappa\lr{(\check{D}_\kappa{\lr{<\bm{\zeta}>}}_{\check{\partial}^\mu},{\big\langle\check{D}^\mu\bm{\zeta}\big\rangle}_{\check{\partial}_\kappa})}]}
					\end{split}\nonumber\\
					&=\check{D}^\mu\mathcal{H}\lr{(\lr{<\bm{\zeta}>})}-\db{\mathcal{H}}\lr{(\big\langle\lr{<\bm{\zeta}>}\big\rangle)}\lr{[\check{D}^\mu\lr{<\bm{\zeta}>}]}\nonumber\\
					&=\check{D}^\mu\mathcal{H}_I\lr{(\lr{<\bm{\zeta}>})}-\db{\mathcal{H}_I}\lr{(\big\langle\lr{<\bm{\zeta}>}\big\rangle)}\lr{[\check{D}^\mu\lr{<\bm{\zeta}>}]}\eqpc\label{implementation.stresstensor.locconservexacterror}
				\end{align}
				where
				\begin{equation*}
					\mathcal{H}_I=\mathcal{H}-\mathcal{H}_Q\eqpc
				\end{equation*}
				is the non-quadratic part of $\mathcal{H}$ ($\mathcal{H}_Q$ is the quadratic part of the Hamiltonian density).
				
				Accordingly, the \textsc{msilcc} scheme exactly preserves the local conservation of the stress-energy tensor for any linear Hamiltonian \textsc{pde}.
				When applied on a non-linear problem, the \textsc{msilcc} method breaks the conservation of the stress-energy tensor only because the chain rule does not hold on the discrete space-time.
				Nevertheless, if the sampling is good enough we expect the \textsc{msilcc} integrator not generate large violations of this conservation law.
				
				One can remark that there is no longer any second derivative in \cref{implementation.stresstensor.locconservexacterror}.
				Hence, let us approximate $\check{D}_\nu\check{T}^{\mu\nu}$ by removing the lowest level average value:
				\begin{equation}
					\check{D}_\nu\check{T}^{\mu\nu}\simeq\check{D}^\mu\mathcal{H}_I\lr{(\bm{\zeta})}-\db{\mathcal{H}_I}\lr{(\lr{<\bm{\zeta}>})}\lr{[\check{D}^\mu\bm{\zeta}]}\eqpd\label{implementation.stresstensor.locconservapproxerror}
				\end{equation}
				Obviously, this operation is strictly forbidden.
				Nevertheless, \cref{implementation.stresstensor.locconservapproxerror} is a very good estimator of \cref{implementation.stresstensor.locconservexacterror}.
				This can be understood if we remember that \cref{implementation.stresstensor.locconservexacterror} mainly evaluates how much the chain rule is violated for non-quadratic functions on the lattice.
				Therefore, increasing the averaging is not an essential element.
				
				In practice, on the example of the $\lambda\,\phi^4$ theory in $1+1$ dimensions the difference between \cref{implementation.stresstensor.locconservexacterror} (or explicitly \cref{implementation.phi4.stresstensor.locconservexacterror}) and \cref{implementation.stresstensor.locconservapproxerror} (explicitly \cref{implementation.phi4.stresstensor.locconservapproxerror}) is negligible and it is almost impossible to distinguish the two on the numerical results.
				
				The substantial advantage of the estimator (\ref{implementation.stresstensor.locconservapproxerror}) is that it is simpler to compute, but first and foremost, that it is more local (it involves only the current cell).
				Thus, the accuracy of the integration can be checked regardless of the neighbouring cells.
				This ensures a better scalability of the method by reducing the number of communications.
			\subsubsection{The \texorpdfstring{$\lambda\,\phi^4$}{lambda phi4} theory in \texorpdfstring{$1+1$}{1+1} dimensions}
				In the case of the $\lambda\,\phi^4$ theory in $1+1$ dimensions, and assuming that the extra field $\gamma$ is free (\emph{i.e.} used as a control parameter), \cref{implementation.stresstensor.locconservexacterror} explicitly becomes
				\begin{align}
					\epsilon^\pm&=\check{D}_\nu\check{T}^{\mu\nu}\nonumber\\
						&=\frac{1}{4\,\delta}\lr{[{\lr{(\frac{\phi_U+\phi_\pm}{2})}}^4-{\lr{(\frac{\phi_D+\phi_\mp}{2})}}^4]}\label{implementation.phi4.stresstensor.locconservexacterror}\\*
						&\;-\frac{1}{2\,\delta}\lr{(\phi_U+\phi_\pm-\phi_D-\phi_\mp)}{\lr{(\frac{\phi_U+\phi_\pm+\phi_D+\phi_\mp}{4})}}^3\eqpc\nonumber
				\end{align}
				where $\pm$ selects $\mu=0$ or $1$ and
				\begin{align*}
					\phi_U&=\frac{1}{4}\lr{[{\phi\,}_n^j+{\phi\,}_{n+1}^j+{\phi\,}_{n+1}^{j-\sigma_n}+{\phi\,}_{n+2}^j]}\eqpc\\
					\phi_\pm&=\frac{1}{4}\lr{[{\phi\,}_n^{j+\frac{\sigma_n\pm1}{2}}+{\phi\,}_n^{j-\frac{\sigma_n\mp1}{2}}+{\phi\,}_{n-1}^{j-\frac{\sigma_n\mp1}{2}}+{\phi\,}_{n+1}^{j-\frac{\sigma_n\mp1}{2}}]}\eqpc\\
					\phi_D&=\frac{1}{4}\lr{[{\phi\,}_n^j+{\phi\,}_{n-1}^j+{\phi\,}_{n-1}^{j-\sigma_n}+{\phi\,}_{n-2}^j]}\eqpc
				\end{align*}
				while \cref{implementation.stresstensor.locconservapproxerror} explicitly becomes
				\begin{align}
						&\epsilon^\pm\simeq\nonumber\\*
						&\quad\frac{1}{4\,\delta}\lr{[{\lr{(\frac{{\phi\,}_n^{j+\frac{\sigma_n\pm1}{2}}+{\phi\,}_{n+1}^j}{2})}}^4-{\lr{(\frac{{\phi\,}_{n-1}^j+{\phi\,}_n^{j+\frac{\sigma_n\mp1}{2}}}{2})}}^4]}\nonumber\\*
						&\qquad-\frac{1}{2\,\delta}\lr{({\phi\,}_n^{j+\frac{\sigma_n\pm1}{2}}+{\phi\,}_{n+1}^j-{\phi\,}_{n-1}^j-{\phi\,}_n^{j+\frac{\sigma_n\mp1}{2}})}\nonumber\\*
						&\qquad\quad\,\,\times{\lr{(\frac{{\phi\,}_n^j+{\phi\,}_n^{j+\sigma_n}+{\phi\,}_{n-1}^j+{\phi\,}_{n+1}^j}{4})}}^3\eqpd\label{implementation.phi4.stresstensor.locconservapproxerror}
				\end{align}
			\subsubsection{Note on the possibility of an exact conservation of the stress-energy tensor}
				It is actually possible to obtain an exact conservation of the stress-energy tensor.
				Let us remember the idea of the \textsc{BDdV} method: the discretization rules were applied on the energy instead of the equation of motion, then the constrains on the conservation of the energy were used as an equation of motion.
				These two procedures are equivalent in the continuum limit, but they are not on the lattice since the rules of differential calculus are no longer fulfilled in the latter setting.
				
				One can imagine here to proced in the same way by applying the discretization rules on the stress-energy tensor and then use its conservation as an equation of motion (hence an exact conservation of the stress-energy tensor).
				However, it would become necessary to evaluate the error committed on the original equation of motion.
				This would leads to evaluate the quantity: $\check{D}^\mu\mathcal{H}_I\lr{(\lr{<\bm{\zeta}>})}-\db{\mathcal{H}_I}\lr{(\big\langle\lr{<\bm{\zeta}>}\big\rangle)}\lr{[\check{D}^\mu\lr{<\bm{\zeta}>}]}$ (\emph{i.e.} \cref{implementation.stresstensor.locconservexacterror}).
				
				The $\db{\mathcal{H}_I}\lr{(\big\langle\lr{<\bm{\zeta}>}\big\rangle)}\lr{[\check{D}^\mu\lr{<\bm{\zeta}>}]}$ term is the right hand side of the equation of motion while $\check{D}^\mu\mathcal{H}_I\lr{(\lr{<\bm{\zeta}>})}$ arises with the derivatives of the stress-energy tensor.
				So, whether the discretization is performed on the equation of motion or on the stress-energy tensor, to estimate the quality of the approximation we have to evaluate how much $\db{\mathcal{H}_I}\lr{(\big\langle\lr{<\bm{\zeta}>}\big\rangle)}\lr{[\check{D}^\mu\lr{<\bm{\zeta}>}]}$ differs from $\check{D}^\mu\mathcal{H}_I\lr{(\lr{<\bm{\zeta}>})}$ in both cases.
		\subsection{Motivation to use the light-cone coordinates}
			Solving the numerical equation of motion requires the simultaneous solution of a set of algebraic equations.
			This is equivalent to finding the root of a vector function.
			It is generally preferable to see this problem as the minimization of the square norm of a vector function since algorithms for this kind of optimization are more robust and more diversified than the ones for finding roots.
			Indeed, this can be achieved by using standard optimization methods such as the \textsc{Levenberg}~-- \textsc{Marquardt}, \textsc{Powell}'s Dog Leg, etc.~\cite{NocedalWright2006}.
			
			We alluded to this feature earlier, but we now want to stress the importance of the lattice.
			It has been chosen such that in each cell there is only one point at the latest time.
			Thus, in each cell we have as many algebraic equations as unknowns.
			The method is well defined locally.
			Usually, the centered box scheme is implemented on a hypercubic lattice which is indeed simpler but leads to more unknowns than equations in each cell (except in dimension $D=0+1$).
			The method is still globally well defined since each unknown is involved in the equations of the neighbouring cells.
			However, at each time step, it requires to solve the whole system in one block.
			Therefore, if we want to dispatch the problem on several process units a huge number of communications are needed (known to be a bottleneck for high performance computations).
			
			The main advantage of the \textsc{msilcc} method, is that it restores the locality of the algorithm while most of the expressions (the equation of motion, the conservation of the stress-energy tensor, \dots) remain quite simple as we have shown through the example of the $\lambda\,\phi^4$ theory.
			
			We finally want to make a remark concerning the initial conditions: the lattice of the \textsc{msilcc} method is such that a cell involves three levels of time.
			Therefore, at the initial time, in each cell, we have two unknowns for only one equation.
			The idea to solve this tricky problem is to assume (only at the initial point) that the average in space is equal to the average in time (\emph{i.e.} the average over all the points of the cell at $t=0$ is equal to the average of the two points at $t=0\pm\delta^C$).
			In this way we have removed the superfluous unknowns.
			Nevertheless, it requires that the equation of motion contains a derivative along time of all the fields of the vector state.
			Hence the necessity to work with a formulation of the problem that will not lead to a degeneracy of the multi-symplectic structure.
		\subsection{Higher dimensions}
			In the present section we discuss the limits of the \textsc{msilcc} method and possible ways of improvement.
			
			The lattice $M$, defined in \cref{implementation.lattice.lattice}, is an attempt to generalize to higher dimensions the one introduced in \cref{msilccvsother.bddv.latticesection} for $D=1+1$ (see \cref{msilccvsother.bddv.lattice}).
			However, we have experienced some instabilities of the method in dimension $D>1+1$.
			This change in the behavior of our method when the dimension of space-time increases is a bit astonishing.
			We suspect two reasons for that.
			When the dimension of space-time becomes higher than $1+1$:
			\begin{enumerate}[label = \emph{\roman{*}}., labelindent = 0em, leftmargin = *, widest* = 2, nosep]
				\item On the one hand, the ensemble of the cells of the lattice is no longer a tessellation of the space-time manifold (\emph{i.e.} there are points in space-time which are not contained in any cell).
				\item On the other hand, the oriented average, ${\lr{<\cdots>}}_{\check{\partial}_\mu}$ (see \cref{implementation.Leibniz.orientedaverage}), no longer coincides with the full average $\lr{<\cdots>}$ (\cref{implementation.rules.average}).
			\end{enumerate}
			
			One can imagine another generalization of the lattice introduced in \cref{msilccvsother.bddv.latticesection} which avoids the two problems mentioned above.
			This is a hypercubic lattice, oriented in such a way that there is only one unknown in each cell (cells are now hypercubes).
			It consists in starting with another light-cone coordinate system:
			\begin{equation*}
				\check{\partial}_\mu=\frac{\partial x^\nu}{\partial\check{x}^\mu}\,\partial_\nu=R_\mu^{\phantom{\mu}\nu}\,\partial_\nu\eqpc
			\end{equation*}
			where $R$ is a rotation matrix (\emph{i.e.} $R\in SO\lr{(D)}$) such that the direction $\lr{(1,1,\cdots,1)}$ is mapped to $\lr{(1,0,\cdots,0)}$ (in that way, each cell will only contain one unknown).
			So,
			\begin{equation*}
				R_0^{\phantom{0}\nu}=\frac{1}{\sqrt{D}}\eqpd
			\end{equation*}
			Then, it remains to orthogonalize the remaining rows of $R$, which can be achieved by defining
			\begin{align*}
				R_\mu^{\phantom{\mu}\nu}&=\frac{1}{\sqrt{\mu\lr{(\mu+1)}}}&\nu<\mu\eqpc\\
				R_\mu^{\phantom{\mu}\nu}&=-\sqrt{\frac{\mu}{\mu+1}}&\nu=\mu\eqpc\\
				R_\mu^{\phantom{\mu}\nu}&=0&\nu>\mu\eqpc
			\end{align*}
			where $\mu\in\lr{\llbracket1,d\rrbracket}$ and $\nu\in\lr{\llbracket0,d\rrbracket}$.
			
			The space-time is now discretized using an hypercubic lattice rotated by $R$, and the discretization rules are simple concatenations of midpoint rules as introduced in~\cite{BridgesReich2001}.
			
			Such a method would have the same property of multi-symplecticity and locality while, hopefully, remaining stable in any dimension.
			But we leave the complete development of it for future work.
	\sectionstar{Conclusions}
		The purpose of this paper was to introduce a new numerical method to integrate partial differential equations stemming from the Hamiltonian dynamics of field theories.
		The method is a centered box scheme, implemented on the light-cone coordinates, in such a way to restore the locality of the algorithm without losing its multi-symplectic properties.
		
		Our method has \emph{local} conservation properties (and therefore global conservation properties as well) in agreement with what is generally achieved by multi-symplectic integrators.
		The errors committed do not strongly accumulate, remaining very small over very long periods of time.
		This is important in applications in which the long-time limit of evolution should be reached with good confidence.
		
		In the process of comparing the performance of our algorithm to other ones in the literature we showed that exact global conservation properties, as the ones imposed in the \textsc{BDdV} technique, do not necessarily guarantee small errors in the local conservation laws.
		
		We highlighted the link between the \textsc{De Donder}~-- \textsc{Weyl} formalism of field theories and the multi-symplectic structure of phase space, and we treated the latter on a rigorous geometric way.
		We developed the construction of the stress-energy tensor in the Hamiltonian formalism.
		We showed that it is exactly conserved in the continuum and we derived the error committed by the algorithm in its discrete implementation.
		In particular, we showed that it is exactly preserved for a linear equation.
		
		Interestingly, depending on the model that we considered, the multi-symplectic structure was found to be degenerate in spatial dimension larger than zero.
		We showed how to solve this problem in any dimension using the particular case of the wave equation as an example.
		The generalization to other field equations should follow similar steps.
	\sectionstar{Acknowledgements}
		We thank G. \textsc{Biroli} and A. \textsc{Ciocchetta} for useful discussions.
		L.~F.~C. is a member of Institut Universitaire de France.
	\bibliography{bibliography}
\end{document}